\documentclass[smallextended,envcountsect]{svjour3}

\smartqed  
\usepackage{amsmath}
\usepackage{amssymb}
\usepackage{booktabs}
\usepackage[mathscr]{eucal}
\usepackage{graphicx}
\usepackage{color}
\usepackage{empheq}
\usepackage{graphics}
\usepackage{epstopdf}
\usepackage{array}
\usepackage{multirow}
\usepackage{makecell}
\usepackage[justification=centering]{caption}
\usepackage[latin1]{inputenc}
\usepackage[T1]{fontenc}
\usepackage{graphicx}
\usepackage{subcaption}
\usepackage{hyperref}
\usepackage{xcolor}
\usepackage{enumitem}
\usepackage{algpseudocode}
\hypersetup{
    colorlinks=true,    
    linkcolor=black,    
    citecolor=black,    
    urlcolor=black      
}
\numberwithin{equation}{section}

\newtheorem{algorithm}{\sc\bf Algorithm}[section]
\begin{document}
\title{ Single-loop variance reduction methods in Bregman setups for finite-sum structured variational inequalities}

\author{ Zhong-bao Wang$^{1,2}$, Zhong-cheng Zhang$^3$ }

\institute{
$^1$ Department of Mathematics, Southwest Jiaotong University, Chengdu, Sichuan 611756,  China\\
$^2$ National Engineering Laboratory of Integrated Transportation Big Data Application Technology, Chengdu, Sichuan 611756, China\\
$^3$ School of Electrical Engineering, Southwest Jiaotong
 Univeristy, Chengdu, Sichuan 611756,  China\\
$^*$Corresponding author: zhongbaowang@hotmail.com}

\date{Received: date / Accepted: date}

\maketitle
\noindent{\bf Abstract}
\par  In this paper, we address variational inequalities (VI) with a finite sum structure by proposing a novel single-loop variance-reduced algorithm that incorporates the Bregman distance. Under the monotone setting, we establish the almost sure convergence of the proposed algorithm and prove that it achieves the optimal complexity of $\mathcal{O}\left(\frac{\sqrt{M}}{\varepsilon }\right)$ for finding an $\varepsilon$-gap. Furthermore, under the non-monotone setting, we derive a complexity of $\mathcal{O}\left(\frac{1}{\varepsilon^2 }\right)$ of the algorithm. Our proposed method yields complexity results that either match or improve the state-of-the-art complexity bounds reported in existing literature. Notably, this work is the first to rigorously establish the linear convergence rate of the algorithm for solving finite-sum variational inequalities in Bregman setups. Finally, we report two numerical experiments to validate the effectiveness and practical performance of our method. \vspace{12pt}


\section{Introduction}
\par Let $\mathcal{X}$ be a finite dimensional spaces equipped with the inner product $\left\langle { \cdot , \cdot } \right\rangle $ and the norm $\left\|  \cdot  \right\|$, whose dual norm is represented by $\left\|  \cdot \right\|_*$. Throughout the paper, we consider the following variational inequalities (VI): find a point $x^*\in K$ such that
\begin{equation}
\label{1.1}
\langle F(x^*), y-x^*\rangle +g(y)-g(x^*)\geq 0,~\forall~ y\in K,
\end{equation}
where $K \subset \mathcal{X}$ is a nonempty closed convex set, $g:K \to (-\infty,+\infty]$ is a proper convex lower semicontinuous function and $F: {\rm{dom}}g=\{x\in\mathcal{X}:g(x)<+\infty\} \to {{\mathcal{X}}}$ is a operator. Optimization problems and saddle-point problems can be formulated as VI\eqref{1.1} and its special cases, see \cite{FP}.
\par Variational inequalities was initially proposed by Fichera \cite{FG,FG+1} in 1963 for addressing the Signorini problem. Subsequently, Stampacchia \cite{StaG} studied variational inequalities to resolve partial differential equations involving unilateral boundary conditions and elliptic-type free boundary value problems. Later, variational inequalities are applied to classical problems such as equilibrium theory, economics, game theory, control theory and differential equations, see \cite{FP,K+1,KS,LXJ,LXJ1,LXJ2,PT3}. The advancement of variational inequalities is currently dominated by the machine learning community, wherein fundamental problems                                                                                                                                                              like reinforcement learning, adversarial training, generative models can all be reformulated as variational inequalities through first-order optimality conditions to facilitate in-depth theoretical analysis, see \cite{Gid,Sid}.

 \par  The simplest approach to solve variational inequalities is the projected gradient method. However, this method converges only when the cost operator $F$ is strongly monotone. Korpelevich \cite{KGM} established the extragradient method (EGM) and prove its convergence under the conditions that is $F$ monotone and Lipschitz continuous. Unfortunately, EGM necessitates computing the proximal operators of the mapping $g$ twice per iteration. When the proximal operators of the mapping $g$ is non-explicit or complex, this significantly increases the computational burden. Additionally, EGM requires two evaluations of the cost operator, which can compromise the method's efficiency, particularly in contexts where evaluating the cost operator $F$ is computationally expensive. Reducing the number of operator calls and proximal operators (or projections) used per iteration become one of the best important direction of developing algorithms for solving variational inequalities, for example  forward-backward-forward (FBF) \cite{TPA}, Golden ratio algorithm(GRA) \cite{LXJ,YM+1}, reflected gradient method \cite{MY1} and forward-reflected-backward (FoRB) \cite{MY1,MY2} among others.
 In 1967, Bregman \cite{Bregman(1967)} introduced a dissimilarity measure that describes the difference between two points on a differentiable convex function. This measure, known as the Bregman distance or Bregman divergence, has been shown to enhance algorithmic flexibility across diverse contexts, yielding notable computational efficiencies in specific scenarios, see \cite{IC,NA,ZBW}. To solve variational inequalities, Izuchukwu et al. \cite{IC} introduced one-step Bregman projection methods that require only one proximal operator and one evaluation of the cost operator per iteration.
It was proved that the methods converge weakly when the cost operator is pseudomonotone, and strongly when it is strongly pseudomonotone.
Additionally, they pointed out that "Part of our future projects is to study the rate of convergence of the proposed methods of this paper."
Establishing the convergence rate for the one-step Bregman projection method constitutes the first motivation of this work.

\par A common structure of large-scale setting arises in VI\eqref{1.1} is that the operator $F$ can be written as a finite sum $F = \frac{1}{M}\sum\limits_{i = 1}^M {{F_i}} $, for concrete examples,  see \cite{AA,AA1,AZ,PA}.
For notational convenience, when the operator $F$ adopts the aforementioned finite sum structure, we denote this particular class of variational inequalities simply as FSVI \eqref{1.1}. In addition,  the formulation $F = \frac{1}{M}\sum\limits_{i = 1}^M {{F_i}} $  can be regarded as an approximation of the expected value of a random function. Therefore, FSVI \eqref{1.1} also is called to be stochastic finite-sum variational
inequalities, see \cite{PA1,PA}. If $M$ is large,  then the computational cost of the operator
$F$ may escalate to a prohibitive level.  Variance reduction techniques leverage this specific structure to mitigate the computational complexity of deterministic methods in optimization problems, achieving this by furnishing unbiased gradient estimators while simultaneously minimizing the variance of estimation errors, see \cite{Gow,Sch}. Recently, there has been a growing interest in developing variance-reduced techniques for saddle
point problems and FSVI\eqref{1.1}. In \cite{AA2}, the Forward-reflected-backward method with variance
reduction (VR-FoRB) was introduced for solving FSVI\eqref{1.1}, and later extended to the general Bregman case. In the monotone case, the method was the first to be shown to achieve almost sure convergence. When $F$ is strongly monotone, the method exhibits linear convergence without requiring knowledge of the strong monotonicity constant. However,   it only gets a complexity of  $\mathcal{O}\left(\frac{M}{\varepsilon }\right)$ and thus lacks the complexity improvements offered by deterministic methods. Consequently, in the conclusion of \cite{AA2}, the authors posed an open problem: to develop a method that (i) converges under the same minimal set of assumptions as their algorithm, and (ii) provides improved complexity guarantees compared to deterministic methods. In Euclidean setup, Alacaoglu and Malitsky \cite{AA} answered this open problem by proposing EGM with variance reduction and established the convergence, the improved complexity of $\mathcal{O}\left(\frac{\sqrt{M}}{\varepsilon }\right)$ to find an $\varepsilon$-gap and the linear convergence for the method. Using the same technique,  variance-reduced versions of different algorithms including forward-backward-forward and forward-reflected-backward methods for monotone inclusion are studied. In Bregman setting, they also introduced a variance-reduced extragradient method with a double-loop structure, that is Algorithm 2 of \cite{AA}, and achieved the same complexity of $\mathcal{O}\left(\frac{\sqrt{M}}{\varepsilon }\right)$. However, they did not establish convergence guarantees (including linear convergence) in this setting. The second motivation of this work is to introduce a new variance reduction method in Bregman setups,  prove its almost sure convergence and obtain the same complexity of $\mathcal{O}\left(\frac{\sqrt{M}}{\varepsilon }\right)$, thereby under Bregman setups addressing the open problem proposed in \cite{AA2}.

\par By integrating the stochastic path-integrated differential estimator (SPIDER) technique proposed in \cite{Fang}, Alizadeh et al. \cite{AZ} for the first time developed a single-loop method with Bregman distance functions, referred to as the stochastic variance-reduced forward reflected moving average backward method (VR-FoRMAB), for solving both monotone and non-monotone FSVI\eqref{1.1}. VR-FoRMAB achieves an $\varepsilon$-gap with complexities of $\mathcal{O}\left(\frac{\sqrt{M}}{\varepsilon }\right)$ and $\mathcal{O}\left(\frac{1}{\varepsilon^2 }\right)$ for the monotone and nonmonotone scenarios, respectively.  In Bregman frameworks, research on developing single-loop  variance reduction methods for solving  variational inequalities and saddle point problems is still in its infancy. This is presumably because, as noted by the authors of \cite{AZ}, combing the Bregman distance with variance reduction techniques introduces an extra layer of complexity to algorithmic analysis. Unfortunately, Alizadeh et al.\cite{AZ} did not provide any proof for the almost sure convergence and linear convergence of VR-FoRMAB, regardless of whether the underlying mappings are monotone or non-monotone.
 The third motivation for this work is to develop a novel single-loop variance reduction algorithm within the context of the Bregman distance for solving FSVI\eqref{1.1}, alongside proving the  almost sure convergence and linear convergence of the algorithm.

\par Next, we elaborate on other related works pertinent to our research. Several notable studies explore variance reduction techniques for addressing monotone FSVI\eqref{1.1} in \cite{PA,PA1,SC}.
In the Euclidean setting, Pichugin et al. \cite{PA1} proposed an algorithm named Optimistic Method with Momentum and Batching (OM-MB*).  This method supports batching and exhibits a complexity of $\mathcal{O}\left(\frac{\sqrt{M}}{\varepsilon }\right)$. Under Bregman setups, Pichugin et al. \cite{PA} introduced another optimistic method with momentum and batching (OM-MB) which features a double-loop structure and also achieves a complexity of $\mathcal{O}\left(\frac{\sqrt{M}}{\varepsilon }\right)$. Song and Diakonikolas \cite{SC} proposed a cyclic coordinate dual averaging method with extrapolation (CODER). Furthermore, they integrated variance reduction techniques into the CODER algorithm, thereby developing a three-loop algorithm termed VR-CODER, that achieves a computational complexity of $\mathcal{O}\left(\frac{\sqrt{M}}{\varepsilon }\right)$. Exploring variational inequalities with non-monotone classes including
pseudo monotonicity, quasi monotonicity, cohypomonotonicity, and
weak Minty variational inequality (WMVI) are studied in \cite{DCD,PT1,PT2,Wang,Wang1}. The weak Minty variational inequality (WMVI) is an important nonmonotone class, see \cite {DCD,PT1,PT2}.  For solving nonlinear inclusion problems with WMVI condition,  Pethick et al. \cite{PT1} proposed a bias-corrected stochastic extragradient(BC-SEG+) algorithm with a complexity of $\mathcal{O}\left(\frac{1}{\varepsilon^4 }\right)$; Alacaoglu et al. \cite{AA1} introduced a multi-loop inexact variant of Krasnoselski-Mann (KM) iterations and derived a complexity of $\mathcal{O}\left(\frac{1}{\varepsilon^4 }\right)$  for stochastic nonlinear inclusion problems with WMVI condition.
For a further comparison of the existing methods, please see Table \ref{table1} and Table \ref{table2}.
\begin{table}[h!]
\caption{Compression of some methods for solving monotone FSVI\eqref{1.1}}
\label{table1}
\centering
\begin{tabular}{|c|c|c|c|}
\hline
Ref & Bregman & Loop & Complexity \\ \hline
VR-MP \cite{AA}       & Yes       & 2 &  $\mathcal{O}\left(\frac{\sqrt{M}}{\varepsilon }\right)$  \\ \hline
VR-FoRB \cite{AA2}       & No       & 1 &  $\mathcal{O}\left(\frac{M}{\varepsilon }\right)$  \\ \hline
OM-MB \cite{PA}     & Yes       & 2    &$\mathcal{O}\left(\frac{\sqrt{M}}{\varepsilon }\right)$   \\ \hline
OM-MB$^*$ \cite{PA1}       & No       & 1    & $\mathcal{O}\left(\frac{\sqrt{M}}{\varepsilon }\right)$  \\ \hline
VR-CODER \cite{SC}       & Yes       & 3     &$\mathcal{O}\left(\frac{\sqrt{M}}{\varepsilon }\right)$  \\ \hline
This paper      & Yes       & 1     &$\mathcal{O}\left(\frac{\sqrt{M}}{\varepsilon }\right)$  \\ \hline
\end{tabular}
\end{table}
\begin{table}[h]
\caption{Compression of some methods for solving non-monotone problem}
\label{table2}
\centering
\begin{tabular}{|c|c|c|c|c|}
\hline
Ref &Assumption& Bregman & Loop & Complexity \\ \hline
Inexact-KM \cite{AA1}  & WMVI& No& multi-loop &  $\mathcal{O}\left(\frac{1}{\varepsilon^4 }\right)$  \\ \hline
VR-FoRMAB \cite{AZ}    & WMVI&Yes & 1 &  $\mathcal{O}\left(\frac{1}{\varepsilon^2 }\right)$  \\ \hline
BC-SEG+ \cite{PT1}&WMVI & No & 1    &$\mathcal{O}\left(\frac{1}{\varepsilon^4 }\right)$   \\ \hline
This paper    &WMVI  & Yes       & 1     &$\mathcal{O}\left(\frac{1}{\varepsilon^2 }\right)$  \\ \hline
\end{tabular}

\end{table}

 This paper mainly focuses on proposing a novel single-loop variance-reduced algorithm that embeds Bregman distance function for solving FSVI\eqref{1.1}. Our proposed method requires only a single computation of the Bregman proximal mapping onto the feasible set per iteration, while incorporating both inertial acceleration and a batch size selection. We emphasize that our algorithm is fundamentally distinct from VR-FoRMAB, as it does not employ the SPIDER technique. In the monotone setting, we establish the almost sure convergence of this algorithm and derive its complexity of $\mathcal{O}\left(\frac{\sqrt{M}}{\varepsilon }\right)$. Thus, under Bregman frameworks, we resolve the open problem posed in \cite{AA2}. In the non-monotone settings where FSVI\eqref{1.1} admits weak Minty solutions, we establish a complexity of $\mathcal{O}\left(\frac{1}{\varepsilon^2 }\right)$ for our algorithm. It is particularly noteworthy that, to the best of our knowledge, we first rigorously show the linear convergence of algorithms for solving FSVI\eqref{1.1} under the Bregman framework. Some numerical experiments show that the proposed algorithms outperform some existing ones.
\par The structure of the paper is as follows: Section 2 introduces the necessary definitions and lemmas for this article. In Section 3, we present the main results of the proposed algorithm. Section 4 includes several numerical experiments to showcase the performance of our algorithm. Finally, Section 5 provides concluding remarks.
\section{Preliminaries}
 In this section, we recapitulate some notations, definitions, and results that will be utilized in subsequent parts. The notations $\mathbb{E}[\xi]$ and $\mathbb{E}[\xi|{\mathcal{F}}]$  denote the expectation of a random variable $\xi$  and the conditional expectation of $\xi$  with respect to a $\sigma$-algebra $\mathcal{F}$, respectively. Meanwhile, $\sigma({\xi}_1,...,{\xi}_k)$ represents the $\sigma$ -algebra generated by the random variables $\left\{ {{\xi _i}} \right\}_{i = 1}^k$.  For a set $C \in \mathcal{X}$, the indicator function of $C$ is defined as
\[{{\rm{I}}_C}\left( x \right) = \left\{ \begin{array}{l}
0,x \in C,\\
 + \infty ,otherwise.
\end{array} \right.\]
Given a function $f:K \to \left( { - \infty , + \infty } \right]$, the domain of $f$ denotes  ${\rm{dom}} f: = \left\{ {x \in K:f\left( x \right) < \infty } \right\}$ and the subdifferential of $f$ at $x \in {\rm{dom}} f$ is defined by
\[\partial f\left( x \right) = \left\{ {v \in {\mathcal{X}}:f\left( x \right) - f\left( y \right) - \left\langle {v,x - y} \right\rangle  \le 0,\forall y \in {\mathcal{X}}} \right\},\]
and the effective domain of $\partial f$ is denoted as ${\rm{dom}}\left(\partial f\right) = \left\{ {x \in {\rm{dom}}f:\partial f\left( x \right) \ne \emptyset } \right\}$.
\par Let $f:K \to ( - \infty , + \infty ]$ be a proper convex lower semi-continuous function such that $\mathop{\rm{dom}} g \subset \mathop{\rm{dom}} f$, $f$ is differentiable on ${\rm{dom}}\left(\partial f\right)$ and $f$ is 1-strongly convex on $\mathop{\rm{dom}} f$, that is,
 \[f\left( x \right) \ge f\left( y \right) + \left\langle {\nabla f\left( y \right),x - y} \right\rangle  + \frac{1 }{2}{\left\| {x - y} \right\|^2},~\forall x \in \mathop{\rm{dom}} f,y \in {\rm{dom}}\left(\partial f\right).\]
 The Bregman distance $D:\mathop{\rm{dom}} g~\times~{\rm{dom}}\left(\partial f\right) \to [0,+\infty)$ associated with $f$ is defined as
 \[D\left( {x,y} \right): = f\left( x \right) - f\left( y \right) - \left\langle {\nabla f\left( y \right),x - y} \right\rangle.\]

 Due to  1-strongly convex property  of  $f$, we have
 \begin{equation}\label{2.1}
 D\left( {x,y} \right) \ge \frac{1}{2}{\left\| {x - y} \right\|^2},~\forall x \in \mathop{\rm{dom}} g,~y \in {\rm{dom}}\left(\partial f\right).
 \end{equation}

 From the definition of Bregman distance, it is easy to obtain that
  \begin{equation}\label{l2.1}
 D\left( {x,y} \right) = D\left( {x,z} \right) + D\left( {z,y} \right) + \left\langle {\nabla f\left( y \right) - \nabla f\left( z \right),z - x} \right\rangle,~\forall y,z \in {\rm{dom}}\left(\partial f\right),~x \in \mathop{\rm{dom}} f.
  \end{equation}

 \par The generalized Bregman proximal mapping is defined by
 \[{P_x}\left( {y,\alpha } \right): = \mathop {\arg \min }\limits_{z \in K} \left\{ {\alpha g\left( z \right) + \left\langle {\alpha y,x - z} \right\rangle  + D\left( {z,x} \right)} \right\},~x \in {\rm{dom}}f,~\alpha  > 0,~y \in \mathcal{X}.\]
 Next, we give the following lemmas.
 \begin{lemma} \cite{AA}
 \label{l2}
 Let ${z^ + } = \mathop {\arg \min }\limits_{z \in K} \left\{ {\alpha g\left( z \right) + \alpha \left\langle {y,z} \right\rangle  + \gamma D\left( {z,{z_1}} \right) + \left( {1 - \gamma } \right)D\left( {z,{z_2}} \right)} \right\}$, then, for $\alpha>0$, $0 \le \gamma  \le 1$, $z_1,~z_2 \in {\rm{dom}}\left(\partial f\right)$ and $y \in \mathcal{X}$, we have ${x^+} \in {\rm{dom}}\left(\partial f\right)$ and
\begin{equation}\label{2.2}
\begin{aligned}
\alpha g\left( z \right) - \alpha g\left( {{z^ + }} \right) + \alpha \left\langle {y,z-{z^ + }} \right\rangle  &\ge D\left( {z,{z^ + }} \right) + \gamma \left( {D\left( {{z^ + },{z_1}} \right) - D\left( {z,{z_1}} \right)} \right)\\
 &+ \left( {1 - \gamma } \right)\left( {D\left( {{z^ + },{z_2}} \right) - D\left( {z,{z_2}} \right)} \right),\forall z \in K.
\end{aligned}
\end{equation}
 \end{lemma}
 The proof of Lemma \ref{l2} see Lemma 2.4 from \cite{AA}.
\begin{lemma}\cite{RH}
\label{l2'}
Let$\{ y_s \},~ \{ u_s \},~\{ a_s \}$ and $\{ b_s \}$ be nonnegative random variables adapted to the filtration $\{ F_s \}$ such that almost surely $\sum\limits_s {{a_s}} < \infty$, $\sum\limits_s {{b_s}} < \infty$, and for all $s$,
\[E\left[ {\left. {{y_{s + 1}}} \right|{F_s}} \right] \le \left( {1 + {a_s}} \right){y_s} - {u_s} + {b_s}.\]
Then, almost surely $\{ y_s \}$ converges and $\sum\limits_s {{u_s}}  < \infty $.
\end{lemma}
Finally, we introduce the following definitions will be used.
\begin{definition}
\cite{AA}
The operator $F$  has a stochastic oracle $F_\xi$ that is variable ${\bar L}$-Lipschitz in mean: for any $x,y \in {\rm{dom}} g$ there exists a distribution $Q_{x,y}$ such that
\begin{itemize}
  \item [$(i)$] $F_\xi$ is unbiased: $F\left( z \right) = {E_{\xi \sim {Q_{x,y}}}}\left[ {{F_\xi}\left( z \right)} \right],\forall z \in {\rm{dom}}g;$
  \item [$(ii)$] ${E_{\xi \sim {Q_{x,y}}}}\left[ {{{\left\| {{F_\xi}\left( x \right) - {F_\xi}\left( y \right)} \right\|}_*^2}} \right] \le {{\bar L}^2}{\left\| {x - y} \right\|^2}.$
\end{itemize}
\end{definition}
\begin{remark}
It is crucial to emphasize that the second condition applies exclusively for given $x$, $y$. In contrast, the constant
${\bar L}$ remains invariant for all $x$, $y$.
 The alteration of $x$ and $y$ inherently results in a transformation of the distribution, which justifies the term "variable."
\end{remark}

\begin{definition}
\cite{Wang1}
The operator $F$ is said to be $g$-pseudomonotone if
\[\left\langle {F\left( x \right),y - x} \right\rangle  + g\left( y \right) - g\left( x \right) \ge 0 \Rightarrow \left\langle {F\left( y \right),y - x} \right\rangle  + g\left( y \right) - g\left( x \right) \ge 0, \forall ~x,y \in K.\]
\end{definition}
\begin{definition}
\cite{OO}
The operator $F$ is said to be $\beta$-strongly pseudomonotone if there exists some constant $\beta>0$ satisfying
\[\left\langle {F(x),y - x} \right\rangle  \ge 0 \Rightarrow \left\langle {F(y),y - x} \right\rangle  \ge \beta {\left\| {x - y} \right\|^2},\forall x,y \in K.\]
\end{definition}

\begin{definition}
\label{def3.4} \cite{PT1}
For a mapping $T:{\mathcal{X}} \to {2^{{\mathcal{X}}}}$, if there exists $\rho >0$ such that \[\left\langle {u - v,x - y} \right\rangle  \ge  - \rho {\left\| {u - v} \right\|^2},\forall (x,u),(y,v) \in {\rm{gra}}(T),\]
then $T$ is called $\rho$-cohypomonotone, where ${\rm{gra}}(T) = \left\{ {\left( {x,y} \right) \in {\mathcal{X}} \times {\mathcal{X}}:y \in Tx} \right\}$, and ${2^{{\mathcal{X}}}}$ is the set of all subsets of $\mathcal{X}$.
\end{definition}
\section{Algorithm and convergence analysis}
The notation \[\mathbb{E}\left[ { \cdot \left| {\sigma \left( {{x_0},{x_1},\dots,{x_{s}},{w_0},{w_1},\dots,{w_{s-1}}} \right)} \right.} \right] = {\mathbb{E}_s}\left[  \cdot  \right]\] is used to denote the conditional expectation.
Now we introduce our algorithm.
\renewcommand{\thealgorithm}{1}
\begin{algorithm}
\label{Alg1}
\mbox{}\\
\begin{algorithmic}[1]
\State \textbf{Parameters:} Choose $0< \gamma  \le 1$, $0< p \le 1$ stepsize $\alpha>0$, batch size $b \in \{1,...,M\}.$
\State \textbf{Initialization:} Take $x_0=x_{-1}=w_0=w_{-1} \in K$.
\For{$s=0,1,2,...,S-1$}
\par \State Fix distribution $Q_{x_s,w_{s-1}}$ and sample set $\{1,...,M\}$ according to it;
\par \State Draw samples $\xi_1^s$, $\cdots$, $\xi_b^s$ independently and uniformly from $\{1,...,M\}$;
\par \State $B^s=\{\xi_1^s,...,\xi_b^s\};$
\par \State ${\bar x_s} = \nabla {f^{ - 1}}\left( {\gamma \nabla f\left( {{x_s}} \right) + \left( {1 - \gamma } \right)\nabla f\left( {{x_{s - 1}}} \right)} \right)$;\label{1}
\par \State ${\Delta ^s} = \frac{1}{b}\sum\limits_{\xi \in {B^s}} {F\left( {{w_s}} \right) + } {F_\xi}\left( {{x_s}} \right) - {F_\xi}\left( {{w_{s - 1}}} \right);$\label{7}
\par \State ${x_{s + 1}} = {P_{\bar x_s}}\left( {\Delta^s,\alpha } \right)= \mathop {\arg \min }\limits_{x \in K} \left\{ {\alpha g\left( x \right) + \left\langle {\alpha {\Delta ^s},x } \right\rangle  +  D\left( {x,{\bar x_s}} \right) } \right\},$\label{8}
\par \State ${w_{s + 1}} = \left\{ \begin{array}{l}
{x_{s + 1}},~ with ~probability ~p,\\
{x_{s}},~~~ with ~probability ~1 - p.
\end{array} \right.$
\EndFor
\end{algorithmic}
\end{algorithm}

In line \ref{1}, the item ${\bar x_s} = \nabla {f^{ - 1}}\left( {\gamma \nabla f\left( {{x_s}} \right) + \left( {1 - \gamma } \right)\nabla f\left( {{x_{s - 1}}} \right)} \right)$ is referred to as the inertial item, which can enhance the convergence rate of algorithms for variational inequalities and optimization problems, see \cite{ZBW}. Convex combination updates analogous to the inertial item are also widely employed, for instance, the golden ratio technique in \cite{LXJ,TMK} and the negative momentum method in \cite{PA}. Motivated by \cite{PA1} and \cite{PA}, our algorithm supports batching. Specifically, choosing an appropriate batch size can enhance both the convergence rate and convergence complexity of the algorithm, as illustrated in the proof of Corollary 3.1 and Table 3.
Compared to VR-FoRB method \cite{AA2}
\[\begin{array}{l}
{x_{s + 1}} = \mathop {\arg \min }\limits_z \left\{ {g\left( z \right) + \left\langle {F\left( {{w_s}} \right) + {F_{{i_s}}}\left( {{x_s}} \right) - {F_{{i_s}}}\left( {{w_{s - 1}}} \right),z - {x_s}} \right\rangle  + \frac{1}{\tau }D\left( {z,{x_s}} \right)} \right\}\\
{w_{s + 1}} = \left\{ \begin{array}{l}
{x_{s + 1}},~ with ~probability ~p,\\
{w_{s}},~~~ with ~probability ~1 - p,
\end{array} \right.
\end{array}\]
where $\tau$ is step size, $i_s$ is an index  independently and uniformly from $\{1,...,M\}$.
\par In the update of $w_{s+1}$, we use the most recent iterate $x_s$  instead of $w_s$ in VR-FoRB method.
This distinction and selection of batch sizes serve as crucial factors in ensuring that Algorithm \ref{Alg1} achieves a complexity of $\mathcal{O}\left(\frac{\sqrt{M}}{\varepsilon }\right)$ for finding an $\varepsilon$-gap. In addition, when $\gamma=p=M=1$ and $g \equiv 0$, then Algorithm \ref{Alg1} reduces to the one-step Bregman projection method that is Algorithm 3.2 of \cite{IC}.

To further analyze the problem, we introduce the following assumptions:\\
  $\mathbf{Assumption}$ (A1) The set $S_1 \neq\emptyset$, where $S_1$ is the solution set of FSVI (\ref{1.1});\\
  $\mathbf{Assumption}$ (A2) The operator $F$ is monotone, that is,  $\left\langle {F\left( x \right) - F\left( y \right),x - y} \right\rangle  \ge 0,~\forall~x,y\in K$;\\
  $\mathbf{Assumption}$ (A3) The operator $F$ has a stochastic oracle $F_\xi$ that is variable ${\bar L}$-Lipschitz in mean;\\
  $\mathbf{Assumption}$ (A4) The operator $F$ is $L$-Lipschitz continuous, that is, ${\left\| {F\left( x \right) - F\left( y \right)} \right\|_*} \le L\left\| {x - y} \right\|,~\forall x,y \in K$.\\
Firstly, we estimate the variance of ${\Delta}^s$.
\begin{lemma}
\label{l3.1}
Assume that Assumption (A3) holds. If $\left\{ {{x_{s}}} \right\}$  is a sequence generated by Algorithm \ref{Alg1}, then the following inequality holds:
\[{\mathbb{E}}\left[ {{{\left\| {{\Delta ^s} - {E_s}\left[ {{\Delta ^s}} \right]} \right\|}_*^2}} \right] \le \frac{{{\bar L}^2}}{b}{\mathbb{E}}\left[ {{{\left\| {{x_s} - {w_{s - 1}}} \right\|}^2}} \right].\]
\end{lemma}
\begin{proof}
Since stochastic oracle $F_{\xi}$ are unbiased, we can get
\begin{equation}\label{l3.1.e1}
{{\mathbb{E}}_s}\left[ {{\Delta ^s}} \right] = F\left( {{w_s}} \right)+ F\left( {{x_s}} \right) - F\left( {{w_{s - 1}}} \right).
\end{equation}
By the definition of ${\Delta}^s$, we can get
\begin{equation}
\label{l3.1.e2}
\begin{aligned}
{{\mathbb{E}}_s}\left[ {{{\left\| {{\Delta ^s} - {{\mathbb{E}}_s}\left[ {{\Delta ^s}} \right]} \right\|}_*^2}} \right] &= {{\mathbb{E}}_s}\left[ {{{\left\| {\frac{1}{b}\sum\limits_{\xi \in {B^s}} {F\left( {{w_s}} \right) + } {F_\xi}\left( {{x_s}} \right) - {F_\xi}\left( {{w_{s - 1}}} \right) - \left( {F\left( {{w_s}} \right)+ F\left( {{x_s}} \right) - F\left( {{w_{s - 1}}} \right)} \right)} \right\|}_*^2}} \right]\\
 &= {{\mathbb{E}}_s}\left[ {{{\left\| {\frac{1}{b}\sum\limits_{\xi \in {B^s}} {{F_\xi}\left( {{x_s}} \right) - {F_\xi}\left( {{w_{s - 1}}} \right) - \left( {F\left( {{x_s}} \right) - F\left( {{w_{s - 1}}} \right)} \right)} } \right\|}_*^2}} \right].
\end{aligned}
\end{equation}
Since $\xi_1^s,\cdots,\xi_b^s$ in $B^s$ are uniform and independent, for any $1\leq i < l \le b$, we know that
\begin{equation}\label{el13.e1}
{\mathbb{E}_s}\left[ {\left\langle {{F_{\xi_i^s}}\left( {{x_s}} \right) - {F_{\xi_i^s}}\left( {{w_{s - 1}}} \right) - \left( {F\left( {{x_s}} \right) - F\left( {{w_{s - 1}}} \right)} \right),{F_{\xi_l^s}}\left( {{x_s}} \right) - {F_{\xi_l^s}}\left( {{w_{s - 1}}} \right) - \left( {F\left( {{x_s}} \right) - F\left( {{w_{s - 1}}} \right)} \right)} \right\rangle } \right] = 0.
\end{equation}
Hence, we have
\begin{equation}\label{el13.e2}
{{\mathbb{E}}_s}\left[ {{{\left\| {{\Delta ^s} - {{\mathbb{E}}_s}\left[ {{\Delta ^s}} \right]} \right\|}_*^2}} \right]
= \frac{1}{b^2}{{\mathbb{E}}_s}\left[ {\sum\limits_{\xi \in {B^s}} {{{\left\| {{F_\xi}\left( {{x_s}} \right) - {F_\xi}\left( {{w_{s - 1}}} \right) - \left( {F\left( {{x_s}} \right) - F\left( {{w_{s - 1}}} \right)} \right)} \right\|}_*^2}} } \right].
\end{equation}
Using the fact that ${\mathbb{E}}\left[ {{{\left\| {X - {\mathbb{E}}X} \right\|}^2}} \right] = {\mathbb{E}}\left[ {{{\left\| X \right\|}^2}} \right] - {\left\| {{\mathbb{E}}X} \right\|^2}$, we have
\begin{equation}\label{l3.1.e3}
\begin{array}{l}
{{\mathbb{E}}_s}\left[ {\sum\limits_{\xi \in {B^s}} {{{\left\| {{F_\xi}\left( {{x_s}} \right) - {F_\xi}\left( {{w_{s - 1}}} \right) - \left( {F\left( {{x_s}} \right) - F\left( {{w_{s - 1}}} \right)} \right)} \right\|}_*^2}} } \right]\\
 \le {{\mathbb{E}}_s}\left[ {\sum\limits_{\xi \in {B^s}} {{{\left\| {{F_\xi}\left( {{x_s}} \right) - {F_\xi}\left( {{w_{s - 1}}} \right)} \right\|}_*^2}} } \right].
\end{array}
\end{equation}
Since $\xi_1^s,...,\xi_b^s$ in $B^s$ are uniform and independent, and Assumption (A3) holds, we can get
\begin{equation}\label{l3.1.e4}
\begin{array}{l}
\frac{1}{b^2}{{\mathbb{E}}_s}\left[ {\sum\limits_{\xi \in {B^s}} {{{\left\| {{F_\xi}\left( {{x_s}} \right) - {F_\xi}\left( {{w_{s - 1}}} \right)} \right\|}_*^2}} } \right] \\
= \frac{1}{b}{{\mathbb{E}}_s}\left[ {{{\mathbb{E}}_{\xi{\sim} Q_{x_s,w_{s-1}}}}\left[ {{{\left\| {{F_\xi}\left( {{x_s}} \right) - {F_\xi}\left( {{w_{s - 1}}} \right)} \right\|}_*^2}} \right]} \right]\\
 \le \frac{{{\bar L}^2}}{b}{{\mathbb{E}}_s}\left[ {{{\left\| {{x_s} - {w_{s - 1}}} \right\|}^2}} \right].
\end{array}
\end{equation}
In summary of the above inequalities, we have
\begin{equation}\label{l3.1.e5}
{{\mathbb{E}}_s}\left[ {{{\left\| {{\Delta ^s} - {{\mathbb{E}}_s}\left[ {{\Delta ^s}} \right]} \right\|}_*^2}} \right] \le \frac{{{\bar L}^2}}{b}{{\mathbb{E}}_s}\left[ {{{\left\| {{x_s} - {w_{s - 1}}} \right\|}^2}} \right].
\end{equation}
Taking the full expectation of both parts, we can get
\[{\mathbb{E}}\left[ {{{\left\| {{\Delta ^s} - {E_s}\left[ {{\Delta ^s}} \right]} \right\|}_*^2}} \right] \le \frac{{{\bar L}^2}}{b}{\mathbb{E}}\left[ {{{\left\| {{x_s} - {w_{s - 1}}} \right\|}^2}} \right].\]
This proof is finished.

\end{proof}

Foremost, we delve into almost surely convergence of the sequence generated by Algorithm \ref{Alg1}.
\begin{lemma}
\label{eta3.1}
Suppose that Assumptions $(A1)$-$(A3)$ hold and the sequence $\left\{ {{x_{s}}} \right\}$ is from Algorithm \ref{Alg1}. If $\alpha \le \gamma$ and $0<p<1$, then we have
  \[\begin{array}{*{20}{l}}
{{\mathbb{E}_s}\left[ {D\left( {x',{x_{s + 1}}} \right) + \left( {1 - \gamma } \right)D\left( {x',{x_s}} \right) + \alpha \left\langle {F\left( {{x_{s + 1}}} \right) - F\left( {{w_s}} \right),x' - {x_{s + 1}}} \right\rangle } \right]}\\
{ + {\mathbb{E}_s}\left[ \frac{1}{2}\left( {\frac{{\gamma  - p}}{{1 - p}} - {\alpha }} \right){\left\| {{x_{s + 1}} - {x_{s }}} \right\|^2}+ \frac{{1 - \gamma }}{{2\left( {1 - p} \right)}}{\left\| {{x_{s + 1}} - {w_s}} \right\|^2} \right]}\\
{ \le D\left( {x',{x_s}} \right) + \left( {1 - \gamma } \right)D\left( {x',{x_{s - 1}}} \right) + \alpha \left\langle {F\left( {{x_s}} \right) - F\left( {{w_{s - 1}}} \right),x' - {x_s}} \right\rangle }\\
{ +  \frac{\alpha{{{\bar L}^2}}}{2}{\left\| {{x_s} - {w_{s - 1}}} \right\|^2}.}
\end{array}\]
where ${x'}\in S_1$.
\end{lemma}
\begin{proof}
By the definition of ${\bar x}_{s}$ and $x_{s+1}$, we have
\begin{equation}\label{3.8}
\begin{aligned}
{x_{s + 1}} &= \mathop {\arg \min }\limits_{x \in K} \left\{ {\alpha g\left( x \right) + \left\langle {\alpha {\Delta ^s},x - {x_s}} \right\rangle  + D\left( {x,{{\bar x}_s}} \right)} \right\}\\
 &= \mathop {\arg \min }\limits_{x \in K} \left\{ {\alpha g\left( x \right) + \left\langle {\alpha {\Delta ^s},x - {x_s}} \right\rangle  + f\left( x \right) - f\left( {{{\bar x}_s}} \right) - \left\langle {\nabla f\left( {{{\bar x}_s}} \right),x - {{\bar x}_s}} \right\rangle } \right\}\\
 &= \mathop {\arg \min }\limits_{x \in K} \left\{ {\alpha g\left( x \right) + \left\langle {\alpha {\Delta ^s},x - {x_s}} \right\rangle  + f\left( x \right) - \left\langle {\nabla f\left( {{{\bar x}_s}} \right),x} \right\rangle } \right\}\\
 &= \mathop {\arg \min }\limits_{x \in K} \left\{ {\alpha g\left( x \right) + \left\langle {\alpha {\Delta ^s},x} \right\rangle  + \gamma D\left( {x,{x_s}} \right) + \left( {1 - \gamma } \right)D\left( {x,{x_{s - 1}}} \right)} \right\}.
\end{aligned}
\end{equation}
Plugging the definition of $x_{s+1}$ into (\ref{2.2}) of Lemma \ref{l2}, we have

\begin{equation}
\label{etaaa1}
\begin{aligned}
\alpha g\left( u \right) - \alpha g\left( {{x_{s + 1}}} \right) + \alpha \left\langle {{\Delta ^s},u - {x_{s + 1}}} \right\rangle  &\ge D\left( {u,{x_{s + 1}}} \right) + \gamma \left( {D\left( {{x_{s + 1}},{x_s}} \right) - D\left( {u,{x_s}} \right)} \right)\\
 &+ \left( {1 - \gamma } \right)\left( {D\left( {{x_{s + 1}},{x_{s-1 }}} \right) - D\left( {u,{x_{s-1 }}} \right)} \right),\forall u \in  K.
\end{aligned}
\end{equation}
Using (\ref{etaaa1}), ${x'}\in K$ and the definition of $\Delta^s$, we have
\[\begin{array}{l}
\alpha \left\langle {F\left( {{w_s}} \right) - F\left( {{x_{s + 1}}} \right),x{\rm{'}} - {x_{s + 1}}} \right\rangle  - \alpha \left( {\left\langle {F\left( {{x_{s + 1}}} \right),{x_{s + 1}} - x{\rm{'}}} \right\rangle  + g\left( {{x_{s + 1}}} \right) - g\left( {x'} \right)} \right)\\
 \ge D\left( {x',{x_{s + 1}}} \right) + \gamma \left( {D\left( {{x_{s + 1}},{x_s}} \right) - D\left( {x',{x_s}} \right)} \right) + \left( {1 - \gamma } \right)\left( {D\left( {{x_{s + 1}},{x_{s-1}}} \right) - D\left( {x',{x_{s-1}}} \right)} \right)\\
 + \frac{\alpha }{b}\sum\limits_{\xi  \in {B^s}} {\left\langle {{F_\xi }\left( {{w_{s - 1}}} \right) - {F_\xi }\left( {{x_s}} \right),x' - {x_s}} \right\rangle }  + \frac{\alpha }{b}\sum\limits_{\xi  \in {B^s}} {\left\langle {{F_\xi }\left( {{w_{s - 1}}} \right) - {F_\xi }\left( {{x_s}} \right),{x_s} - {x_{s + 1}}} \right\rangle } .
\end{array}\]
It follows from Assumption $(A2)$ and ${x'}\in S_1$ that
\begin{equation}\label{3.3.1}
\begin{array}{l}
\frac{\alpha }{b}\sum\limits_{\xi  \in {B^s}} {\left\langle {{F_\xi }\left( {{x_s}} \right) - {F_\xi }\left( {{w_{s - 1}}} \right),x' - {x_s}} \right\rangle }  + \frac{\alpha }{b}\sum\limits_{\xi  \in {B^s}} {\left\langle {{F_\xi }\left( {{x_s}} \right) - {F_\xi }\left( {{w_{s - 1}}} \right),{x_s} - {x_{s + 1}}} \right\rangle } \\
 \ge D\left( {x',{x_{s + 1}}} \right) + \gamma \left( {D\left( {{x_{s + 1}},{x_s}} \right) - D\left( {x',{x_s}} \right)} \right) + \left( {1 - \gamma } \right)\left( {D\left( {{x_{s + 1}},{x_{s-1}}} \right) - D\left( {x',{x_{s-1}}} \right)} \right)\\
 + \alpha \left\langle {F\left( {{x_{s + 1}}} \right) - F\left( {{w_s}} \right),x{\rm{'}} - {x_{s + 1}}} \right\rangle
\end{array}
\end{equation}
Using the Young's inequality, we can get
\[\begin{aligned}
\frac{\alpha }{b}\sum\limits_{\xi  \in {B^s}} {\left\langle {{F_\xi }\left( {{x_s}} \right) - {F_\xi }\left( {{w_{s - 1}}} \right),{x_s} - {x_{s + 1}}} \right\rangle }
  \le \sum\limits_{\xi  \in {B^s}} {\frac{\alpha}{{2b}}{{\left\| {{x_s} - {x_{s + 1}}} \right\|}^2} + \frac{{{\alpha}}}{2b}{{\left\| {{F_\xi }\left( {{x_s}} \right) - {F_\xi }\left( {{w_{s - 1}}} \right)} \right\|}^2}} \\
 \le \frac{\alpha}{2}{\left\| {{x_s} - {x_{s + 1}}} \right\|^2} + \frac{{{\alpha }}}{2b}\sum\limits_{\xi  \in {B^s}} {{{\left\| {{F_\xi }\left( {{x_s}} \right) - {F_\xi }\left( {{w_{s - 1}}} \right)} \right\|}^2}} .
\end{aligned}\]
Substituting the above inequality into \eqref{3.3.1}, by (\ref{2.1}), we can conclude that
\begin{equation}\label{3.10'}
\begin{array}{l}
\frac{\alpha }{b}\sum\limits_{\xi  \in {B^s}} {\left\langle {{F_\xi }\left( {{x_s}} \right) - {F_\xi }\left( {{w_{s - 1}}} \right),x' - {x_s}} \right\rangle }  + \frac{{{\alpha }}}{2b}\sum\limits_{\xi  \in {B^s}} {{{\left\| {{F_\xi }\left( {{x_s}} \right) - {F_\xi }\left( {{w_{s - 1}}} \right)} \right\|}^2}} \\
 \ge \left(\gamma-\alpha\right)D\left( {{x_{s + 1},x_s}} \right)+D\left( {x',{x_{s + 1}}} \right) - \gamma D\left( {x',{x_s}} \right) + \left( {1 - \gamma } \right)\left( {D\left( {{x_{s + 1}},{x_{s-1}}} \right) - D\left( {x',{x_{s-1}}} \right)} \right)\\
 + \alpha \left\langle {F\left( {{x_{s + 1}}} \right) - F\left( {{w_s}} \right),x{\rm{'}} - {x_{s + 1}}} \right\rangle .
\end{array}
\end{equation}
Using the definition of $w_{s+1}$, we have
\begin{equation}\label{3.11}
\begin{aligned}
{\mathbb{E}_s}\left[ {\left\| {{x_{s + 1}} - {w_s}} \right\|^2} \right] &={\mathbb{E}_s}\left[ {\mathbb{E}\left[ {{\left\| {{x_{s + 1}} - {w_s}} \right\|^2}\left| {\sigma \left( {{x_0},{x_1}, \cdots ,{x_s},{x_{s + 1}},{w_0},{w_1}, \cdots ,{w_{s - 1}}} \right)} \right.} \right]} \right]\\
&= {\mathbb{E}_s}\left[ {p{\left\| {{x_{s + 1}} - {x_s}} \right\|^2} + \left( {1 - p} \right){\left\| {{x_{s + 1}} - {x_{s-1}}} \right\|^2}} \right].
\end{aligned}
\end{equation}
 Taking conditional expectation on both sides of \eqref{3.10'} and using  (\ref{l3.1.e4}), \eqref{3.11} and Assumption (A3), we can obtain
 \[\begin{array}{*{20}{l}}
{{\mathbb{E}_s}\left[ {D\left( {x',{x_{s + 1}}} \right) + \left( {1 - \gamma } \right)D\left( {x',{x_s}} \right) + \alpha \left\langle {F\left( {{x_{s + 1}}} \right) - F\left( {{w_s}} \right),x' - {x_{s + 1}}} \right\rangle } \right]}\\
{ + {\mathbb{E}_s}\left[ {\left( {\gamma - \alpha } \right)D\left( {{x_{s + 1}},{x_s}} \right)-\frac{{p\left( {1 - \gamma } \right)}}{{2\left( {1 - p} \right)}}{\left\| {{x_{s + 1}} - {x_s}} \right\|^2}} \right]}\\
{ + {\mathbb{E}_s}\left[   \left( {1 - \gamma } \right)\left( {D\left( {{x_{s + 1}},{x_{s - 1}}} \right) - \frac{1}{2}{{\left\| {{x_{s + 1}} - {x_{s - 1}}} \right\|}^2}} \right)+ \frac{{1 - \gamma }}{{2\left( {1 - p} \right)}}{\left\| {{x_{s + 1}} - {w_s}} \right\|^2}\right]}\\
{ \le D\left( {x',{x_s}} \right) + \left( {1 - \gamma } \right)D\left( {x',{x_{s - 1}}} \right) + \alpha \left\langle {F\left( {{x_s}} \right) - F\left( {{w_{s - 1}}} \right),x' - {x_s}} \right\rangle }\\
{ +  \frac{\alpha{{{\bar L}^2}}}{2}{\left\| {{x_s} - {w_{s - 1}}} \right\|^2}.}
\end{array}\]
Owing to (\ref{2.1}), $0<p<1$ and $0<\gamma \le 1$, it yields
 \[\begin{array}{*{20}{l}}
{{\mathbb{E}_s}\left[ {D\left( {x',{x_{s + 1}}} \right) + \left( {1 - \gamma } \right)D\left( {x',{x_s}} \right) + \alpha \left\langle {F\left( {{x_{s + 1}}} \right) - F\left( {{w_s}} \right),x' - {x_{s + 1}}} \right\rangle } \right]}\\
{ + {\mathbb{E}_s}\left[ \frac{1}{2}\left( {\frac{{\gamma  - p}}{{1 - p}} - {\alpha }} \right){\left\| {{x_{s + 1}} - {x_{s }}} \right\|^2}+ \frac{{1 - \gamma }}{{2\left( {1 - p} \right)}}{\left\| {{x_{s + 1}} - {w_s}} \right\|^2} \right]}\\
{ \le D\left( {x',{x_s}} \right) + \left( {1 - \gamma } \right)D\left( {x',{x_{s - 1}}} \right) + \alpha \left\langle {F\left( {{x_s}} \right) - F\left( {{w_{s - 1}}} \right),x' - {x_s}} \right\rangle }\\
{ +  \frac{\alpha{{{\bar L}^2}}}{2}{\left\| {{x_s} - {w_{s - 1}}} \right\|^2}.}
\end{array}\]
We have completed this proof.
\end{proof}
\begin{theorem}
\label{Theorem monotone}
Assume that Assumptions $(A1)$-$(A4)$ hold, $\nabla f$ is continuous and ${{F_\xi }}$ is continuous for all $\xi$. If $0<p<1$,  $p<\gamma<1 - {{\bar L}^2}\left( {1 - p} \right)  $, $\alpha  \le \min \left\{ {\frac{{\gamma  - p}}{{1 - p}} ,\frac{{1 - \gamma }}{{{L^2}\left( {1 - p} \right)}}} \right\}$ and $\alpha<\frac{1-\gamma}{\bar{L}^2\left(1-p\right)}  $,  then, the sequence $\left\{ {{x_s}} \right\}$ generated by Algorithm \ref{Alg1} converges to $\tilde{x} \in S_1$ almost surely.
\end{theorem}
\begin{proof}
Fixing any $x'\in S_1$, Define
\[{c_{s}} = {D\left( {x',{x_s}} \right) + \left( {1 - \gamma } \right)D\left( {x',{x_{s - 1}}} \right) + \alpha \left\langle {F\left( {{x_s}} \right) - F\left( {{w_{s - 1}}} \right),x' - {x_s}} \right\rangle  + \frac{{1 - \gamma }}{{2\left( {1 - p} \right)}}{\left\| {{x_{s }} - {w_{s-1}}} \right\|^2}}.\]
In view of the Young's inequality, Assumption (A4) and (\ref{2.1}) imply that
\[\begin{aligned}
\alpha \left\langle {F\left( {{x_{s}}} \right) - F\left( {{w_{s-1}}} \right),{x'} - {{x_s}}} \right\rangle  &\ge  - \frac{\alpha}{2}{\left\| {F\left( {{w_{s - 1}}} \right) - F\left( {{x_s}} \right)} \right\|_*^2} - \frac{\alpha}{2}{\left\| {{x_s} - {{x'}}} \right\|^2}\\
 &\ge  - \frac{{\alpha }{L^2}}{2}{\left\| {{x_s} - {w_{s - 1}}} \right\|^2}- \alpha D\left( {{{x'}},{x_s}} \right).
\end{aligned}
\]
Hence, it follows from $\alpha  \le \min \left\{ {\frac{{\gamma  - p}}{{1 - p}}=\gamma-\frac{p\left(1-\gamma\right)}{1-p} ,\frac{{1 - \gamma }}{{{L^2}\left( {1 - p} \right)}}} \right\}$ and $p < \gamma  < 1$ that
\begin{equation}\label{etraaa1}
{c_s} \ge \left(1-\alpha\right)D\left( {x',{x_s}} \right) +{\left( 1-\gamma \right)D\left( {x',{x_{s-1}}} \right)}+ \frac{1}{2}\left( {\frac{1-\gamma}{1-p}-\alpha L^2} \right){\left\| {{x_s} - {w_{s - 1}}} \right\|^2}  \ge 0.
\end{equation}
From Lemma \ref{eta3.1} and $\alpha  \le \frac{{\gamma  - p}}{{1 - p}}$, we infer that
\[{{\mathbb{E}}_s}\left[ {{c_{s + 1}}} \right] \le {c_s} + \frac{1}{2}\left( {\alpha{{\bar L}^2} - \frac{{1 - \gamma }}{{1 - p}}} \right){\left\| {{x_s} - {w_{s - 1}}} \right\|^2}.\]
Due to $\alpha<\frac{1-\gamma}{\bar{L}^2\left(1-p\right)}  $ and Lemma \ref{l2'}, we know that
$\{c_s\}$ converges almost surely and $\left\| {{x_s} - {w_{s - 1}}} \right\|$ converges to $0$  almost surely. Furthermore, due to (\ref{etraaa1}), it is easy to see that $\{x_s\}$ is bounded almost surely. Therefore, using the definition of $c_s$ and the continuity of $F$, we can conclude that $D\left( {x',{x_s}} \right) + \left( {1 - \gamma } \right)D\left( {x',{x_{s - 1}}} \right)$ converges almost surely.
\par Owing to the fact that $\left\| {{x_s} - {w_{s - 1}}} \right\|$ converges to $0$ almost surely, there exists probability $1$ set $\Omega $ with $\left\| {{x_s}(w) - {w_{s - 1}}}(w) \right\|\rightarrow 0,~\forall w \in \Omega $. Similarly, there exists probability $1$ set $\Omega'$ such that  $D\left( {x',{x_s(w')}} \right) + \left( {1 - \gamma } \right)D\left( {x',{x_{s - 1}}(w')} \right)$ converges,$~\forall w' \in \Omega'$.
\par Picking $\bar{w} \in \Omega  \cap \Omega '$, we have $\left\| {{x_s}(\bar{w}) - {w_{s - 1}}}(\bar{w}) \right\|\rightarrow 0$, $D\left( {x',{x_s(\bar{w})}} \right) + \left( {1 - \gamma } \right)D\left( {x',{x_{s - 1}}(\bar{w})} \right)$ converges and $\{x_s(\bar{w})\}$ is bounded. Let $\tilde x(\bar{w})$ a cluster point of the bounded sequence $\{x_s(\bar{w})\}$. Using \eqref{l2.1}, \eqref{etaaa1} and the definition of $\Delta^s$, we have
\begin{equation*}
\begin{array}{l}
0 \le \alpha \left\langle {F\left( {{x_{s + 1}}\left( \bar{w} \right)} \right),u - {x_{s + 1}}\left( \bar{w} \right)} \right\rangle  + \alpha g\left( u \right) - \alpha g\left( {{x_{s + 1}}\left( \bar{w} \right)} \right)\\
 + \frac{\alpha }{b}\sum\limits_{\xi  \in {B^s}} {\left\langle {{F_\xi }\left( {{x_s}\left( \bar{w} \right)} \right) - {F_\xi }\left( {{w_{s - 1}}\left( \bar{w} \right)} \right),u - {x_{s + 1}}\left( \bar{w} \right)} \right\rangle }  + \alpha \left\langle {F\left( {{w_s}\left( \bar{w} \right)} \right) - F\left( {{x_{s + 1}}\left( \bar{w} \right)} \right),u - {x_{s + 1}}\left( \bar{w} \right)} \right\rangle \\
+\gamma \left( {D\left( {u,{x_s}\left( {\bar w} \right)} \right) - D\left( {u,{x_{s + 1}}\left( {\bar w} \right)} \right) - D\left( {{x_{s + 1}}\left( {\bar w} \right),{x_s}\left( {\bar w} \right)} \right)} \right)\\
+\left( {1 - \gamma } \right)\left( {D\left( {u,{x_{s - 1}}\left( {\bar w} \right)} \right) - D\left( {u,{x_{s + 1}}\left( {\bar w} \right)} \right) - D\left( {{x_{s + 1}}\left( {\bar w} \right),{x_{s - 1}}\left( {\bar w} \right)} \right)} \right)\\
=\alpha \left\langle {F\left( {{x_{s + 1}}\left( \bar{w} \right)} \right),u - {x_{s + 1}}\left( \bar{w} \right)} \right\rangle  + \alpha g\left( u \right) - \alpha g\left( {{x_{s + 1}}\left( \bar{w} \right)} \right)\\
 + \frac{\alpha }{b}\sum\limits_{\xi  \in {B^s}} {\left\langle {{F_\xi }\left( {{x_s}\left( \bar{w} \right)} \right) - {F_\xi }\left( {{w_{s - 1}}\left( \bar{w} \right)} \right),u - {x_{s + 1}}\left( \bar{w} \right)} \right\rangle }  + \alpha \left\langle {F\left( {{w_s}\left( \bar{w} \right)} \right) - F\left( {{x_{s + 1}}\left( \bar{w} \right)} \right),u - {x_{s + 1}}\left( \bar{w} \right)} \right\rangle \\
+\gamma \left\langle {\nabla f\left( {{x_s}\left( {\bar w} \right)} \right) - \nabla f\left( {{x_{s + 1}}\left( {\bar w} \right)} \right),{x_{s + 1}}\left( {\bar w} \right) - u} \right\rangle \\
+\left( {1 - \gamma } \right)\left( \left\langle {\nabla f\left( {{x_{s-1}}\left( {\bar w} \right)} \right) - \nabla f\left( {{x_{s + 1}}\left( {\bar w} \right)} \right),{x_{s + 1}}\left( {\bar w} \right) - u} \right\rangle \right),\forall u \in K.
\end{array}
\end{equation*}
By extracting the subsequence of $\{x_s(\bar{w})\}$ if needed, taking the limit along that subsequence,  the
lower semicontinuity of $g$ and the continuity of, $\nabla f$, $F$ and ${{F_\xi }}$ imply that
\begin{equation*}
\begin{array}{l}
\left\langle {F\left( {{\tilde x}\left( \bar{w} \right)} \right),u - {\tilde x(\bar{w})}} \right\rangle  +  g\left( u \right) -  g\left( {{\tilde x(\bar{w})}} \right)\geq 0,\forall u \in K.
\end{array}
\end{equation*}
Hence we can infer that
$\tilde x(\bar{w})\in  S_1$. Since $\tilde x(\bar{w})$ a cluster point of the bounded sequence $\{x_s(\bar{w})\}$, there exist a subsequence  $\{x_{s_j}(\bar{w})\}$ of  $\{x_s(\bar{w})\}$ such that $\lim_{j\to\infty}\{x_{s_j}(\bar{w})\}=\tilde x(\bar{w})$. By the definition of Bregman distance and the boundedness of $\{x_s(\bar{w})\}$, we can show that  $D\left( {\tilde x(\bar{w}),{x_{s_j}(\bar{w})}} \right) + \left( {1 - \gamma } \right)D\left( {\tilde x(\bar{w}),{x_{s_j - 1}}(\bar{w})} \right)$ converges to $0$. Note that for any $x'\in S_1$, $D\left( {x',{x_s(\bar{w})}} \right) + \left( {1 - \gamma } \right)D\left( {x',{x_{s - 1}}(\bar{w})} \right)$ converges. Therefore,  $D\left( {\tilde x(\bar{w}),{x_s(\bar{w})}} \right) + \left( {1 - \gamma } \right)D\left( {\tilde x(\bar{w}),{x_{s - 1}}(\bar{w})} \right)\rightarrow 0,~\forall ~\bar{w} \in \Omega  \cap \Omega '$. Due to \eqref{2.1}, we can get $${\left\| {\tilde x\left( {\bar w} \right) - {x_s}\left( {\bar w} \right)} \right\|^2} + \left( {1 - \gamma } \right){\left\| {\tilde x\left( {\bar w} \right) - {x_{s - 1}}\left( {\bar w} \right)} \right\|^2} \to 0.$$ So, the sequence $\left\{ {{x_s}} \right\}$ converges to $\tilde{x} \in S_1$ almost surely.
 The proof is completed.
\end{proof}
 \par When $p=1$, that is, $w_s=x_s$, we can also establish the sequence $\left\{ {{x_s}} \right\}$ generated by Algorithm 1 converges almost surely. For details, see the following theorem.
 \begin{theorem}
\label{Theorem monotone'}
Assume that Assumptions $(A1)$-$(A4)$ hold, ${{F_\xi }}$ is continuous for all $\xi$ and the gradient $\nabla f$ of $f$ is continuous. If $p=1$, $\alpha  < \frac{{\gamma }}{{1+{\bar L}^2}}$ and $\alpha  \le \frac{{\sqrt {1 + 8\gamma {L^2}}  - 1}}{{4{L^2}}}$, then the sequence $\left\{ {{x_s}} \right\}$ from Algorithm \ref{Alg1} converges to $\tilde{x} \in S_1$ almost surely.
\end{theorem}
\begin{proof}
Using the fact $w_s=x_s$ in \eqref{3.10'}, we can obtain
\[\begin{array}{l}
D\left( {x',{x_{s + 1}}} \right) + \left( {1 - \gamma } \right)D\left( {x',{x_s}} \right) + \left( {\gamma  - \alpha} \right)D\left( {{x_{s + 1}},{x_s}} \right) + \alpha \left\langle {F\left( {{x_s}} \right) - F\left( {{x_{s + 1}}} \right),{x_{s + 1}} - x'} \right\rangle \\
 \le D\left( {x',{x_s}} \right) + \left( {1 - \gamma } \right)D\left( {x',{x_{s - 1}}} \right) + \frac{\alpha }{b}\sum\limits_{\xi  \in {B^s}} {\left\langle {{F_\xi }\left( {{x_s}} \right) - {F_\xi }\left( {{x_{s - 1}}} \right),x' - {x_s}} \right\rangle } \\
 + \frac{{{\alpha}}}{2b}\sum\limits_{\xi  \in {B^s}} {{{\left\| {{F_\xi }\left( {{x_s}} \right) - {F_\xi }\left( {{x_{s - 1}}} \right)} \right\|}^2_*}}.
\end{array}\]
 Taking conditional expectation on both sides of the above inequality, and using (\ref{2.1}), (\ref{l3.1.e4}) and Assumption (A3), we can infer that
 \begin{equation}\label{eta3.1'}
\begin{array}{l}
{\mathbb{E}_s}\left[ {D\left( {x',{x_{s + 1}}} \right) + \left( {1 - \gamma } \right)D\left( {x',{x_s}} \right) + \left( {\gamma  - \alpha} \right)D\left( {{x_{s + 1}},{x_s}} \right) + \alpha \left\langle {F\left( {{x_{s + 1}}} \right) - F\left( {{x_s}} \right),x' - {x_{s + 1}}} \right\rangle } \right]\\
 \le D\left( {x',{x_s}} \right) + \left( {1 - \gamma } \right)D\left( {x',{x_{s - 1}}} \right) + \alpha \left\langle {F\left( {{x_s}} \right) - F\left( {{x_{s - 1}}} \right),x' - {x_s}} \right\rangle  + {\alpha }{{\bar L}^2}D\left( {{x_s},{x_{s - 1}}} \right).
\end{array}
 \end{equation}
Fix any $x'\in S_1$ and Set
\[{c'_s} = {D\left( {x',{x_{s }}} \right) + \left( {1 - \gamma } \right)D\left( {x',{x_{s-1}}} \right) + \left( {\gamma  - \alpha} \right)D\left( {{x_{s }},{x_{s-1}}} \right) + \alpha \left\langle {F\left( {{x_{s}}} \right) - F\left( {{x_{s-1}}} \right),x' - {x_{s}}} \right\rangle }.\]
In view of the Young's inequality, Assumption (A4) and (\ref{2.1}) imply that
\[\begin{aligned}
\alpha \left\langle {F\left( {{x_{s}}} \right) - F\left( {{x_{s-1}}} \right),{x'} - {{x_s}}} \right\rangle  &\ge  - {\alpha ^2}{\left\| {F\left( {{x_{s - 1}}} \right) - F\left( {{x_s}} \right)} \right\|_*^2} - \frac{1}{4}{\left\| {{x_s} - {{x'}}} \right\|^2}\\
 &\ge  - 2{\alpha ^2}{L^2}D\left( {{x_s},{x_{s - 1}}} \right) - \frac{1}{2}D\left( {{{x'}},{x_s}} \right).
\end{aligned}
\]
Hence, it follows from $\alpha  \le \frac{{\sqrt {1 + 8\gamma {L^2}}  - 1}}{{4{L^2}}}$ and $0<\gamma \le 1$ that
\[
{c'_s} \ge \frac{1}{2}D\left( {x',{x_s}} \right) +{\left( {1 - \gamma } \right)D\left( {x',{x_{s-1}}} \right)}+ \left( {\gamma-\alpha - 2{\alpha ^2}{L^2}} \right)D\left( {{x_s},{x_{s - 1}}} \right) \ge \frac{1}{2}D\left( {x',{x_s}} \right) \ge 0.
\]
From \eqref{eta3.1'}, we infer that
\[{{\mathbb{E}}_s}\left[ {{c'_{s + 1}}} \right] \le {c'_s} - \left( {\gamma-\alpha -{{{\alpha }{\bar L}^2}}} \right)D\left( {{x_s},{x_{s - 1}}} \right).\]
 Using the similar arguments as in the proof of Theorem \ref{Theorem monotone}, we can show that the sequence $\left\{ {{x_s}} \right\}$ converges to $\tilde{x} \in S_1$ almost surely. Hence, the proof is completed.
\end{proof}
\begin{remark}
If the operator $F$ is $g$-pseudomonotone other than  monotone, then using the similar arguments as in Theorem \ref{Theorem monotone} and Theorem \ref{Theorem monotone'},  we can also assure that the sequence $\left\{ {{x_s}} \right\}$ generated by Algorithm \ref{Alg1} converges to $\tilde{x} \in S_1$ almost surely.
\end{remark}
We now study the complexity of Algorithm \ref{Alg1} when $F$ is monotone.
In order to obtain a ergodic convergence rate of Algorithm 1, Define the following gap function
\[{\rm{Gap}}\left( z \right): = \mathop {\max }\limits_{u \in C} \left\{ {\left\langle {F\left( u \right),z - u} \right\rangle +g\left( z \right) - g\left( u \right)} \right\},\]
where $C$ is a compact subset of $K$. This gap function has been used in \cite{AA,AA2,AZ,PA1,PA}.
\begin{lemma}
\label{al2.1}
(see Lemma 2.4 from \cite{AA})
Let $\mathcal{F}=({\mathcal{F}}_s)_{s \ge 0}$ be a filtration, the stochastic process $\{u_s\}$ be adapted to $\mathcal{F}$  and  ${\mathbb{E}}[u_{s+1}|{\mathcal{F}}_s]=0$. Then for any $S \in \mathbb{N}$, $x_0 \in K$, and any compact set $C \subset K$,
\[{\mathbb{E}}\left[ {\mathop {\max }\limits_{x \in C} \sum\limits_{s = 0}^{S - 1} {\left\langle {{u_{s + 1}},x} \right\rangle } } \right] \le \mathop {\max }\limits_{x \in C} \frac{1}{2}{\left\| {{x_0} - x} \right\|^2} + \frac{1}{2}\sum\limits_{s = 0}^{S - 1} {{\mathbb{E}}\left[ {{{\left\| {{u_{s + 1}}} \right\|}^2}} \right]}. \]
\end{lemma}

Now, we present the theorem to state the convergence rate of Algorithm \ref{Alg1}.
\begin{theorem}
\label{T3.1}
Let Assumptions $(A1)-(A4)$ hold, the sequence $\left\{ {{x_s}} \right\}$ be from Algorithm \ref{Alg1}, $p \in \left(0,1\right)$, $\gamma \in \left(p,1\right)$ and $\alpha  \le \min \left\{ {\frac{{\gamma  - p}}{{2\left( {1 - p} \right)}},\frac{{\left( {1 - \gamma} \right)b}}{{\left( {1 - p } \right)\left( {2{{\bar L}^2} + b{L^2}} \right)}}} \right\}$.
Thus
\[{\mathbb{E}}\left[ {{\rm{Gap}}\left( {{z^S}} \right)} \right] \le \frac{\left(2+\alpha-\gamma\right)}{{\alpha S}}\mathop {\max }\limits_{u \in C} D\left( {u,{x_0}} \right),\]
where ${z^S} = \frac{1}{S}\sum\limits_{s = 0}^{S - 1} {{x_{s + 1}}} $.
\end{theorem}
\begin{proof}
Thanks to \eqref{etaaa1}, we have
\[\begin{aligned}
\alpha g\left( u \right) - \alpha g\left( {{x_{s + 1}}} \right) + \alpha \left\langle {{\Delta ^s},u - {x_{s + 1}}} \right\rangle  &\ge D\left( {u,{x_{s + 1}}} \right) + \gamma \left( {D\left( {{x_{s + 1}},{x_s}} \right) - D\left( {u,{x_s}} \right)} \right)\\
 &+ \left( {1 - \gamma } \right)\left( {D\left( {{x_{s + 1}},{x_{s-1 }}} \right) - D\left( {u,{x_{s-1 }}} \right)} \right),\forall u \in C \subset K.
\end{aligned}\]
Rearranging the order, it yields
\begin{equation}
\label{etaa1}
\begin{aligned}
D\left( {u,{x_{s + 1}}} \right) &\le \alpha g\left( u \right) - \alpha g\left( {{x_{s + 1}}} \right) + \alpha \left\langle {{\Delta ^s},u - {x_{s + 1}}} \right\rangle \\
 &+ \gamma \left( {D\left( {u,{x_s}} \right) - D\left( {{x_{s + 1}},{x_s}} \right)} \right) + \left( {1 - \gamma } \right)\left( {D\left( {u,{x_{s-1 }}} \right)-D\left( {{x_{s + 1}},{x_{s-1}}} \right)} \right),\forall u \in C \subset K.
\end{aligned}
\end{equation}
From \eqref{l3.1.e1} and \eqref{etaa1}, it follows that
\[\begin{aligned}
D\left( {u,{x_{s + 1}}} \right) &\le  - \alpha \left\langle {{{\mathbb{E}}_s}\left[ {{\Delta ^s}} \right],{x_{s + 1}} - u} \right\rangle  + \alpha \left\langle {{{\mathbb{E}}_s}\left[ {{\Delta ^s}} \right] - {\Delta ^s},{x_{s + 1}} - {x_s}} \right\rangle \\
 &+ \alpha \left\langle {{{\mathbb{E}}_s}\left[ {{\Delta ^s}} \right] - {\Delta ^s},{x_s} - u} \right\rangle   + \alpha g\left( u \right) - \alpha g\left( {{x_{s + 1}}} \right)\\
  &+ \gamma \left( {D\left( {u,{x_s}} \right) - D\left( {{x_{s + 1}},{x_s}} \right)} \right) + \left( {1 - \gamma } \right)\left( {D\left( {u,{x_{s-1 }}} \right)-D\left( {{x_{s + 1}},{x_{s-1 }}} \right)} \right)\\
 &=  - \alpha \left\langle {F\left( {{w_s}} \right) + F\left( {{x_s}} \right) - F\left( {{w_{s - 1}}} \right),{x_{s + 1}} - u} \right\rangle  + \alpha \left\langle {{{\mathbb{E}}_s}\left[ {{\Delta ^s}} \right] - {\Delta ^s},{x_{s + 1}} - {x_s}} \right\rangle \\
 & + \gamma \left( {D\left( {u,{x_s}} \right) - D\left( {{x_{s + 1}},{x_s}} \right)} \right) + \left( {1 - \gamma } \right)\left( {D\left( {u,{x_{s-1 }}} \right)-D\left( {{x_{s + 1}},{x_{s-1 }}} \right)} \right)\\
 & + \alpha \left\langle {{{\mathbb{E}}_s}\left[ {{\Delta ^s}} \right] - {\Delta ^s},{x_s} - u} \right\rangle+ \alpha g\left( u \right) - \alpha g\left( {{x_{s + 1}}} \right)\\
 &=  - \alpha \left\langle {F\left( {{w_s}} \right) - F\left( {{x_{s + 1}}} \right) + F\left( {{x_s}} \right) - F\left( {{w_{s - 1}}} \right),{x_{s + 1}} - u} \right\rangle  - \alpha \left\langle {F\left( {{x_{s + 1}}} \right),{x_{s + 1}} - u} \right\rangle \\
 &+ \alpha \left\langle {{{\mathbb{E}}_s}\left[ {{\Delta ^s}} \right] - {\Delta ^s},{x_{s + 1}} - {x_s}} \right\rangle  + \alpha \left\langle {{{\mathbb{E}}_s}\left[ {{\Delta ^s}} \right] - {\Delta ^s},{x_s} - u} \right\rangle+ \alpha g\left( u \right) - \alpha g\left( {{x_{s + 1}}} \right) \\
 &+ \gamma \left( {D\left( {u,{x_s}} \right) - D\left( {{x_{s + 1}},{x_s}} \right)} \right) + \left( {1 - \gamma } \right)\left( {D\left( {u,{x_{s-1 }}} \right)-D\left( {{x_{s + 1}},{x_{s-1 }}} \right)} \right).
\end{aligned}\]
After some rearrangements, we derive
\begin{equation}\label{t1.1}
\begin{array}{l}
\alpha g\left( {{x_{s + 1}}} \right) - \alpha g\left( u \right) + \alpha \left\langle {F\left( {{x_{s + 1}}} \right),{x_{s + 1}} - u} \right\rangle \\
 \le  - \alpha \left\langle {F\left( {{w_s}} \right) - F\left( {{x_{s + 1}}} \right) + F\left( {{x_s}} \right) - F\left( {{w_{s - 1}}} \right),{x_{s + 1}} - u} \right\rangle
 \\+ \gamma \left( {D\left( {u,{x_s}} \right) - D\left( {{x_{s + 1}},{x_s}} \right)} \right) + \left( {1 - \gamma } \right)\left( {D\left( {u,{x_{s-1 }}} \right)-D\left( {{x_{s + 1}},{x_{s-1 }}} \right)} \right)  \\
 + \alpha \left\langle {{{\mathbb{E}}_s}\left[ {{\Delta ^s}} \right] - {\Delta ^s},{x_{s + 1}} - {x_s}} \right\rangle  + \alpha \left\langle {{{\mathbb{E}}_s}\left[ {{\Delta ^s}} \right] - {\Delta ^s},{x_s} - u} \right\rangle-D\left( {u,{x_{s + 1}}} \right) \\
 =  - \alpha \left\langle {F\left( {{w_s}} \right) - F\left( {{x_{s + 1}}} \right),{x_{s + 1}} - u} \right\rangle  + \alpha \left\langle {F\left( {{w_{s-1}}} \right) - F\left( {{x_{s }}} \right),{x_s} - u} \right\rangle \\
 - \alpha \left\langle {F\left( {{x_s}} \right) - F\left( {{w_{s - 1}}} \right),{x_{s + 1}} - {x_s}} \right\rangle  + \alpha \left\langle {{{\mathbb{E}}_s}\left[ {{\Delta ^s}} \right] - {\Delta ^s},{x_{s + 1}} - {x_s}} \right\rangle + \alpha \left\langle {{{\mathbb{E}}_s}\left[ {{\Delta ^s}} \right] - {\Delta ^s},{x_s} - u} \right\rangle\\
  + \gamma \left( {D\left( {u,{x_s}} \right) - D\left( {{x_{s + 1}},{x_s}} \right)} \right) + \left( {1 - \gamma } \right)\left( {D\left( {u,{x_{s-1 }}} \right)-D\left( {{x_{s + 1}},{x_{s-1 }}} \right)} \right) -D\left( {u,{x_{s + 1}}} \right).
\end{array}
\end{equation}
By Young's inequality and Assumption (A4), we have
\begin{equation}\label{t1.2}
\begin{aligned}
 - \alpha \left\langle {F\left( {{x_s}} \right) - F\left( {{w_{s - 1}}} \right),{x_{s + 1}} - {x_s}} \right\rangle  &\le \frac{{\alpha}}{2}{\left\| {F\left( {{x_s}} \right) - F\left( {{w_{s - 1}}} \right)} \right\|_*^2} + \frac{{\alpha}}{2}{\left\| {{x_{s + 1}} - {x_s}} \right\|^2}\\
 &\le \frac{{\alpha L^2}}{2}{\left\| {{x_{s}} - {w_{s-1}}} \right\|^2} + \frac{{\alpha}}{2}{\left\| {{x_{s + 1}} - {x_s}} \right\|^2}.
\end{aligned}
\end{equation}
Substituting (\ref{t1.2}) into (\ref{t1.1}), we can get
\begin{equation}
\label{3.10}
\begin{array}{l}
\alpha g\left( {{x_{s + 1}}} \right) - \alpha g\left( u \right) + \alpha \left\langle {F\left( {{x_{s + 1}}} \right),{x_{s + 1}} - u} \right\rangle \\
 \le  - \alpha \left\langle {F\left( {{w_s}} \right) - F\left( {{x_{s + 1}}} \right),{x_{s + 1}} - u} \right\rangle  +\alpha \left\langle {F\left( {{w_{s-1}}} \right) - F\left( {{x_{s }}} \right),{x_s} - u} \right\rangle \\
 + \frac{{\alpha}L^2}{2}{\left\| {{x_{s}} - {w_{s-1}}} \right\|^2} + \frac{{\alpha}}{2}{\left\| {{x_{s + 1}} - {x_s}} \right\|^2} + \alpha \left\langle {{{\mathbb{E}}_s}\left[ {{\Delta ^s}} \right] - {\Delta ^s},{x_{s + 1}} - {x_s}} \right\rangle + \alpha \left\langle {{{\mathbb{E}}_s}\left[ {{\Delta ^s}} \right] - {\Delta ^s},{x_s} - u} \right\rangle \\
  + \gamma \left( {D\left( {u,{x_s}} \right) - D\left( {{x_{s + 1}},{x_s}} \right)} \right) + \left( {1 - \gamma } \right)\left( {D\left( {u,{x_{s-1 }}} \right)-D\left( {{x_{s + 1}},{x_{s-1 }}} \right)} \right)-D\left( {u,{x_{s + 1}}} \right) .
\end{array}
\end{equation}
Owing to $x_0=w_{-1}$, we have
\begin{equation}\label{extral1}
\begin{array}{l}
\sum\limits_{s = 0}^{S - 1} { - \alpha \left\langle {F\left( {{w_s}} \right) - F\left( {{x_{s + 1}}} \right),{x_{s + 1}} - u} \right\rangle  + \alpha \left\langle {F\left( {{w_{s-1}}} \right) - F\left( {{x_{s}}} \right),{x_s} - u} \right\rangle } \\
 = \alpha \left\langle {F\left( {{x_S}} \right) - F\left( {{w_{S - 1}}} \right),{x_S} - u} \right\rangle.
\end{array}
\end{equation}
Thanks to \eqref{3.10} and  \eqref{extral1}, we know that
\begin{equation}
\label{t1.3}
\begin{array}{l}
\alpha \sum\limits_{s = 0}^{S - 1} {\left[ {g\left( {{x_{s + 1}}} \right) - g\left( u \right) + \left\langle {F\left( {{x_{s + 1}}} \right),{x_{s + 1}} - u} \right\rangle } \right]} \\
 \le \alpha \left\langle {F\left( {{x_S}} \right) - F\left( {{w_{S - 1}}} \right),{x_S} - u} \right\rangle   \\
 + \sum\limits_{s = 0}^{S - 1} {\left[ {\frac{{\alpha}L^2}{2}{\left\| {{x_{s}} - {w_{s-1}}} \right\|^2} + \frac{{\alpha}}{2}{{\left\| {{x_{s + 1}} - {x_s}} \right\|}^2} + \alpha \left\langle {{{\mathbb{E}}_s}\left[ {{\Delta ^s}} \right] - {\Delta ^s},{x_{s + 1}} - {x_s}} \right\rangle+\alpha \left\langle {{{\mathbb{E}}_s}\left[ {{\Delta ^s}} \right] - {\Delta ^s},{x_s} - u} \right\rangle   } \right]} \\
 + \sum\limits_{s = 0}^{S - 1} {\left[ { \gamma \left( {D\left( {u,{x_s}} \right) - D\left( {{x_{s + 1}},{x_s}} \right)} \right) + \left( {1 - \gamma } \right)\left( {D\left( {u,{x_{s-1 }}} \right)-D\left( {{x_{s + 1}},{x_{s-1 }}} \right)} \right)-D\left( {u,{x_{s + 1}}} \right) } \right]}
\end{array}
\end{equation}
 Taking maximum and expectation on both sides of \eqref{t1.3}, we can obtain
\begin{equation}\label{t1.4}
\begin{array}{*{20}{l}}
{\mathbb{E}\left[ {\mathop {\max }\limits_{u \in C} \left\{ {\alpha \sum\limits_{s = 0}^{S - 1} {\left[ {g\left( {{x_{s + 1}}} \right) - g\left( u \right) + \left\langle {F\left( {{x_{s + 1}}} \right),{x_{s + 1}} - u} \right\rangle } \right]} } \right\}} \right]}\\
\begin{array}{l}
 \le \mathbb{E}\left[ {\mathop {\max }\limits_{u \in C} \left\{ {\sum\limits_{s = 0}^{S - 1} {\left[ {\alpha \left\langle {{\mathbb{E}_s}\left[ {{\Delta ^s}} \right] - {\Delta ^s},{x_s} - u} \right\rangle } \right]} } \right\}} \right]\\
 +\mathbb{E}\left[ {\mathop {\max }\limits_{u \in C} \left\{ {\sum\limits_{s = 0}^{S - 1} {\left[ \gamma D\left( {u,{x_s}} \right) + \left( {1 - \gamma } \right)\left( {D\left( {u,{x_{s-1}}} \right) - D\left( {{x_{s + 1}},{x_{s-1}}} \right)} \right) - D\left( {u,{x_{s + 1}}} \right) \right]} } \right\}} \right]\\
 + \mathbb{E}\left[ {\mathop {\max }\limits_{u \in C} \left\{ {\alpha \left\langle {F\left( {{x_S}} \right) - F\left( {{w_{S - 1}}} \right),{x_S} - u} \right\rangle } \right\}} \right]
\end{array}\\
{ + \mathbb{E}\left[ {\sum\limits_{s = 0}^{S - 1} {\left[ {\frac{{\alpha}L^2}{2}{\left\| {{x_{s}} - {w_{s-1}}} \right\|^2} + \frac{{\alpha}}{2}{{\left\| {{x_{s + 1}} - {x_s}} \right\|}^2} + \alpha \left\langle {{\mathbb{E}_s}\left[ {{\Delta ^s}} \right] - {\Delta ^s},{x_{s + 1}} - {x_s}} \right\rangle  -\gamma D\left( {{x_{s + 1}},{x_s}} \right)} \right]} } \right].}
\end{array}
\end{equation}
According to properties of the conditional expectation, we have
 \begin{equation}\label{var1}
\mathbb{E}\left[ {\left\langle {{\mathbb{E}_s}\left[ {{\Delta ^s}} \right] - {\Delta ^s},{x_s}} \right\rangle } \right] = \mathbb{E}\left[ {{\mathbb{E}_s}\left[ {\left\langle {{\mathbb{E}_s}\left[ {{\Delta ^s}} \right] - {\Delta ^s},{x_s}} \right\rangle } \right]} \right] = 0,\;\forall \;s = 0,...,S - 1.
 \end{equation}
 This and Lemma \ref{al2.1} imply that
\begin{equation}\label{t1.5}
\begin{aligned}
{\mathbb{E}}\left[ {\mathop {\max }\limits_{u \in C} \left\{ {\sum\limits_{s = 0}^{S - 1} {\left[ {\alpha \left\langle {{{\mathbb{E}}_s}\left[ {{\Delta ^s}} \right] - {\Delta ^s},{x_s} - u} \right\rangle } \right]} } \right\}} \right] &= {\mathbb{E}}\left[ {\mathop {\max }\limits_{u \in C} \left\{ {\sum\limits_{s = 0}^{S - 1} {\left[ {\alpha \left\langle { - \left( {{{\mathbb{E}}_s}\left[ {{\Delta ^s}} \right] - {\Delta ^s}} \right),u} \right\rangle } \right]} } \right\}} \right]\\
 &\le \frac{\alpha}{2}\mathop {\max }\limits_{u \in C} {\left\| {u - {x_0}} \right\|^2} + \frac{{{\alpha }}}{2}\sum\limits_{s = 0}^{S - 1} {{\mathbb{E}}\left[ {{{\left\| {{{\mathbb{E}}_s}\left[ {{\Delta ^s}} \right] - {\Delta ^s}} \right\|}_*^2}} \right]}.
\end{aligned}
\end{equation}
Next, using again the Young's inequality, we have
 \begin{equation}\label{t1.6}
\alpha {\mathbb{E}}\left[ {\sum\limits_{s = 0}^{S - 1} {\left[ {\left\langle {{{\mathbb{E}}_s}\left[ {{\Delta ^s}} \right] - {\Delta ^s},{x_{s + 1}} - {x_s}} \right\rangle } \right]} } \right] \le \frac{{\alpha}}{2}\sum\limits_{s = 0}^{S - 1} {{\mathbb{E}}\left[ {{{\left\| {{{\mathbb{E}}_s}\left[ {{\Delta ^s}} \right] - {\Delta ^s}} \right\|}_*^2}} \right]}  + \frac{{\alpha}}{2}\sum\limits_{s = 0}^{S - 1} {{\mathbb{E}}\left[ {{{\left\| {{x_{s + 1}} - {x_s}} \right\|}^2}} \right]},
 \end{equation}
 and
 \begin{equation}\label{t1.7}
\begin{aligned}
\alpha \left\langle {F\left( {{x_S}} \right) - F\left( {{w_{S - 1}}} \right),{x_S} - u} \right\rangle  &\le \frac{{{\alpha }}}{2}{\left\| {F\left( {{x_S}} \right) - F\left( {{w_{S - 1}}} \right)} \right\|_*^2} + \frac{\alpha}{2}{\left\| {{x_S} - u} \right\|^2}\\
 &\le  \frac{{{\alpha }{L^2}}}{2}{\left\| {{x_S} - w_{S-1}} \right\|^2}+ \frac{\alpha}{2}{\left\| {{x_S} - u} \right\|^2},
\end{aligned}
 \end{equation}
 where the last inequality is from Assumption (A4).\\
Substituting (\ref{t1.5}), (\ref{t1.6}) and (\ref{t1.7}) into (\ref{t1.4}), we have
\begin{equation}\label{t1.8}
\begin{array}{*{20}{l}}
{\mathbb{E}\left[ {\mathop {\max }\limits_{u \in C} \left\{ {\alpha \sum\limits_{s = 0}^{S - 1} {\left[ {g\left( {{x_{s + 1}}} \right) - g\left( u \right) + \left\langle {F\left( {{x_{s + 1}}} \right),{x_{s + 1}} - u} \right\rangle } \right]} } \right\}} \right]}\\
{\begin{array}{*{20}{l}}
{ \le \alpha \mathop {\max }\limits_{u \in C} D\left( {u,{x_0}} \right) + \alpha\mathbb{E}\left[ {\sum\limits_{s = 0}^{S - 1} {\left[ {\left\| {{\mathbb{E}_s}\left[ {{\Delta ^s}} \right] - {\Delta ^s}} \right\|_*^2} \right]} } \right]}\\
{ + \mathbb{E}\left[ {\mathop {\max }\limits_{u \in C} \left\{ {\sum\limits_{s = 0}^{S - 1} {\left[ {D\left( {u,{x_s}} \right) - D\left( {u,{x_{s + 1}}} \right) + \left( {1 - \gamma } \right)\left( {D\left( {u,{x_{s - 1}}} \right) - D\left( {u,{x_s}} \right)} \right)} \right]} } \right\}} \right]}\\
{ + \mathbb{E}\left[ {\mathop {\max }\limits_{u \in C} \left\{ {\frac{{\alpha {L^2}}}{2}{{\left\| {{x_S} - {w_{S - 1}}} \right\|}^2} + \frac{\alpha }{2}{{\left\| {{x_S} - u} \right\|}^2}} \right\}} \right]}
\end{array}}\\
{ + \mathbb{E}\left[ {\sum\limits_{s = 0}^{S - 1} {\left[ {\frac{{\alpha {L^2}}}{2}{{\left\| {{x_s} - {w_{s - 1}}} \right\|}^2} + \alpha {{\left\| {{x_{s + 1}} - {x_s}} \right\|}^2} - \gamma D\left( {{x_{s + 1}},{x_s}} \right)-\left(1-\gamma\right)D\left( {{x_{s + 1}},{x_{s-1}}} \right)} \right]} } \right].}
\end{array}
\end{equation}

  Taking expectation on both sides of \eqref{3.11}, it yields
\begin{equation}\label{t1.9}
\mathbb{E}\left[ {\left\| {{x_{s + 1}} - {w_s}} \right\|^2} \right]  = p\mathbb{E}\left[ {\left\| {{x_{s + 1}} - {x_s}} \right\|^2} \right] + \left( {1 - p} \right)\mathbb{E}\left[ {\left\| {{x_{s + 1}} - {x_{s-1}}} \right\|^2} \right],~\forall s=0,\ldots,S-1.
\end{equation}
From \eqref{t1.8}, Lemma \ref{l3.1}, $x_0=x_{-1}=w_{-1}$, $\alpha  \le \min \left\{ {\frac{{\gamma  - p}}{{2\left( {1 - p} \right)}},\frac{{\left( {1 - \gamma} \right)b}}{{\left( {1 - p } \right)\left( {2{{\bar L}^2} + b{L^2}} \right)}}} \right\}$,  \eqref{2.1} and \eqref{t1.9}, it follows that
\[\begin{array}{*{20}{l}}
{\mathbb{E}\left[ {\mathop {\max }\limits_{u \in C} \left\{ {\alpha \sum\limits_{s = 0}^{S - 1} {\left[ {g\left( {{x_{s + 1}}} \right) - g\left( u \right) + \left\langle {F\left( {{x_{s + 1}}} \right),{x_{s + 1}} - u} \right\rangle } \right]} } \right\}} \right]}\\
{\begin{array}{*{20}{l}}
{ \le \left( {2 + \alpha  - \gamma } \right)\mathop {\max }\limits_{u \in C} D\left( {u,{x_0}} \right) + \left( {\frac{{\alpha {{\bar L}^2}}}{b} + \frac{\alpha {L^2}}{2}} \right)\mathbb{E}\left[ {\sum\limits_{s = 0}^{S - 1} {\left[ {\left\| {{x_s} - {w_{s - 1}}} \right\|^2} \right]} } \right]}\\
{ + \mathbb{E}\left[ {\mathop {\max }\limits_{u \in C} \left\{ { - \left( {1 - \alpha } \right)D\left( {u,{x_S}} \right) - \left( {1 - \gamma } \right)D\left( {u,{x_{S - 1}}} \right)} \right\}} \right]}
\end{array}}\\
\begin{array}{l}
 + \frac{\alpha {L^2}}{2}\mathbb{E}\left[ {\left\| {{x_S} - {w_{S - 1}}} \right\|^2} \right] + \mathbb{E}\left[ {\sum\limits_{s = 0}^{S - 1} {\left[ { - \frac{\left( {\gamma  - 2\alpha } \right)}{2}{\left\| {{x_{s+1}} - {x_{s}}} \right\|^2} - \frac{\left( {1 - \gamma } \right)}{2}{\left\| {{x_{s+1}} - {x_{s - 1}}} \right\|^2}} \right]} } \right]\\
 \leq\left( {2 + \alpha  - \gamma } \right)\mathop {\max }\limits_{u \in C} D\left( {u,{x_0}} \right)+ \frac{\alpha {L^2}}{2}\mathbb{E}\left[ {\left\| {{x_S} - {w_{S - 1}}} \right\|^2} \right] \\
 + \left( {\frac{{\alpha {{\bar L}^2}}}{b} + \frac{\alpha {L^2}}{2}} \right)\mathbb{E}\left[ {\sum\limits_{s = 0}^{S - 1} {\left[ {\left\| {{x_s} - {w_{s - 1}}} \right\|^2} \right]} } \right] - \frac{{1 - \gamma }}{2\left({1 - p}\right)}\mathbb{E}\left[ {\sum\limits_{s = 0}^{S - 1} {\left[ {\left\| {{x_{s+1}} - {w_{s}}} \right\|^2} \right]} } \right]\\
 + \mathbb{E}\left[ {\sum\limits_{s = 0}^{S - 1} {\left[ {\frac{{1 - \gamma }}{2\left({1 - p}\right)}{\left\| {{x_{s + 1}} - {w_s}} \right\|^2} { - \frac{\left( {\gamma  - 2\alpha } \right)}{2}{\left\| {{x_{s+1}} - {x_{s}}} \right\|^2} - \frac{\left( {1 - \gamma } \right)}{2}{\left\| {{x_{s+1}} - {x_{s - 1}}} \right\|^2}}} \right]} } \right]\\
 \le \left( {2 + \alpha  - \gamma } \right)\mathop {\max }\limits_{u \in C} D\left( {u,{x_0}} \right) + \left( {\frac{{\alpha {{\bar L}^2}}}{b} + \frac{\alpha {L^2}}{2}} \right)\mathbb{E}\left[ {\sum\limits_{s = 0}^{S - 1} {\left[ {\left\| {{x_s} - {w_{s - 1}}} \right\|^2} - {\left\| {{x_{s + 1}} - {w_s}} \right\|^2} \right]} } \right]\\
 +\frac{\alpha {L^2}}{2}\mathbb{E}\left[ {\left\| {{x_S} - {w_{S - 1}}} \right\|^2} \right] - \mathbb{E}\left[ {\sum\limits_{s = 0}^{S - 1} {\left[ {\frac{1}{2}\left( {\gamma  - 2\alpha  - \frac{{\left( {1 - \gamma } \right)p}}{{1 - p}}} \right){\left\| {{x_{s + 1}} - {x_s}} \right\|^2}} \right]} } \right]\\
 = \left( {2 + \alpha  - \gamma } \right)\mathop {\max }\limits_{u \in C} D\left( {u,{x_0}} \right) -\frac{\alpha {\bar{L}^2}}{b}\mathbb{E}\left[ {\left\| {{x_S} - {w_{S - 1}}} \right\|^2} \right] - \mathbb{E}\left[ {\sum\limits_{s = 0}^{S - 1} {\left[ {\frac{1}{2}\left( {\frac{\gamma-p}{1-p}  - 2\alpha } \right){\left\| {{x_{s + 1}} - {x_s}} \right\|^2}} \right]} } \right]\\
 \le \left( {2 + \alpha  - \gamma } \right)\mathop {\max }\limits_{u \in C} D\left( {u,{x_0}} \right)
\end{array}
\end{array}\]
The monotonicity of $F$ and the convexity of $g$ imply that
\[\alpha S{\mathbb{E}}\left[ {{\rm{Gap}}\left( {{z^S}} \right)} \right] \le {\mathbb{E}}\left[ {\mathop {\max }\limits_{u \in C} \left\{ {\alpha \sum\limits_{s = 0}^{S - 1} {\left[ {g\left( {{x_{s + 1}}} \right) - g\left( u \right) + \left\langle {F\left( {{x_{s + 1}}} \right),{x_{s + 1}} - u} \right\rangle } \right]} } \right\}} \right].\]
Hence
\[{\mathbb{E}}\left[ {{\rm{Gap}}\left( {{z^S}} \right)} \right] \le \frac{\left(2+\alpha-\gamma\right)}{{\alpha S}}\mathop {\max }\limits_{u \in C} D\left( {u,{x_0}} \right).\]
We compete this proof.
\end{proof}
\begin{corollary}
Under the conditions of Theorem \ref{T3.1},  choose $\alpha  = \min \left\{ {\frac{{\gamma  - p}}{{2\left( {1 - p} \right)}},\frac{{\left( {1 - \gamma} \right)b}}{{\left( {1 -p } \right)\left( {2{{\bar L}^2} + b{L^2}} \right)}}} \right\}$, $p=\frac{1}{M}$, $\gamma=p+\frac{1}{\sqrt{M}}$. Thus ${\mathbb{E}}\left[ {{\rm{Gap}}\left( {{z^S}} \right)} \right] \le {\cal O}\left( {\frac{{ \sqrt{M}}}{S}} \right)$ and the total average
complexity of Algorithm \ref{Alg1} to reach $\epsilon$-gap is $\mathcal{O}\left(\frac{{\sqrt{M}}}{\epsilon}\right)$, where ${z^S} = \frac{1}{S}\sum\limits_{s = 0}^{S - 1} {{x_{s + 1}}} $.
\end{corollary}
\begin{proof}
If $\alpha  = \min \left\{ {\frac{{\gamma  - p}}{{2\left( {1 - p} \right)}},\frac{{\left( {1 - \gamma} \right)b}}{{\left( {1 - p} \right)\left( {2{{\bar L}^2} + b{L^2}} \right)}}} \right\}$, $p=\frac{1}{M}$, $\gamma=p+\frac{1}{\sqrt{M}}$ from Theorem \ref{T3.1}, we have $\frac{1}{\alpha } = \max \left\{ {\frac{{2\left( {M - 1} \right)}}{{\sqrt M }},\frac{{\left( {M - 1} \right)}}{{\left( {M - 1 - \sqrt M } \right)}}\left( {\frac{2{\bar L}^2}{b} + {L^2}} \right)} \right\} ={\cal O}\left( {\sqrt M } \right)$, then
\[\begin{array}{l}
\mathbb{E}\left[ {{\rm{Gap}}\left( {{z^S}} \right)} \right] \le \frac{{\left( {2 + \alpha  - \gamma } \right)}}{{\alpha S}}\mathop {\max }\limits_{u \in C} D\left( {u,{x_0}} \right) \le \frac{{\left( {2 + \frac{{\gamma  - p}}{{2\left( {1 - p} \right)}} - \gamma } \right)}}{{\alpha S}}\mathop {\max }\limits_{u \in C} D\left( {u,{x_0}} \right) = \frac{{\left( {2 - \frac{{p\left( {1 - \gamma } \right)}}{{2\left( {1 - p} \right)}} - \frac{\gamma }{2}} \right)}}{{\alpha S}}\mathop {\max }\limits_{u \in C} D\left( {u,{x_0}} \right)\\
 \le \frac{2}{{\alpha S}}\mathop {\max }\limits_{u \in C} D\left( {u,{x_0}} \right) = O\left( {\frac{{\sqrt M }}{S}} \right).
\end{array}\]
 To get ${\mathbb{E}}\left[ {{\rm{Gap}}\left( {{z^S}} \right)} \right] \le \epsilon$,  the total average complexity in $S$ iterations is $M+2+\left( pM+2 \right)S=\mathcal{O}\left(\frac{{\sqrt M }}{\epsilon}\right)$. The proof is completed.
\end{proof}

\par Next, we  present a convergence rate analysis of Algorithm \ref{Alg1} for the case where the operator
$F$ is non-monotone. To support this extension, we introduce essential assumptions and lemmas that are indispensable for the subsequent analysis.
\begin{definition}
\cite{AZ}
The vector $\bar x \in K$ is an $\epsilon$-approximated solution of FSVI \eqref{1.1} if $$dist\left( {\mathbf{0},F\left( {\bar x} \right) + \partial \left( {g + {{\rm I}_K}} \right)\left( {\bar x} \right)} \right) =\mathop {\min }\limits_{u \in \partial \left( {g + {{\rm I}_K}} \right)\left( {\bar x} \right)} \left\| {F(\bar{x}) + u} \right\|_*\le \epsilon.$$
\end{definition}
 $\mathbf{Assumption}$ (A5)  There exists $x^* \in K$ and some $\rho >0$ such that
\[ - \rho {\left\| w \right\|^2} \le \left\langle {w,x - x^*} \right\rangle ,\forall x \in K,\forall w \in F\left( x \right) + \partial \left( {g + {{\rm I}_K}} \right)\left( x \right).\]

\begin{remark}
Assumption (A5) is called the weak minty variational inequalities which has been employed in the works of \cite{AZ,PT1,PT2}, serving as a crucial element in their respective analyses. Note that $\bar x \in K$ satisfies  $\mathbf{0} \in F\left( {\bar x} \right) + \partial \left( {g + {{\rm I}_K}} \right)\left( {\bar x} \right)$ if and only if $\bar x $ is a solution of FSVI \eqref{1.1}.
According to Definition \ref{def3.4}, it is easy to see that if the mapping $F+\partial \left( {g + {{\rm{I}}_K}} \right)$ is $\rho$-cohypomonotone, then Assumption (A5) holds.
\end{remark}
$\mathbf{Assumption}$ (A6) The gradient of the Bregman distance generating function $f$ is siad to be Lipschitz continuous if  for any
$x,y \in K$ we have $\left\| {\nabla f\left( x \right) - \nabla f\left( y \right)} \right\|_* \le {L_f}\left\| {x - y} \right\|$ for some $L_f>0$.
\begin{lemma}
\label{l2.3}
Assume that $\mathcal{F}=({\mathcal{F}}_s)_{s \ge 0}$ is a filtration and $u_s$ is a stochastic process adapted to $\mathcal{F}_s$ satisfying ${\mathbb{E}}[u_{s+1}|{\mathcal{F}}_s]=0$. Then
\[{\mathbb{E}}\left[ { \sum\limits_{s = 0}^{S - 1} {\left\langle {{u_{s + 1}},x} \right\rangle } } \right] \le \frac{1}{2}{\left\| {{x_0} - x} \right\|^2} + \frac{1}{2}\sum\limits_{s = 0}^{S - 1} {{\mathbb{E}}\left[ {{{\left\| {{u_{s + 1}}} \right\|}^2}} \right]}, \forall~ S \in \mathbb{N}, ~\forall~x,x_0 \in K. \]
\end{lemma}
\begin{proof}
The proof of this result is analogous to that of Lemma 2.4 of \cite{AA}. For the convenience of readers, we provide its detailed proof herein.
Set $x_{s+1}=x_s+u_{s+1}$. It is easy to see that $x_s$ is ${\mathcal{F}}_s$-measurable and
\[{\left\| {{x_{s + 1}} - x} \right\|^2} = {\left\| {{x_s} - x} \right\|^2} + 2\left\langle {{u_{s + 1}},{x_s} - x} \right\rangle  + {\left\| {{u_{s + 1}}} \right\|^2}.\]
Summing over $s=0,...,S-1$, we obtain
\[\sum\limits_{s = 0}^{S - 1} {2\left\langle {{u_{s + 1}},x - {x_s}} \right\rangle }  \le {\left\| {{x_0} - x} \right\|^2} + \sum\limits_{s = 0}^{S - 1} {{{\left\| {{u_{s + 1}}} \right\|}^2}} .\]
Taking expectation on both sides of the above inequality, we have
\[{\mathbb{E}}\left[ {\sum\limits_{s = 0}^{S - 1} {\left\langle {{u_{s + 1}},x} \right\rangle } } \right] \le \frac{1}{2}{\left\| {{x_0} - x} \right\|^2} + \frac{1}{2}\sum\limits_{s = 0}^{S - 1} {{\mathbb{E}}\left[ {{{\left\| {{u_{s + 1}}} \right\|}^2}} \right]}  + \sum\limits_{s = 0}^{S - 1} {{\mathbb{E}}\left[ {\left\langle {{u_{s + 1}},{x_s}} \right\rangle } \right]} .\]
Owing to ${\mathcal{F}}_s$-measurability of $x_s$ and ${\mathbb{E}}[u_{s+1}|{\mathcal{F}}_s]=0$, we can conclude that
\[\sum\limits_{s = 0}^{S - 1} {{\mathbb{E}}\left[ {\left\langle {{u_{s + 1}},{x_s}} \right\rangle } \right]}  = \sum\limits_{s = 0}^{S - 1} {{\mathbb{E}}\left[ {\left\langle {{\mathbb{E}}\left[ {{u_{s + 1}}|{{\mathcal{F}}_s}} \right],{x_s}} \right\rangle } \right]}  = 0.\]
 We have finished this proof.
\end{proof}
For the sake of convenience, define ${\sigma _1} = \frac{{1 - \gamma }}{2\left({1 - p}\right)}\left( {1 - \frac{{8\rho L_f^2}}{\alpha }} \right) - \left( {2{\alpha ^2}{L^2} + 8\alpha \rho {L^2} + \frac{{4\alpha \rho {{\bar L}^2}}}{b} + \frac{{5{\alpha ^2}{{\bar L}^2}}}{2b}} \right)$ and ${\sigma _2} = 8{L^2} + \frac{{4\left( {1 - \gamma } \right)L_f^2}}{{\left( {1 - p} \right){\alpha ^2}}} + \frac{{4{{\bar L}^2}}}{b}$.
\begin{theorem}
\label{theorem2}
  Let the sequence $\left\{ {{x_s}} \right\}$ be generated by Algorithm \ref{Alg1}, and Assumptions (A1), (A3), (A4), (A5) and (A6) hold, $0<p <1$ and $p<\gamma <1$. If ${\sigma _1}>0$ and
  $\frac{{8\left( {\gamma  - p} \right)L_f^2}}{{\left( {1 - p} \right){\alpha ^2}}} - \frac{{{\sigma _2}}}{{{\sigma _1}}}\left( {\left( {\frac{{\gamma  - p}}{{1 - p}}} \right)\left( {1 - \frac{{8\rho L_f^2}}{\alpha }} \right) - \frac{1}{2}} \right) \le 0$, then
  \[
\mathbb{E}\left[ {\sum\limits_{s = 0}^{S - 1} {{{\left\| {{u_s}} \right\|}_*^2}} } \right] \le \frac{{\left(3-\gamma \right){\sigma _2}}}{{{\sigma _1}}}D\left( {{x^*},{x_0}} \right),
\]
   where  ${u_s} = F\left( {{x_{s + 1}}} \right) - {\Delta ^s} - \frac{1}{\alpha }\left( {\nabla f\left( {{x_{s + 1}}} \right) - {\left( {\gamma \nabla f\left( {{x_s}} \right) + \left( {1 - \gamma } \right)\nabla f\left( {{x_{s-1 }}} \right)} \right)}} \right)$.
\end{theorem}
\begin{proof}
Via the optimality condition and  the definition of ${\bar x}_s$ and $x_{s+1}$ we know that
\[\mathbf{0} \in \frac{1}{\alpha }\left( {\nabla f\left( {{x_{s + 1}}} \right)  -{\left( {\gamma \nabla f\left( {{x_s}} \right) + \left( {1 - \gamma } \right)\nabla f\left( {{x_{s-1}}} \right)} \right)}} \right) + {\Delta ^s} + \partial \left( {g + {{\rm I}_K}} \right)\left( {{x_{s + 1}}} \right).\]
 Thus
\[{u_s} \in F\left( {{x_{s + 1}}} \right) + \partial \left( {g + {{\rm I}_K}} \right)\left( {{x_{s + 1}}} \right).\]
In view of  Assumption (A5) and \eqref{l2.1}, there exists $x^* \in K$ and $\rho >0$ such that
\[\begin{aligned}
 - \rho \left\| {{u_s}} \right\|_*^2 &\le \left\langle {F\left( {{x_{s + 1}}} \right) - {\Delta ^s},{x_{s + 1}} - {x^*}} \right\rangle  - \frac{1}{\alpha }\left\langle {\nabla f\left( {{x_{s + 1}}} \right) - \gamma \nabla f\left( {{x_s}} \right) - \left( {1 - \gamma } \right)\nabla f\left( {{x_{s-1}}} \right),{x_{s + 1}} - {x^*}} \right\rangle \\
 &=\left\langle {F\left( {{x_{s + 1}}} \right) - {\Delta ^s},{x_{s + 1}} - {x^*}} \right\rangle  + \frac{\gamma }{\alpha }\left( {D\left( {{x^*},{x_s}} \right) - D\left( {{x^*},{x_{s + 1}}} \right) - D\left( {{x_{s+1}},{x_{s  }}} \right)} \right)\\
 &+ \frac{{1 - \gamma }}{\alpha }\left( {D\left( {{x^*},{x_{s-1}}} \right) - D\left( {{x^*},{x_{s + 1}}} \right) - D\left( {{x_{s + 1}},{x_{s-1}}} \right)} \right).
\end{aligned}\]
Multiplying both sides of the above inequality by $\alpha$, it gets
\[ \begin{aligned}
 - \alpha \rho {\left\| {{u_s}} \right\|_*^2} &\le \alpha \left\langle {F\left( {{x_{s + 1}}} \right) - {\Delta ^s},{x_{s + 1}} - {x^*}} \right\rangle+
\gamma \left( {D\left( {{x^*},{x_s}} \right) - D\left( {{x^*},{x_{s + 1}}} \right) - D\left( {{x_{s+1}},{x_{s }}} \right)} \right)\\
 &+ \left(1 - \gamma\right) \left( {D\left( {{x^*},{x_{s-1}}} \right) - D\left( {{x^*},{x_{s + 1}}} \right) - D\left( {{x_{s + 1}},{x_{s-1}}} \right)} \right).
\end{aligned}
\]
By \eqref{l3.1.e1}, we have
\begin{equation}\label{theorem2.1}
  \begin{aligned}
 - \alpha \rho {\left\| {{u_s}} \right\|_*^2} &\le \alpha \left\langle {F\left( {{x_{s + 1}}} \right),{x_{s + 1}} - {x^*}} \right\rangle  + \alpha \left\langle {{{\mathbb{E}}_s}\left[ {{\Delta ^s}} \right] - {\Delta ^s},{x_{s + 1}} - {x_s}} \right\rangle\\
 & + \alpha \left\langle {{{\mathbb{E}}_s}\left[ {{\Delta ^s}} \right] - {\Delta ^s},{x_s} - {x^*}} \right\rangle  - \alpha \left\langle {{{\mathbb{E}}_s}\left[ {{\Delta ^s}} \right],{x_{s + 1}} - {x^*}} \right\rangle \\
 &+\gamma \left( {D\left( {{x^*},{x_s}} \right) - D\left( {{x^*},{x_{s + 1}}} \right) - D\left( {{x_{s+1}},{x_{s}}} \right)} \right)
 \\&+ \left( {1 - \gamma } \right)\left( {D\left( {{x^*},{x_{s-1}}} \right) - D\left( {{x^*},{x_{s + 1}}} \right) - D\left( {{x_{s + 1}},{x_{s-1}}} \right)} \right)\\
 &= \alpha \left\langle {F\left( {{x_{s + 1}}} \right),{x_{s + 1}} - {x^*}} \right\rangle+ \alpha \left\langle {{{\mathbb{E}}_s}\left[ {{\Delta ^s}} \right] - {\Delta ^s},{x_{s + 1}} - {x_s}} \right\rangle  + \alpha \left\langle {{{\mathbb{E}}_s}\left[ {{\Delta ^s}} \right] - {\Delta ^s},{x_s} - {x^*}} \right\rangle   \\
 &+\gamma \left( {D\left( {{x^*},{x_s}} \right) - D\left( {{x^*},{x_{s + 1}}} \right) - D\left( {{x_{s+1}},{x_{s}}} \right)} \right)
 \\&+ \left( {1 - \gamma } \right)\left( {D\left( {{x^*},{x_{s-1}}} \right) - D\left( {{x^*},{x_{s + 1}}} \right) - D\left( {{x_{s + 1}},{x_{s-1}}} \right)} \right)\\
 &- \alpha \left\langle {F\left( {{w_s}} \right) + F\left( {{x_s}} \right) - F\left( {{w_{s - 1}}} \right),{x_{s + 1}} - {x^*}} \right\rangle \\
 & =- \alpha \left\langle {F\left( {{w_s}} \right) - F\left( {{x_{s + 1}}} \right),{x_{s + 1}} - {x^*}} \right\rangle  - \alpha \left\langle {F\left( {{x_s}} \right) - F\left( {{w_{s - 1}}} \right),{x_s} - {x^*}} \right\rangle \\
 &- \alpha \left\langle {F\left( {{x_s}} \right) - F\left( {{w_{s - 1}}} \right),{x_{s + 1}} - {x_s}} \right\rangle  + \alpha \left\langle {{{\mathbb{E}}_s}\left[ {{\Delta ^s}} \right] - {\Delta ^s},{x_{s + 1}} - {x_s}} \right\rangle \\
 &+ \alpha \left\langle {{{\mathbb{E}}_s}\left[ {{\Delta ^s}} \right] - {\Delta ^s},{x_s} - {x^*}} \right\rangle  +\gamma \left( {D\left( {{x^*},{x_s}} \right) - D\left( {{x^*},{x_{s + 1}}} \right) - D\left( {{x_{s+1}},{x_{s }}} \right)} \right)
 \\&+ \left( {1 - \gamma } \right)\left( {D\left( {{x^*},{x_{s - 1}}} \right) - D\left( {{x^*},{x_{s + 1}}} \right) - D\left( {{x_{s + 1}},{x_{s - 1}}} \right)} \right).
\end{aligned}
\end{equation}
Using Assumptions (A4) and (A6), (\ref{l3.1.e1}) and the convexity of ${\left\|  \cdot  \right\|}^2$, we can get
\begin{equation}\label{theorem2.2}
\begin{aligned}
{\left\| {{u_s}} \right\|_*^2} &= {\left\| {F\left( {{x_{s + 1}}} \right) - {{\mathbb{E}}_s}\left[ {{\Delta ^s}} \right] + {{\mathbb{E}}_s}\left[ {{\Delta ^s}} \right] - {\Delta ^s} - \frac{1}{\alpha }\left( {\nabla f\left( {{x_{s + 1}}} \right) - {\left( {\gamma \nabla f\left( {{x_s}} \right) + \left( {1 - \gamma } \right)\nabla f\left( {{x_{s - 1}}} \right)} \right)}} \right)} \right\|_*^2}\\
 &\le 2{\left\| {F\left( {{x_{s + 1}}} \right) - {{\mathbb{E}}_s}\left[ {{\Delta ^s}} \right]} \right\|_*^2} + 4{\left\| {{{\mathbb{E}}_s}\left[ {{\Delta ^s}} \right] - {\Delta ^s}} \right\|_*^2} \\
 &+ \frac{4}{{{\alpha ^2}}}{\left\| {\nabla f\left( {{x_{s + 1}}} \right) - {\left( {\gamma \nabla f\left( {{x_s}} \right) + \left( {1 - \gamma } \right)\nabla f\left( {{x_{s - 1}}} \right)} \right)}} \right\|_*^2}\\
 &\le 4{\left\| {F\left( {{x_{s + 1}}} \right) - F\left( {{w_s}} \right)} \right\|_*^2} + 4{\left\| {F\left( {{x_s}} \right) - F\left( {{w_{s - 1}}} \right)} \right\|_*^2} + 4{\left\| {{{\mathbb{E}}_s}\left[ {{\Delta ^s}} \right] - {\Delta ^s}} \right\|_*^2}\\
 &+ \frac{{4\gamma }}{{{\alpha ^2}}}\left\| {\nabla f\left( {{x_{s + 1}}} \right) - \nabla f\left( {{x_s}} \right)} \right\|_*^2 + \frac{{4\left( {1 - \gamma } \right)}}{{{\alpha ^2}}}\left\| {\nabla f\left( {{x_{s + 1}}} \right) - \nabla f\left( {{x_{s - 1}}} \right)} \right\|_*^2\\
 &\le 4{L^2}{\left\| {{x_{s + 1}} - {w_s}} \right\|^2} + 4{L^2}{\left\| {{x_s} - {w_{s - 1}}} \right\|^2}+ \frac{{4\gamma L_f^2}}{{{\alpha ^2}}}{\left\| {{x_{s + 1}} - {x_s}} \right\|^2}  \\
 &+ \frac{{4\left( {1 - \gamma } \right)L_f^2}}{{{\alpha ^2}}}{\left\| {{x_{s + 1}} - {x_{s-1}}} \right\|^2}+ 4{\left\| {{{\mathbb{E}}_s}\left[ {{\Delta ^s}} \right] - {\Delta ^s}} \right\|_*^2}.
\end{aligned}
\end{equation}
By Young's inequality and Assumption (A4), we have
\begin{equation}\label{et1.2}
\begin{aligned}
 - \alpha \left\langle {F\left( {{x_s}} \right) - F\left( {{w_{s - 1}}} \right),{x_{s + 1}} - {x_s}} \right\rangle  &\le 2 \alpha^2{\left\| {F\left( {{x_s}} \right) - F\left( {{w_{s - 1}}} \right)} \right\|_*^2} + \frac{{1}}{8}{\left\| {{x_{s + 1}} - {x_s}} \right\|^2}\\
 &\le 2\alpha^2 L^2{\left\| {{x_{s}} - {w_{s-1}}} \right\|^2} + \frac{{1}}{8}{\left\| {{x_{s + 1}} - {x_s}} \right\|^2}.
\end{aligned}
\end{equation}
Substituting \eqref{et1.2}) and \eqref{theorem2.2} into \eqref{theorem2.1},  from \eqref{2.1}, it can conclude that
\begin{equation}
\label{+1}
\begin{aligned}
0 &\le \left( {{2\alpha^2 }{L^2} + 4\alpha \rho {{ L}^2}} \right){\left\| {{x_s} - {w_{s - 1}}} \right\|^2} +\left({ \frac{1}{4} + \frac{{8\gamma \rho L_f^2}}{\alpha }} \right)D\left(x_{s+1},x_s\right)+ 4\alpha \rho {{ L}^2}{\left\| {{x_{s + 1}} - {w_s}} \right\|^2}\\
 &{ + \frac{{8 \left(1-\gamma \right) \rho L_f^2}}{{{\alpha }}}} D\left(x_{s+1},x_{s-1}\right)- \alpha \left\langle {F\left( {{w_s}} \right) - F\left( {{x_{s + 1}}} \right),{x_{s + 1}} - {x^*}} \right\rangle  - \alpha \left\langle {F\left( {{x_s}} \right) - F\left( {{w_{s - 1}}} \right),{x_s} - {x^*}} \right\rangle \\
 &+ \alpha \left\langle {{{\mathbb{E}}_s}\left[ {{\Delta ^s}} \right] - {\Delta ^s},{x_{s + 1}} - {x_s}} \right\rangle  + 4\alpha \rho {\left\| {{{\mathbb{E}}_s}\left[ {{\Delta ^s}} \right] - {\Delta ^s}} \right\|_*^2}\\
 &+ \alpha \left\langle {{{\mathbb{E}}_s}\left[ {{\Delta ^s}} \right] - {\Delta ^s},{x_s} - {x^*}} \right\rangle + \gamma \left( {D\left( {{x^*},{x_s}} \right) - D\left( {{x^*},{x_{s + 1}}} \right) - D\left( {{x_{s+1}},{x_{s }}} \right)} \right)
 \\&+ \left( {1 - \gamma } \right)\left( {D\left( {{x^*},{x_{s - 1}}} \right) - D\left( {{x^*},{x_{s + 1}}} \right) - D\left( {{x_{s + 1}},{x_{s - 1}}} \right)} \right).
\end{aligned}
\end{equation}

 The Young's inequality implies that
 \begin{equation}\label{ext1.6}
\alpha {\mathbb{E}}\left[ {\sum\limits_{s = 0}^{S - 1} {\left[ {\left\langle {{{\mathbb{E}}_s}\left[ {{\Delta ^s}} \right] - {\Delta ^s},{x_{s + 1}} - {x_s}} \right\rangle } \right]} } \right] \le 2\alpha^2\sum\limits_{s = 0}^{S - 1} {{\mathbb{E}}\left[ {{{\left\| {{{\mathbb{E}}_s}\left[ {{\Delta ^s}} \right] - {\Delta ^s}} \right\|}_*^2}} \right]}  + \frac{{1}}{8}\sum\limits_{s = 0}^{S - 1} {{\mathbb{E}}\left[ {{{\left\| {{x_{s + 1}} - {x_s}} \right\|}^2}} \right]},
 \end{equation}
 and
 \begin{equation}\label{ext1.7}
\begin{aligned}
\alpha \left\langle {F\left( {{x_S}} \right) - F\left( {{w_{S - 1}}} \right),{x_S} - u} \right\rangle  &\le \frac{{{\alpha^2 }}}{2}{\left\| {F\left( {{x_S}} \right) - F\left( {{w_{S - 1}}} \right)} \right\|_*^2} + \frac{1}{2}{\left\| {{x_S} - u} \right\|^2}\\
 &\le  \frac{{{\alpha^2 }{L^2}}}{2}{\left\| {{x_S} - w_{S-1}} \right\|^2}+ \frac{1}{2}{\left\| {{x_S} - u} \right\|^2},
\end{aligned}
 \end{equation}
  where the last inequality is from Assumption (A4).\\
Summing over $s=0,...,S-1$, taking expectation on both sides of \eqref{+1}, and using \eqref{2.1}, \eqref{extral1}, \eqref{ext1.6} and \eqref{ext1.7}, we have
\begin{equation}\label{theorem2.3}
\begin{aligned}
0 &\le {\mathbb{E}}\left[ {\sum\limits_{s = 0}^{S - 1} \left( {{2\alpha^2 }{L^2} + 4\alpha \rho {{ L}^2}} \right){\left\| {{x_s} - {w_{s - 1}}} \right\|^2} +\left({ \frac{1}{2} + \frac{{8\gamma \rho L_f^2}}{\alpha }} \right)D\left(x_{s+1},x_s\right)+ 4\alpha \rho {{ L}^2}{\left\| {{x_{s + 1}} - {w_s}} \right\|^2} } \right]\\
 &+ {\mathbb{E}}\left[ {\sum\limits_{s = 0}^{S - 1} {\left[ {\left( {4\alpha \rho  + 2{\alpha ^2}} \right){{\left\| {{{\mathbb{E}}_s}\left[ {{\Delta ^s}} \right] - {\Delta ^s}} \right\|}_*^2} + \alpha \left\langle {{{\mathbb{E}}_s}\left[ {{\Delta ^s}} \right] - {\Delta ^s},{x_s} - {x^*}} \right\rangle } \right]{ + \frac{{8 \left(1-\gamma \right) \rho L_f^2}}{{{\alpha }}}} D\left(x_{s+1},x_{s-1}\right)} } \right]\\
 &+ {\mathbb{E}}\left[ {\sum\limits_{s = 0}^{S - 1} {\left[ \gamma \left( {D\left( {{x^*},{x_s}} \right) - D\left( {{x^*},{x_{s + 1}}} \right) - D\left( {{x_{s+1}},{x_{s }}} \right)} \right)  \right] + \frac{{\alpha^2 {L^2}}}{2}{\left\| {{x_S} - {w_{S - 1}}} \right\|^2} + \frac{1}{2}{\left\| {{x_S} - {x^*}} \right\|^2} } } \right]
 \\
 &+ {\mathbb{E}}\left[ {\sum\limits_{s = 0}^{S - 1} {\left[  \left( {1 - \gamma } \right)\left( {D\left( {{x^*},{x_{s - 1}}} \right) - D\left( {{x^*},{x_{s + 1}}} \right) - D\left( {{x_{s + 1}},{x_{s - 1}}} \right)} \right) \right]} } \right].
\end{aligned}
\end{equation}
It follows from \eqref{var1} and Lemma \ref{l2.3} that
\begin{equation}\label{theorem2.4}
\begin{aligned}
{\mathbb{E}}\left[ {\sum\limits_{s = 0}^{S - 1} {\left[ {\alpha \left\langle {{{\mathbb{E}}_s}\left[ {{\Delta ^s}} \right] - {\Delta ^s},{x_s} - {x^*}} \right\rangle } \right]} } \right] &= {\mathbb{E}}\left[ {\sum\limits_{s = 0}^{S - 1} {\left[ {\alpha \left\langle { - \left( {{{\mathbb{E}}_s}\left[ {{\Delta ^s}} \right] - {\Delta ^s}} \right),{x^*}} \right\rangle } \right]} } \right]\\
 &\le \frac{1}{2}{\left\| {{x^*} - {x_0}} \right\|^2} + \frac{{{\alpha ^2}}}{2}\sum\limits_{s = 0}^{S - 1} {{\mathbb{E}}\left[ {{{\left\| {{{\mathbb{E}}_s}\left[ {{\Delta ^s}} \right] - {\Delta ^s}} \right\|}_*^2}} \right]}.
\end{aligned}
\end{equation}
Thanks to $x_0=x_{-1}=w_{-1}$, we have
\begin{equation}\label{++1}
\begin{array}{l}
\mathbb{E}\left[ {\sum\limits_{s = 0}^{S - 1} {\gamma \left( {D\left( {{x^*},{x_s}} \right) - D\left( {{x^*},{x_{s + 1}}} \right)} \right) + \left( {1 - \gamma } \right)\left( {D\left( {{x^*},{x_{s - 1}}} \right) - D\left( {{x^*},{x_{s + 1}}} \right)} \right)} } \right]\\
 = \mathbb{E}\left[ {\sum\limits_{s = 0}^{S - 1} {D\left( {{x^*},{x_s}} \right) - D\left( {{x^*},{x_{s + 1}}} \right) + \left( {1 - \gamma } \right)\left( {D\left( {{x^*},{x_{s - 1}}} \right) - D\left( {{x^*},{x_s}} \right)} \right)} } \right]\\
 = \left( {2 - \gamma } \right)D\left( {{x^*},{x_0}} \right) + \mathbb{E}\left[ { - D\left( {{x^*},{x_S}} \right) - \left( {1 - \gamma } \right)D\left( {{x^*},{x_{S - 1}}} \right)} \right],
\end{array}
\end{equation}
and
\begin{equation}\label{+++1}
\mathbb{E}\left[ {\sum\limits_{s = 0}^{S - 1} {\left[ {{{\left\| {{x_{s + 1}} - {w_s}} \right\|}^2}} \right]} } \right] = \mathbb{E}\left[ {\sum\limits_{s = 0}^{S - 1} {\left[ {{{\left\| {{x_s} - {w_{s - 1}}} \right\|}^2}} \right] + {{\left\| {{x_S} - {w_{S - 1}}} \right\|}^2}} } \right].
\end{equation}
Substituting \eqref{theorem2.4}, \eqref{++1} and \eqref{+++1} into \eqref{theorem2.3}, by Lemma \ref{l3.1}, \eqref{2.1} and $0<\gamma<1$, we can deduce that
\begin{equation}\label{3.32}
\begin{aligned}
0 &\le {\mathbb{E}}\left[ {\sum\limits_{s = 0}^{S - 1} \left( {{2\alpha^2 }{L^2} + 4\alpha \rho {{ L}^2}} \right){\left\| {{x_s} - {w_{s - 1}}} \right\|^2} +\left({ \frac{1}{2} + \frac{{8\gamma \rho L_f^2}}{\alpha }} \right)D\left(x_{s+1},x_s\right)+ 4\alpha \rho {{ L}^2}{\left\| {{x_{s + 1}} - {w_s}} \right\|^2} } \right]\\
 &+ {\mathbb{E}}\left[ {\sum\limits_{s = 0}^{S - 1} {\left[ {\left( {4\alpha \rho  + \frac{5}{2}{\alpha ^2}} \right){{\left\| {{{\mathbb{E}}_s}\left[ {{\Delta ^s}} \right] - {\Delta ^s}} \right\|}_*^2} } \right]} }{ + \frac{{8 \left(1-\gamma \right) \rho L_f^2}}{{{\alpha }}}} D\left(x_{s+1},x_{s-1}\right) \right]+ \frac{1}{2}{{\left\| {{x^*} - {x_0}} \right\|}^2}\\
 &+ {\mathbb{E}}\left[ {\sum\limits_{s = 0}^{S - 1} {\left[ \gamma \left( {D\left( {{x^*},{x_s}} \right) - D\left( {{x^*},{x_{s + 1}}} \right) - D\left( {{x_{s+1}},{x_{s }}} \right)} \right)  \right]+ \frac{{{\alpha ^2}{L^2}}}{2}{\left\| {{x_S} - {w_{S - 1}}} \right\|^2} + D\left(x^*,x_{S}\right) } } \right]
 \\
 &+ {\mathbb{E}}\left[ {\sum\limits_{s = 0}^{S - 1} {\left[  \left( {1 - \gamma } \right)\left( {D\left( {{x^*},{x_{s - 1}}} \right) - D\left( {{x^*},{x_{s + 1}}} \right) - D\left( {{x_{s + 1}},{x_{s - 1}}} \right)} \right) \right]} } \right]\\
 &\le {\mathbb{E}}\left[ {\sum\limits_{s = 0}^{S - 1} {\left[ {\left( {2{\alpha^2}{L^2} + 4\alpha \rho {{ L}^2}+\frac{4\alpha \rho {{\bar L}^2}}{b} + \frac{5{\alpha ^2}{{\bar L}^2}}{2b}} \right){\left\| {{x_s} - {w_{s - 1}}} \right\|^2}+4\alpha \rho {{ L}^2}{\left\| {{x_{s+1}} - {w_{s}}} \right\|^2} } \right] } } \right]\\
 &+\mathbb{E}\left[ {\sum\limits_{s = 0}^{S - 1} {\left( { - \gamma + \frac{1}{2} + \frac{{8\gamma \rho L_f^2}}{\alpha }} \right)D\left( {{x_{s + 1}},{x_s}} \right)+{{\left(1-\gamma\right)}\left(\frac{{8\rho L_f^2}}{\alpha } - 1\right)D\left( {{x_{s + 1}},{x_{s-1}}} \right)}} } \right]\\
 &+\left(3-\gamma \right)D\left( {{x^*},{x_0}} \right)+ {\mathbb{E}}\left[ {  -\left(1-\gamma\right)D\left( {{x^*},{x_{S-1}}} \right) +\frac{{{\alpha ^2}{L^2}}}{2}{\left\| {{x_S} - {w_{S - 1}}} \right\|^2} } \right]\\
 &\le {\mathbb{E}}\left[ {\sum\limits_{s = 0}^{S - 1} {\left[ {\left( {2{\alpha^2}{L^2} + 8\alpha \rho {{ L}^2}+\frac{4\alpha \rho {{\bar L}^2}}{b} + \frac{5{\alpha ^2}{{\bar L}^2}}{2b}} \right){\left\| {{x_{s + 1}} - {w_s}} \right\|^2} } \right] } } \right]\\
 &+\mathbb{E}\left[ {\sum\limits_{s = 0}^{S - 1} {\left( { - \gamma + \frac{1}{2} + \frac{{8\gamma \rho L_f^2}}{\alpha }} \right)D\left( {{x_{s + 1}},{x_s}} \right)+{{\left(1-\gamma\right)}\left(\frac{{8\rho L_f^2}}{\alpha } - 1\right)D\left( {{x_{s + 1}},{x_{s-1}}} \right)}} } \right]\\
 &+\left(3-\gamma \right)D\left( {{x^*},{x_0}} \right).
 \end{aligned}
\end{equation}
According to $\sigma_1 >0$,  we have ${1 - \frac{{8\rho L_f^2}}{\alpha }}>0$. Therefore, using \eqref{2.1}, \eqref{t1.9} and \eqref{3.32}, we can get
\[
\begin{aligned}
0&\le {\mathbb{E}}\left[ {\sum\limits_{s = 0}^{S - 1} {\left[ {\left( {2{\alpha^2}{L^2} + 8\alpha \rho {{ L}^2}+\frac{4\alpha \rho {{\bar L}^2}}{b} + \frac{5{\alpha ^2}{{\bar L}^2}}{2b}} \right){\left\| {{x_{s + 1}} - {w_s}} \right\|^2} } \right] } } \right]\\
 &+\mathbb{E}\left[ {\sum\limits_{s = 0}^{S - 1} {\left( { - \gamma + \frac{1}{2} + \frac{{8\gamma \rho L_f^2}}{\alpha }} \right)D\left( {{x_{s + 1}},{x_s}} \right)+{{\left(1-\gamma\right)}\left(\frac{{8\rho L_f^2}}{\alpha } - 1\right)D\left( {{x_{s + 1}},{x_{s-1}}} \right)}} } \right]\\
  &+ \mathbb{E}\left[ {\sum\limits_{s = 0}^{S - 1} {\left[ {\frac{{1 - \gamma }}{{2\left( {1 - p} \right)}}\left( {1 - \frac{{8\rho L_f^2}}{\alpha }} \right){{\left\| {{x_{s + 1}} - {w_s}} \right\|}^2}} \right]} } \right] - \mathbb{E}\left[ {\sum\limits_{s = 0}^{S - 1} {\left[ {\frac{{1 - \gamma }}{{2\left( {1 - p} \right)}}\left( {1 - \frac{{8\rho L_f^2}}{\alpha }} \right){{\left\| {{x_{s + 1}} - {w_s}} \right\|}^2}} \right]} } \right] \\
 &+\left(3-\gamma \right)D\left( {{x^*},{x_0}} \right)\\
&={\mathbb{E}}\left[ {\sum\limits_{s = 0}^{S - 1} {\left[ {\left( {2{\alpha^2}{L^2} + 8\alpha \rho {{ L}^2}+\frac{4\alpha \rho {{\bar L}^2}}{b} + \frac{5{\alpha ^2}{{\bar L}^2}}{2b}-\frac{1-\gamma}{2\left(1-p\right)}\left(1-\frac{{8\rho L_f^2}}{\alpha } \right)} \right){\left\| {{x_{s + 1}} - {w_s}} \right\|^2} } \right] } } \right]\\
 &+\mathbb{E}\left[ {\sum\limits_{s = 0}^{S - 1} {\left( {\frac{1}{2}-\gamma \left( {1 - \frac{{8\rho L_f^2}}{\alpha }} \right)} \right)D\left( {{x_{s + 1}},{x_s}} \right)+\frac{{p\left( {1 - \gamma } \right)}}{{2\left( {1 - p} \right)}}\left( {1 - \frac{{8\rho L_f^2}}{\alpha }} \right){\left\| {{x_{s + 1}} - {x_s}} \right\|^2}} } \right]\\
 &+\mathbb{E}\left[ {\sum\limits_{s = 0}^{S - 1} { {\left( {1 - \gamma } \right)\left( {1 - \frac{{8\rho L_f^2}}{\alpha }} \right)} \left({\frac{1}{2}\left\| {{x_{s + 1}} - {x_{s-1}}} \right\|^2}-D\left( {{x_{s + 1}},{x_{s-1}}}\right) \right)} } \right]\\
 &+\left(3-\gamma \right)D\left( {{x^*},{x_0}} \right)\\
 &\le {\mathbb{E}}\left[ {\sum\limits_{s = 0}^{S - 1} {\left[ {\left( {2{\alpha^2}{L^2} + 8\alpha \rho {{ L}^2}+\frac{4\alpha \rho {{\bar L}^2}}{b} + \frac{5{\alpha ^2}{{\bar L}^2}}{2b}-\frac{1-\gamma}{2\left(1-p\right)}\left(1-\frac{{8\rho L_f^2}}{\alpha } \right)} \right){\left\| {{x_{s + 1}} - {w_s}} \right\|^2} } \right] } } \right]\\
 &+ \mathbb{E}\left[ {\sum\limits_{s = 0}^{S - 1} {\left( {\frac{1}{2}-\frac{\gamma-p}{1-p} \left( {1 - \frac{{8\rho L_f^2}}{\alpha }} \right)} \right)D\left( {{x_{s + 1}},{x_s}} \right)} } \right]\\
 &+\left(3-\gamma \right)D\left( {{x^*},{x_0}} \right)
  \end{aligned}
\]
Since $\sigma_1 >0$, we have
\begin{equation}\label{theorem2.5}
\begin{aligned}
{\mathbb{E}}\left[ {\sum\limits_{s = 0}^{S - 1} {\left\| {{x_{s + 1}} - {w_s}} \right\|^2} } \right] &\le \frac{\left(3-\gamma \right)}{{{\sigma _1}}}D\left( {{x^*},{x_0}} \right)\\
&-\mathbb{E}\left[ {\sum\limits_{s = 0}^{S - 1} {\frac{1}{{{\sigma _1}}}\left( {\left( {\frac{{\gamma  - p}}{{1 - p}}} \right)\left( {1 - \frac{{8\rho L_f^2}}{\alpha }} \right)-\frac{1}{2}} \right)D\left( {{x_{s + 1}},{x_{s}}} \right)} } \right].
  \end{aligned}
\end{equation}
For \eqref{theorem2.2}, summing over $s=0,...,S-1$ and taking expectation of both sides,  from Lemma \ref{l3.1}, \eqref{2.1} and \eqref{t1.9}, it follows that
\begin{equation}\label{3.34}
\begin{aligned}
\mathbb{E}\left[ \sum\limits_{s = 0}^{S - 1}\left[{\left\| {{u_s}} \right\|_*^2} \right]\right] &\le \mathbb{E}\left[\sum\limits_{s = 0}^{S - 1}\left[ {4{L^2}{{\left\| {{x_{s + 1}} - {w_s}} \right\|}^2} + 4{L^2}{{\left\| {{x_s} - {w_{s - 1}}} \right\|}^2} + 4\left\| {{\mathbb{E}_s}\left[ {{\Delta ^s}} \right] - {\Delta ^s}} \right\|_*^2} \right]\right]\\
 &+ \mathbb{E}\left[ \sum\limits_{s = 0}^{S - 1}\left[{\frac{{4\gamma L_f^2}}{{{\alpha ^2}}}{{\left\| {{x_{s + 1}} - {x_s}} \right\|}^2} + \frac{{4\left( {1 - \gamma } \right)L_f^2}}{{{\alpha ^2}}}{{\left\| {{x_{s + 1}} - {x_{s - 1}}} \right\|}^2}}\right] \right]\\
 &\le \mathbb{E}\left[\sum\limits_{s = 0}^{S - 1} \left[ {4{L^2}{{\left\| {{x_{s + 1}} - {w_s}} \right\|}^2} + \left( {4{L^2} + \frac{{4{{\bar L}^2}}}{b}} \right){{\left\| {{x_s} - {w_{s - 1}}} \right\|}^2}}\right] \right]\\
 &+ \mathbb{E}\left[\sum\limits_{s = 0}^{S - 1}\left[ {\frac{{4\gamma L_f^2}}{{{\alpha ^2}}}{{\left\| {{x_{s + 1}} - {x_s}} \right\|}^2} + \frac{{4\left( {1 - \gamma } \right)L_f^2}}{{{\alpha ^2}}}{{\left\| {{x_{s + 1}} - {x_{s - 1}}} \right\|}^2}}\right] \right]\\
 &= \mathbb{E}\left[\sum\limits_{s = 0}^{S - 1}\left[ {4{L^2}{{\left\| {{x_{s + 1}} - {w_s}} \right\|}^2} + \left( {4{L^2} + \frac{{4{{\bar L}^2}}}{b}} \right){{\left\| {{x_s} - {w_{s - 1}}} \right\|}^2}}\right] \right]\\
 &+ \mathbb{E}\left[\sum\limits_{s = 0}^{S - 1}\left[ {\frac{{4\left( {1 - \gamma } \right)L_f^2}}{{\left( {1 - p} \right){\alpha ^2}}}{{\left\| {{x_{s + 1}} - {w_s}} \right\|}^2}}\right] \right] - \mathbb{E}\left[\sum\limits_{s = 0}^{S - 1}\left[ {\frac{{4\left( {1 - \gamma } \right)L_f^2}}{{\left( {1 - p} \right){\alpha ^2}}}{{\left\| {{x_{s + 1}} - {w_s}} \right\|}^2}}\right] \right]\\
 &+ \mathbb{E}\left[ \sum\limits_{s = 0}^{S - 1}\left[{\frac{{4\gamma L_f^2}}{{{\alpha ^2}}}{{\left\| {{x_{s + 1}} - {x_s}} \right\|}^2} + \frac{{4\left( {1 - \gamma } \right)L_f^2}}{{{\alpha ^2}}}{{\left\| {{x_{s + 1}} - {x_{s - 1}}} \right\|}^2}}\right] \right]\\
 &= {\mathbb{E}}\left[ {\sum\limits_{s = 0}^{S - 1} {\left[ {\left(4L^2+\frac{{4\left( {1 - \gamma } \right)L_f^2}}{{\left( {1 - p} \right){\alpha ^2}}}\right){\left\| {{x_{s + 1}} - {w_s}} \right\|^2} +\left( {4L^2+\frac{4{\bar L}^2}{b}}\right){\left\| {{x_{s }} - {w_{s-1}}} \right\|^2}} \right]} } \right]\\
 &+{\mathbb{E}}\left[ {\sum\limits_{s = 0}^{S - 1} {\left[\frac{{4\left( { \gamma-p } \right)L_f^2}}{{\left( {1 - p} \right){\alpha ^2}}}{\left\| {{x_{s + 1}} - {x_s}} \right\|^2}\right]} } \right].
\end{aligned}
\end{equation}
Since $x_0=w_{-1}$, we have $\sum\limits_{s = 0}^{S - 1} {\left\| {{x_s} - {w_{s - 1}}} \right\|^2}  \le \sum\limits_{s = 0}^{S - 1} {\left\| {{x_{s + 1}} - {w_s}} \right\|^2} $. Therefore, using \eqref{2.1}, we can obtain
\begin{equation}\label{theorem2.6}
{\mathbb{E}}\left[ {\sum\limits_{s = 0}^{S - 1} {{{\left\| {{u_s}} \right\|}_*^2}} } \right]
 \le {\mathbb{E}}\left[ {\sum\limits_{s = 0}^{S - 1} {\left[ {\sigma_2 {\left\| {{x_{s + 1}} - {w_s}} \right\|^2}}+{ {  \frac{{8 \left(\gamma-p \right) L_f^2}}{{{\left(1-p\right)\alpha ^2}}}} D\left(x_{s+1},x_{s}\right) } \right]} } \right]
\end{equation}
Combining \eqref{theorem2.5} with \eqref{theorem2.6}, we have
\[
\mathbb{E}\left[ {\sum\limits_{s = 0}^{S - 1} {{{\left\| {{u_s}} \right\|}_*^2}} } \right] \le \frac{{\left(3-\gamma \right){\sigma _2}}}{{{\sigma _1}}}D\left( {{x^*},{x_0}} \right),
\]
where we use the fact that $\frac{{8\left( {\gamma  - p} \right)L_f^2}}{{\left( {1 - p} \right){\alpha ^2}}} - \frac{{{\sigma _2}}}{{{\sigma _1}}}\left( {\left( {\frac{{\gamma  - p}}{{1 - p}}} \right)\left( {1 - \frac{{8\rho L_f^2}}{\alpha }} \right) - \frac{1}{2}} \right) \le 0$.\\
 The proof is finished.
\end{proof}

It is interesting to observe that for $p=1$, that is $w_s=x_s$, we have the following result.
\begin{theorem}
\label{theor3.5}
  Let the sequence $\left\{ {{x_s}} \right\}$ be generated by Algorithm \ref{Alg1}, and Assumptions (A1), (A3), (A4), (A5) and (A6) hold. If $p=1$, ${\sigma' _1}>0$ and
  $\left( {1 - \gamma } \right)\left( {\frac{{{\sigma' _1}}}{{{\sigma' _2}}}\left( {\frac{{8\rho L_f^2}}{\alpha } - 1} \right) + \frac{{8L_f^2}}{{{\alpha ^2}}}} \right) \le 0$, then
  \[
\mathbb{E}\left[ {\sum\limits_{s = 0}^{S - 1} {{{\left\| {{u_s}} \right\|}_*^2}} } \right] \le \frac{{\left(3-\gamma \right){\sigma' _2}}}{{{\sigma' _1}}}D\left( {{x^*},{x_0}} \right),
\]
   where $\sigma '_1 =\left(\gamma-\frac{1}{2}\right)-\left( {4{\alpha ^2}{L^2} +16\alpha \rho {{ L}^2}+ \frac{8\alpha \rho {{\bar L}^2}}{b} + \frac{5{\alpha ^2}{{\bar L}^2}}{b} + \frac{{8\gamma \rho L_f^2}}{\alpha }}\right), {\sigma' _2} = 16{L^2} + \frac{{8\gamma L_f^2}}{{{\alpha ^2}}} + \frac{8{{\bar L}^2}}{b}$, ${u_s} = F\left( {{x_{s + 1}}} \right) - {\Delta ^s} - \frac{1}{\alpha }\left( {\nabla f\left( {{x_{s + 1}}} \right) - {\left( {\gamma \nabla f\left( {{x_s}} \right) + \left( {1 - \gamma } \right)\nabla f\left( {{x_{s - 1}}} \right)} \right)}} \right)$.
\end{theorem}
\begin{proof}
  Using the fact $w_s=x_s$, \eqref{2.1} and \eqref{3.32}, we have
\[
\begin{aligned}
 0&\le \left(3-\gamma \right)D\left( {{x^*},{x_0}} \right) +\mathbb{E}\left[ {\sum\limits_{s = 0}^{S - 1} {{\left(1-\gamma\right)}\left(\frac{{8\rho L_f^2}}{\alpha } - 1\right)D\left( {{x_{s + 1}},{x_{s-1}}} \right)} } \right]\\
 &+ {\mathbb{E}}\left[ {\sum\limits_{s = 0}^{S - 1} {\left[ {\left( {4{\alpha ^2}{L^2} +16\alpha \rho {{ L}^2}+ \frac{8\alpha \rho {{\bar L}^2}}{b} + \frac{5{\alpha ^2}{{\bar L}^2}}{b} - \left(\gamma-\frac{1}{2}\right) + \frac{{8\gamma \rho L_f^2}}{\alpha }} \right)D\left( {{x_{s + 1}},{x_s}} \right)} \right]} } \right].
 \end{aligned}
 \]
 Since $\sigma'_1 >0$, we have
\begin{equation}\label{theorem2.5'}
{\mathbb{E}}\left[ {\sum\limits_{s = 0}^{S - 1} {D\left( {{x_{s + 1}},{x_s}} \right)} } \right] \le \frac{\left(3-\gamma \right)}{{{\sigma' _1}}}D\left( {{x^*},{x_0}} \right)+\mathbb{E}\left[ {\sum\limits_{s = 0}^{S - 1} {\frac{\left(1-\gamma\right)}{\sigma'_1}\left(\frac{{8\rho L_f^2}}{\alpha } - 1\right)D\left( {{x_{s + 1}},{x_{s-1}}} \right)} } \right].
\end{equation}
Leveraging the fact $w_s=x_s$, \eqref{2.1}, \eqref{theorem2.2} and Lemma \ref{l3.1}, we can get
\[\begin{aligned}
{\mathbb{E}}\left[ {\sum\limits_{s = 0}^{S - 1} {{{\left\| {{u_s}} \right\|}_*^2}} } \right]
 &\le \mathbb{E}\left[ {\sum\limits_{s = 0}^{S - 1} {\left[ {\left( {4{L^2} + \frac{{4\gamma L_f^2}}{{{\alpha ^2}}}} \right){{\left\| {{x_{s + 1}} - {x_s}} \right\|}^2} + 4{L^2}{{\left\| {{x_s} - {x_{s - 1}}} \right\|}^2}} \right]} } \right]\\
 &+ \mathbb{E}\left[ {\sum\limits_{s = 0}^{S - 1} {\left[ {\frac{{4\left( {1 - \gamma } \right)L_f^2}}{{{\alpha ^2}}}{{\left\| {{x_{s + 1}} - {x_{s - 1}}} \right\|}^2} + 4\left\| {{\mathbb{E}_s}\left[ {{\Delta ^s}} \right] - {\Delta ^s}} \right\|_*^2} \right]} } \right]\\
  &+{\mathbb{E}}\left[ {\sum\limits_{s = 0}^{S - 1} {\left[ {\left( {8{L^2} + \frac{{8\gamma L_f^2}}{{{\alpha ^2}}}} \right)D\left( {{x_{s + 1}},{x_s}} \right) +\left( {8L^2+\frac{8{\bar L}^2}{b}}\right)D\left( {{x_s},{x_{s - 1}}} \right)} \right]} } \right]\\
 &+{\mathbb{E}}\left[ {\sum\limits_{s = 0}^{S - 1} {\left[ {  \frac{{8 \left(1-\gamma \right) L_f^2}}{{{\alpha ^2}}}} D\left(x_{s+1},x_{s-1}\right) \right]} } \right].
\end{aligned}\]
Since $x_0=x_{-1}$, we have $\sum\limits_{s = 0}^{S - 1} {D\left( {{x_s},{x_{s - 1}}} \right)}  \le \sum\limits_{s = 0}^{S - 1} {D\left( {{x_{s + 1}},{x_s}} \right)} $.\\
Therefore
\begin{equation}\label{theorem2.6'}
\mathbb{E}\left[ {\sum\limits_{s = 0}^{S - 1} {{{\left\| {{u_s}} \right\|}_*^2}} } \right] \le \mathbb{E}\left[ {\sum\limits_{s = 0}^{S - 1} {{\sigma' _2}D\left( {{x_{s + 1}},{x_s}} \right)}{ + \frac{{8 \left(1-\gamma \right) L_f^2}}{{{\alpha ^2}}}} D\left(x_{s+1},x_{s-1}\right) } \right].
\end{equation}
Combining \eqref{theorem2.5'} with \eqref{theorem2.6'}, we have
\[
\mathbb{E}\left[ {\sum\limits_{s = 0}^{S - 1} {{{\left\| {{u_s}} \right\|}_*^2}} } \right] \le \frac{{\left(3-\gamma \right){\sigma' _2}}}{{{\sigma' _1}}}D\left( {{x^*},{x_0}} \right),
\]
where we use the fact that $\left( {1 - \gamma } \right)\left( {\frac{{{\sigma' _1}}}{{{\sigma' _2}}}\left( {\frac{{8\rho L_f^2}}{\alpha } - 1} \right) + \frac{{8L_f^2}}{{{\alpha ^2}}}} \right) \le 0$.\\
 The proof is finished.
\end{proof}

\begin{corollary}
Under the conditions of Theorem \ref{theorem2} or Theorem \ref{theor3.5}, there is $s \in \left\{ {1, \cdots ,S } \right\}$ such that
\[\mathbb{E}\left[ {dist\left( {\mathbf{0},F\left( {{x_{s}}} \right) + \partial \left( {g + {{\rm I}_K}} \right)\left( {{x_{s}}} \right)} \right)} \right] \le \epsilon, \]
within $S=\mathcal{O}\left(\frac{1}{ \epsilon^2}\right)$ iterations, that is, to achieve an $\epsilon$-approximated solution, the total complexity in $S$ iterations is $\mathcal{O}\left(\frac{1}{\epsilon^2}\right)$.
\end{corollary}
\begin{proof}
Defining $t = \mathop {\arg \min }\limits_{s = \left\{ {0,1, \cdots ,S - 1} \right\}} \left\{ {\left\| {{u_s}} \right\|_*^2} \right\}$ and using Theorem \ref{theorem2}, we can get
\begin{equation}\label{theorem2.6+1}
  \mathbb{E}\left[ {\left\| {{u_t}} \right\|_*^2} \right] \le \frac{{\left(3-\gamma\right){\sigma _2}}}{{S{\sigma _1}}}D\left( {{x^*},{x_0}} \right), ~\forall~ S \ge 1.
\end{equation}
Let ${A_{s + 1}} = F\left( {{x_{s + 1}}} \right) + \partial \left( {g + {{\rm I}_K}} \right)\left( {{x_{s + 1}}} \right)$, it is easy to see that ${{u_s}}\in {A_{s + 1}}$. We use Jensen's inequality
\begin{equation}\label{theorem2.7}
\mathbb{E}\left[ {\left\| {{u_t}} \right\|_*^2} \right] \ge {\left( {\mathbb{E}\left[ {{{\left\| {{u_t}} \right\|}_*}} \right]} \right)^2} \ge {\left( {\mathbb{E}\left[ {\mathop {\min }\limits_{u \in {A_{t + 1}}} {{\left\| u \right\|}_*}} \right]} \right)^2} = {\left( {\mathbb{E}\left[ {dist\left( {\mathbf{0},{A_{t + 1}}} \right)} \right]} \right)^2}.
\end{equation}
 Using \eqref{theorem2.6+1} and \eqref{theorem2.7}, we have
\[\mathbb{E}\left[ {dist\left( {\mathbf{0},{A_{t + 1}}} \right)} \right] \le \sqrt {\frac{{\left(3-\gamma\right){\sigma _2}}}{{S{\sigma _1}}}D\left( {{x^*},{x_0}} \right)}. \]
This lead to the desired result. Similarly, we can prove that under the conditions of Theorem \ref{theor3.5}, the conclusion also holds.
This proof is completed.
\end{proof}

In order to study the linearly convergence rate of Algorithm 1, we give the following assumption.\\
$\mathbf{Assumption}$ (A7): For some $\beta>0$, the operator $F$ satisfies
\[
\left\langle {F(x),y - x} \right\rangle  + g\left( y \right) - g\left( x \right) \ge 0 \Rightarrow \left\langle {F(y),y - x} \right\rangle  + g\left( y \right) - g\left( x \right) \ge \beta {D}\left( {x,y} \right),\forall x \in {\rm{dom}}f,y \in {\rm{dom}}\left(\partial f\right).
\]
\begin{remark}
If $f\left( x \right) = \frac{1}{2}{\left\| x \right\|^2}$ and $g\left( x \right) \equiv 0$, then this
 assumption can reduce to the case that $F$ is a $\beta$-strongly pseudomonotone mapping. In addition, taking $f\left( x \right) = x\log x$, $g\left( x \right) \equiv 0$ and $F\left( x \right) = \left( {x + 1} \right){e^x},~x>0$, thus $F$ satisfies Assumption (A7) with $\beta=1$.

  A sequence $\left\{ {{x_n}} \right\} \subset \mathcal{X}$ is $Q$-linearly convergent to some point $p \in K$ if and only if there exists a constant $q \in (0,1)$ such that for all $n$ sufficiently large, the inequality $\left\| {{x_{n + 1}} - p} \right\| \le q\left\| {{x_n} - p} \right\|$ holds. A sequence $\left\{ {{x_n}} \right\} \subset \mathcal{X}$ is $R$-linearly convergent if and only if for all $n$ sufficiently large, $\left\| {{x_n} - p} \right\| \le {b_n}$ , where $\{ {b_n}\}  \subset \mathbb{R}$ converges $Q$-linearly to zero.

\end{remark}
\begin{theorem}
\label{T3.2}
Suppose that Assumptions (A1), (A3), (A4) and (A7) hold,  $\frac{{\frac{1}{2} + \gamma }}{\beta } < \alpha  < \frac{{\sqrt {{b^2}{L^4} + 8b{{\bar L}^2}\left( {\frac{{1 - \gamma }}{{1 - p}}} \right)}-b{L^2} }}{{4{{\bar L}^2}}}$, $\alpha  \le \frac{{\gamma  - p}}{{1 - p}}$, $0<p<1$, $p<\gamma < 1$, $x'\in S_1$ and the sequence $\left\{ {{x_s}} \right\}$ is from Algorithm \ref{Alg1}. Then
\[\mathbb{E}\left[ {\left\| {x' - {x_s}} \right\|} \right] \le \sqrt {\frac{4}{{\left(1+\gamma-\alpha\right){\theta ^{s}}}}D\left( {x',{x_0}} \right)}, \]
where $1 < \theta  \le \min \left\{ {\frac{{\frac{1}{2} + \alpha \beta  + \frac{\alpha }{2}}}{{1 + \gamma  + \frac{\alpha }{2}}},\frac{{\frac{{1 - \gamma }}{{1 - p}} + 2\alpha {L^2}}}{{\frac{{2{\alpha ^2}{{\bar L}^2}}}{b} + 3\alpha {L^2}}},\frac{1}{{1 - \gamma }}} \right\}$.
\end{theorem}
\begin{proof}
Owing to ${x'}\in S_1$, it yields
\[\left\langle {F({{x'}}),y - {{x'}}} \right\rangle  + g\left( y \right) - g\left( {{{x'}}} \right) \ge 0,~\forall y \in K.\]
Combing this and Assumption (A7), we have
\begin{equation}\label{T3.2.1}
\left\langle {F(y),y - {{x'}}} \right\rangle  + g\left( y \right) - g\left( {{{x'}}} \right) \ge \beta{D}\left( {{x'},y} \right),~\forall y \in K.
\end{equation}
Letting $u={x'}$ in \eqref{t1.1} and $y=x_{s+1}$ in \eqref{T3.2.1}, we have
\begin{equation}\label{3.44}
  \begin{aligned}
\alpha \beta D\left( {x',{x_{s + 1}}} \right) &\le  - \alpha \left\langle {F\left( {{w_s}} \right) - F\left( {{x_{s + 1}}} \right),{x_{s + 1}} - x'} \right\rangle  + \alpha \left\langle {F\left( {{w_{s - 1}}} \right) - F\left( {{x_s}} \right),{x_s} - x'} \right\rangle \\
 &- \alpha \left\langle {F\left( {{x_s}} \right) - F\left( {{w_{s - 1}}} \right),{x_{s + 1}} - {x_s}} \right\rangle  + \alpha \left\langle {{\mathbb{E}_s}\left[ {{\Delta ^s}} \right] - {\Delta ^s},{x_{s + 1}} - {x_s}} \right\rangle \\
 &+ \alpha \left\langle {{\mathbb{E}_s}\left[ {{\Delta ^s}} \right] - {\Delta ^s},{x_s} - x'} \right\rangle  + \gamma \left( {D\left( {x',{x_s}} \right) - D\left( {{x_{s + 1}},{x_s}} \right)} \right)\\
&+ \left( {1 - \gamma } \right)\left( {D\left( {x',{x_{s - 1}}} \right) - D\left( {{x_{s + 1}},{x_{s - 1}}} \right)} \right) - D\left( {x',{x_{s + 1}}} \right).
\end{aligned}
\end{equation}
Rearranging the order of \eqref{3.44} and taking expectation on both sides of \eqref{3.44}, we can get
\begin{equation}\label{t3.2.3}
\begin{array}{*{20}{l}}
{\begin{array}{*{20}{l}}
{\mathbb{E}\left[ {\left( {\alpha \beta  + 1} \right)D\left( {x',{x_{s + 1}}} \right) + \alpha \left\langle {F\left( {{w_s}} \right) - F\left( {{x_{s + 1}}} \right),{x_{s + 1}} - x'} \right\rangle } \right]}\\
\begin{array}{l}
 \le \mathbb{E}\left[ {\gamma D\left( {x',{x_s}} \right) + \alpha \left\langle {F\left( {{w_{s - 1}}} \right) - F\left( {{x_s}} \right),{x_s} - x'} \right\rangle } \right]\\
 + \mathbb{E}\left[ { - \alpha \left\langle {F\left( {{x_s}} \right) - F\left( {{w_{s - 1}}} \right),{x_{s + 1}} - {x_s}} \right\rangle  + \alpha \left\langle {{\mathbb{E}_s}\left[ {{\Delta ^s}} \right] - {\Delta ^s},{x_{s + 1}} - {x'}} \right\rangle } \right]
\end{array}
\end{array}}\\
{\begin{array}{*{20}{c}}
{ + \mathbb{E}\left[ {  - \gamma D\left( {{x_{s + 1}},{x_s}} \right) + \left( {1 - \gamma } \right)\left( {D\left( {x',{x_{s - 1}}} \right) - D\left( {{x_{s + 1}},{x_{s - 1}}} \right)} \right)} \right]}\\
{}
\end{array}}
\end{array}
\end{equation}
Using \eqref{t1.2} and \eqref{2.1}, we have
\begin{equation}\label{t3.2.4}
 - \alpha \left\langle {F\left( {{x_s}} \right) - F\left( {{w_{s - 1}}} \right),{x_{s + 1}} - {x_s}} \right\rangle  \le \frac{{\alpha L^2}}{2}{\left\| {{x_{s}} - {w_{s-1}}} \right\|^2} + \alpha D\left( {{x_{s + 1}},{x_s}} \right).
\end{equation}
Using the Young's inequality, \eqref{2.1} and Lemma \ref{l3.1}, we have
\begin{equation}\label{t3.2.5}
\begin{aligned}
{\mathbb{E}}\left[ {\alpha \left\langle {{{\mathbb{E}}_s}\left[ {{\Delta ^s}} \right] - {\Delta ^s},{x_{s + 1}} - {{x'}}} \right\rangle } \right] &\le {\alpha ^2}{\mathbb{E}}\left[ {{{\left\| {{{\mathbb{E}}_s}\left[ {{\Delta ^s}} \right] - {\Delta ^s}} \right\|}_*^2}} \right] + \frac{1}{4}{\mathbb{E}}\left[ {{{\left\| {{x_{s + 1}} - {{x'}}} \right\|}^2}} \right]\\
 &\le \frac{{\alpha ^2}{{\bar L}^2}}{b}{\mathbb{E}}\left[ {{{\left\| {{x_s} - {w_{s - 1}}} \right\|}^2}} \right] + \frac{1}{4}{\mathbb{E}}\left[ {{{\left\| {{x_{s + 1}} - {{x'}}} \right\|}^2}} \right]\\
 &\le \frac{{\alpha ^2}{{\bar L}^2}}{b}{\mathbb{E}}\left[ {{{\left\| {{x_s} - {w_{s - 1}}} \right\|}^2}} \right] + \frac{1}{2}{\mathbb{E}}[D\left( {{{x'}},{x_{s + 1}}} \right)].
\end{aligned}
\end{equation}
Substituting \eqref{t3.2.4} and \eqref{t3.2.5} into \eqref{t3.2.3}, we can get
\begin{equation}\label{t3.2.6}
\begin{array}{*{20}{l}}
\begin{array}{l}
\mathbb{E}\left[ {\left( {\frac{1}{2} + \alpha \beta } \right)D\left( {x',{x_{s + 1}}} \right) + D\left( {x',{x_s}} \right) + \alpha \left\langle {F\left( {{w_s}} \right) - F\left( {{x_{s + 1}}} \right),{x_{s + 1}} - x'} \right\rangle } \right]\\
 + \mathbb{E}\left[ {\left( {\gamma  - \alpha } \right)D\left( {{x_{s + 1}},{x_s}} \right) + \left( {1 - \gamma } \right)D\left( {{x_{s + 1}},{x_{s - 1}}} \right)} \right]
\end{array}\\
\begin{array}{l}
 \le \mathbb{E}\left[ {\left( {1 + \gamma } \right)D\left( {x',{x_s}} \right) + \left( {1 - \gamma } \right)D\left( {x',{x_{s - 1}}} \right) + \alpha \left\langle {F\left( {{w_{s - 1}}} \right) - F\left( {{x_s}} \right),{x_s} - x'} \right\rangle } \right]\\
+\mathbb{E}\left[ {\left( {\frac{{{\alpha ^2}{{\bar L}^2}}}{b} + \frac{\alpha {L^2}}{2}} \right){{{\left\| {{x_s} - {w_{s - 1}}} \right\|}^2}}} \right].
\end{array}
\end{array}
\end{equation}
Define
\[\begin{array}{l}
{{a}_s}{\rm{ = }}\left( {1 + \gamma } \right)D\left( {x',{x_s}} \right) + \left( {1 - \gamma } \right)D\left( {x',{x_{s - 1}}} \right) + {\left( {\frac{{{\alpha ^2}{{\bar L}^2}}}{b} + \frac{\alpha {L^2}}{2}} \right){{{\left\| {{x_s} - {w_{s - 1}}} \right\|}^2}}}\\
 + \alpha \left\langle {F\left( {{w_{s - 1}}} \right) - F\left( {{x_s}} \right),{x_s} - x'} \right\rangle.
\end{array}\]
Using the Young's inequality, Assumption (A4) and (\ref{2.1}), we can get
\[
\begin{aligned}
\alpha \left\langle {F\left( {{w_{s - 1}}} \right) - F\left( {{x_s}} \right),{x_s} - {{x'}}} \right\rangle  &\ge  - \frac{\alpha}{2}{\left\| {F\left( {{w_{s - 1}}} \right) - F\left( {{x_s}} \right)} \right\|_*^2} - \frac{\alpha}{2}{\left\| {{x_s} - {{x'}}} \right\|^2}\\
 &\ge  -\frac{{\alpha L^2}}{2}{\left\| {{x_{s}} - {w_{s-1}}} \right\|^2} -\alpha D\left( {{{x'}},{x_s}} \right).
\end{aligned}
\]
On the one hand, using this, $p<\gamma<1$ and $\alpha  \le \frac{{\gamma  - p}}{{1 - p}}=\gamma-\frac{p\left(1-\gamma\right)}{1-p}$, we have
\begin{equation}\label{t3.2.7}
\begin{aligned}
{a_s} &\ge \left( {1 + \gamma  - \alpha } \right)D\left( {x',{x_s}} \right) + \left( {1 - \gamma } \right)D\left( {x',{x_{s - 1}}} \right) + \frac{{{\alpha ^2}{{\bar L}^2}}}{b}{\left\| {{x_{s}} - {w_{s-1}}} \right\|^2}\\
 &\ge \left( {1 + \gamma  - \alpha } \right)D\left( {x',{x_s}} \right) \ge 0.
\end{aligned}
\end{equation}
The Young's inequality,  Assumption (A4) and (\ref{2.1}) imply that
\begin{equation}\label{theorem4.4}
\begin{aligned}
\alpha \left\langle {F\left( {{w_{s}}} \right) - F\left( {{x_{s+1}}} \right),{x_{s+1}} - {{x'}}} \right\rangle  &\le {\alpha }{\left\| {F\left( {{w_{s}}} \right) - F\left( {{x_{s+1}}} \right)} \right\|_*^2} + \frac{\alpha}{4}{\left\| {{x_{s+1}} - {{x'}}} \right\|^2}\\
 &\le {\alpha }{L^2}{\left\| {{x_{s + 1}} - {w_s}} \right\|^2} + \frac{\alpha}{2}D\left( {{{x'}},{x_{s+1}}} \right).
\end{aligned}
\end{equation}
Since $\frac{{\frac{1}{2} + \gamma }}{\beta } < \alpha  < \frac{ \sqrt {{b^2}{L^4} + 8b{{\bar L}^2}\left( {\frac{{1 - \gamma }}{{1 - p}}} \right)} -{b{L^2}}}{{4{{\bar L}^2}}}$, $\alpha  \le \frac{{\gamma  - p}}{{1 - p}}$ and $p < \gamma  < 1$, there exists a constant $\theta$ such that
$1 < \theta  \le \min \left\{ {\frac{{\frac{1}{2} + \alpha \beta  + \frac{\alpha }{2}}}{{1 + \gamma  + \frac{\alpha }{2}}},\frac{{\frac{{1 - \gamma }}{{1 - p}} + 2\alpha {L^2}}}{{\frac{{2{\alpha ^2}{{\bar L}^2}}}{b} + 3\alpha {L^2}}},\frac{1}{{1 - \gamma }}} \right\}$. Combining this, \eqref{2.1}, \eqref{t1.9} and \eqref{theorem4.4}, we can obtain
\begin{equation}\label{t3.2.8}
\begin{array}{*{20}{l}}
\begin{array}{l}
\mathbb{E}\left[ {\left( {\frac{1}{2} + \alpha \beta } \right)D\left( {x',{x_{s + 1}}} \right) + D\left( {x',{x_s}} \right) + \alpha \left\langle {F\left( {{w_s}} \right) - F\left( {{x_{s + 1}}} \right),{x_{s + 1}} - x'} \right\rangle  - \theta {a_{s + 1}}} \right]\\
 + \mathbb{E}\left[ {\left( {\gamma  - \alpha } \right)D\left( {{x_{s + 1}},{x_s}} \right) + \left( {1 - \gamma } \right)D\left( {{x_{s + 1}},{x_{s - 1}}} \right)} \right]
\end{array}\\
{ = \mathbb{E}\left[ {\left( {\frac{1}{2} + \alpha \beta  - \theta \left( {1 + \gamma } \right)} \right)D\left( {x',{x_{s + 1}}} \right) - \theta \left( {\frac{{{\alpha ^2}{{\bar L}^2}}}{b} + \frac{\alpha {L^2}}{2}} \right){\left\| {{x_{s + 1}} - {w_s}} \right\|^2}} \right]}\\
\begin{array}{l}
 - \mathbb{E}\left[ {\alpha \left( {\theta  - 1} \right)\left\langle {F\left( {{w_s}} \right) - F\left( {{x_{s + 1}}} \right),{x_{s + 1}} - x'} \right\rangle  + \left( {1 - \theta \left( {1 - \gamma } \right)} \right)D\left( {x',{x_s}} \right)} \right]\\
 + \mathbb{E}\left[ {\left( {\gamma  - \alpha } \right)D\left( {{x_{s + 1}},{x_s}} \right) + \left( {1 - \gamma } \right)D\left( {{x_{s + 1}},{x_{s - 1}}} \right)} \right]
\end{array}\\
\begin{array}{l}
 \ge \mathbb{E}\left[ {\left( {\frac{1}{2} + \alpha \beta  - \theta \left( {1 + \gamma } \right) - \frac{{\theta  - 1}}{2}\alpha } \right)D\left( {x',{x_{s + 1}}} \right) + \left( {1 - \theta \left( {1 - \gamma } \right)} \right)D\left( {x',{x_s}} \right)} \right]\\
 - \mathbb{E}\left[ {\left( {\theta \left( {\frac{{{\alpha ^2}{{\bar L}^2}}}{b} + \frac{3\alpha {L^2}}{2}} \right) - \alpha {L^2}} \right){\left\| {{x_{s + 1}} - {w_s}} \right\|^2}} \right]\\
 + \mathbb{E}\left[ {\left( {\gamma  - \alpha } \right)D\left( {{x_{s + 1}},{x_s}} \right) + \left( {1 - \gamma } \right)D\left( {{x_{s + 1}},{x_{s - 1}}} \right)} \right]\\
 =
\mathbb{E}\left[ {\left( {\frac{1}{2} + \alpha \beta  - \theta \left( {1 + \gamma } \right) - \frac{{\theta  - 1}}{2}\alpha } \right)D\left( {x',{x_{s + 1}}} \right) + \left( {1 - \theta \left( {1 - \gamma } \right)} \right)D\left( {x',{x_s}} \right)} \right]\\
 + \mathbb{E}\left[ {\frac{{1 - \gamma }}{{2\left( {1 - p} \right)}}{{\left\| {{x_{s + 1}} - {w_s}} \right\|}^2}} \right] - \mathbb{E}\left[ {\frac{{1 - \gamma }}{{2\left( {1 - p} \right)}}{{\left\| {{x_{s + 1}} - {w_s}} \right\|}^2}} \right]\\
 - \mathbb{E}\left[ {\left( {\theta \left( {\frac{{{\alpha ^2}{{\bar L}^2}}}{b} + \frac{{3\alpha {L^2}}}{2}} \right) - \alpha {L^2}} \right){{\left\| {{x_{s + 1}} - {w_s}} \right\|}^2}} \right]
\end{array}\\
{ + \mathbb{E}\left[ {\left( {\gamma  - \alpha } \right)D\left( {{x_{s + 1}},{x_s}} \right) + \left( {1 - \gamma } \right)D\left( {{x_{s + 1}},{x_{s - 1}}} \right)} \right]}\\
\begin{array}{l}
 = \mathbb{E}\left[ {\left( {\frac{1}{2} + \alpha \beta  - \theta \left( {1 + \gamma } \right) - \frac{{\theta  - 1}}{2}\alpha } \right)D\left( {x',{x_{s + 1}}} \right) + \left( {1 - \theta \left( {1 - \gamma } \right)} \right)D\left( {x',{x_s}} \right)} \right]\\
 + \mathbb{E}\left[ {\left( {\frac{{1 - \gamma }}{{2\left(1 - p\right)}} + \alpha {L^2} - \theta \left( {\frac{{{\alpha ^2}{{\bar L}^2}}}{b} + \frac{3\alpha {L^2}}{2}} \right)} \right){\left\| {{x_{s + 1}} - {w_s}} \right\|^2}} \right]\\
 + \mathbb{E}\left[ {\left( {\gamma - \alpha } \right)D\left( {{x_{s + 1}},{x_s}} \right)-\frac{{p\left( {1 - \gamma } \right)}}{{2\left( {1 - p} \right)}}{\left\| {{x_{s + 1}} - {x_s}} \right\|^2}} \right]\\
 +\mathbb{E}\left[ {\left( {1 - \gamma } \right)\left( {D\left( {{x_{s + 1}},{x_{s - 1}}} \right) - \frac{1}{2}{{\left\| {{x_{s + 1}} - {x_{s - 1}}} \right\|}^2}} \right)} \right] \\
 \ge
 \mathbb{E}\left[ {\left( {\frac{1}{2} + \alpha \beta  - \theta \left( {1 + \gamma } \right) - \frac{{\theta  - 1}}{2}\alpha } \right)D\left( {x',{x_{s + 1}}} \right) + \left( {1 - \theta \left( {1 - \gamma } \right)} \right)D\left( {x',{x_s}} \right)} \right]\\
 + \mathbb{E}\left[ {\left( {\frac{{1 - \gamma }}{{2\left(1 - p\right)}} + \alpha {L^2} - \theta \left( {\frac{{{\alpha ^2}{{\bar L}^2}}}{b} + \frac{3\alpha {L^2}}{2}} \right)} \right){\left\| {{x_{s + 1}} - {w_s}} \right\|^2}} \right]\\
 + \frac{1}{2}\mathbb{E}\left[ {\left( {\frac{{\gamma  - p}}{{1 - p}} - \alpha } \right){{\left\| {{x_{s + 1}} - {x_s}} \right\|}^2}} \right] \ge 0.
\end{array}
\end{array}
\end{equation}
Combining (\ref{t3.2.6}) and (\ref{t3.2.8}), we have
\[{\mathbb{E}}\left[ {{a_{s + 1}}} \right] \le \frac{1}{\theta }{\mathbb{E}}\left[ {{a_s}} \right] \le  \cdots  \le \frac{1}{{{\theta ^{s + 1}}}}{\mathbb{E}}\left[ {{a_0}} \right] = \frac{{2 }}{{{\theta ^{s + 1}}}}{D\left( {{{x'}},{x_0}} \right)}.\]
In view of \eqref{t3.2.7}, we can get
\[{\mathbb{E}}\left[ {D\left( {{{x'}},{x_s}} \right)} \right] \le \frac{2}{{\left(1+\gamma-\alpha\right){\theta ^{s }}}}{D\left( {{{x'}},{x_0}} \right)}.\]
Using \eqref{2.1} and Jensen's inequality, we can get
\[\mathbb{E}\left[ {\left\| {x' - {x_s}} \right\|} \right] \le \sqrt {\frac{4}{{\left(1+\gamma-\alpha\right){\theta ^{s}}}}D\left( {x',{x_0}} \right)}. \]
We have completed this proof.
\end{proof}
When $p=1$, that is $w_s=x_s$, we can get the following result.
\begin{theorem}
\label{T3.2'}
Let Assumptions (A1), (A3), (A4) and (A7) hold, $p=1$, $\frac{{\frac{1}{2} + \gamma }}{\beta } < \alpha  < \frac{{\sqrt {{b^2}{{\left( {{L^2} + 1} \right)}^2} + 8b\gamma {{\bar L}^2}}  - b\left( {{L^2} + 1} \right)}}{{4{{\bar L}^2}}}$ and the sequence $\left\{ {{x_s}} \right\}$ is from Algorithm \ref{Alg1}. Then
\[\mathbb{E}\left[ {\left\| {x' - {x_s}} \right\|} \right] \le \sqrt {\frac{4}{{\left(1+\gamma-\alpha\right){\varsigma ^{s}}}}D\left( {x',{x_0}} \right)}, \]
where $1 < \varsigma  \le \min \left\{ {\frac{{\frac{1}{2} + \alpha \beta  + \frac{\alpha }{2}}}{{1 + \gamma  + \frac{\alpha }{2}}},\frac{{\gamma  - \alpha  + 2\alpha {L^2}}}{{\frac{{2\alpha^2 {{\bar L}^2}}}{b} + 3\alpha {L^2}}}} \right\}$ and $1-\left( {1 - \gamma } \right)\varsigma \ge 0$, and $x'\in S_1$.
\end{theorem}
\begin{proof}
Using the fact $w_s=x_s$ and \eqref{2.1} in \eqref{t3.2.6}, we can obtain
\begin{equation}\label{t3.2.6'}
\begin{array}{*{20}{l}}
\begin{array}{l}
\mathbb{E}\left[ {\left( {\frac{1}{2} + \alpha \beta } \right)D\left( {x',{x_{s + 1}}} \right) + D\left( {x',{x_s}} \right) + \alpha \left\langle {F\left( {{x_s}} \right) - F\left( {{x_{s + 1}}} \right),{x_{s + 1}} - x'} \right\rangle } \right]\\
 + \mathbb{E}\left[ {\left( {\gamma  - \alpha } \right)D\left( {{x_{s + 1}},{x_s}} \right) + \left( {1 - \gamma } \right)D\left( {{x_{s + 1}},{x_{s - 1}}} \right)} \right]
\end{array}\\
\begin{array}{l}
 \le \mathbb{E}\left[ {\left( {1 + \gamma } \right)D\left( {x',{x_s}} \right) + \left( {1 - \gamma } \right)D\left( {x',{x_{s - 1}}} \right) + \alpha \left\langle {F\left( {{x_{s - 1}}} \right) - F\left( {{x_s}} \right),{x_s} - x'} \right\rangle } \right]\\
+\mathbb{E}\left[ {\left( {\frac{{2{\alpha ^2}{{\bar L}^2}}}{b} + \alpha {L^2}} \right)D\left( {{x_s},{x_{s - 1}}} \right)} \right].
\end{array}
\end{array}
\end{equation}
Define
\[\begin{array}{l}
{a'_s}{\rm{ = }}\left( {1 + \gamma } \right)D\left( {x',{x_s}} \right) + \left( {1 - \gamma } \right)D\left( {x',{x_{s - 1}}} \right) + \left( {\frac{{2{\alpha ^2}{{\bar L}^2}}}{b} + \alpha {L^2}} \right)D\left( {{x_s},{x_{s - 1}}} \right)\\
 + \alpha \left\langle {F\left( {{x_{s - 1}}} \right) - F\left( {{x_s}} \right),{x_s} - x'} \right\rangle.
\end{array}\]
Using the Young's inequality, Assumption (A4) and (\ref{2.1}), we can get
\[
\begin{aligned}
\alpha \left\langle {F\left( {{x_{s - 1}}} \right) - F\left( {{x_s}} \right),{x_s} - {{x'}}} \right\rangle  &\ge  - \frac{\alpha}{2}{\left\| {F\left( {{x_{s - 1}}} \right) - F\left( {{x_s}} \right)} \right\|_*^2} - \frac{\alpha}{2}{\left\| {{x_s} - {{x'}}} \right\|^2}\\
 &\ge  - {\alpha}{L^2}D\left( {{x_s},{x_{s - 1}}} \right) -\alpha D\left( {{{x'}},{x_s}} \right).
\end{aligned}
\]
Thanks to $\alpha  < \frac{{\sqrt {{b^2}{{\left( {{L^2} + 1} \right)}^2} + 8b\gamma {{\bar L}^2}}  - b\left( {{L^2} + 1} \right)}}{{4{{\bar L}^2}}}$, we have $\gamma  - \alpha  > \frac{{2{\alpha ^2}{{\bar L}^2}}}{b} + \alpha {L^2} >0$. Combing this and $0 < \gamma  \le 1$,  we imply that
\begin{equation}\label{t3.2.7'}
\begin{aligned}
{a'_s} &\ge \left( {1 + \gamma  - \alpha } \right)D\left( {x',{x_s}} \right) + \left( {1 - \gamma } \right)D\left( {x',{x_{s - 1}}} \right) + \frac{{2{\alpha ^2}{{\bar L}^2}}}{b}D\left( {{x_s},{x_{s - 1}}} \right)\\
 &\ge \left( {1 + \gamma  - \alpha } \right)D\left( {x',{x_s}} \right) \ge 0.
\end{aligned}
\end{equation}
 Using the Young's inequality, Assumption (A4) and \eqref{2.1}, we obtain
\begin{equation}\label{theorem4.4'}
\begin{aligned}
\alpha \left\langle {F\left( {{x_{s}}} \right) - F\left( {{x_{s+1}}} \right),{x_{s+1}} - {{x'}}} \right\rangle  &\le {\alpha }{\left\| {F\left( {{x_{s}}} \right) - F\left( {{x_{s+1}}} \right)} \right\|_*^2} + \frac{\alpha}{4}{\left\| {{x_{s+1}} - {{x'}}} \right\|^2}\\
 &\le 2{\alpha }{L^2}D\left( {{x_{s+1}},{x_{s}}} \right) + \frac{\alpha}{2}D\left( {{{x'}},{x_{s+1}}} \right).
\end{aligned}
\end{equation}
Since $\frac{{\frac{1}{2} + \gamma }}{\beta } < \alpha  < \frac{{\sqrt {{b^2}{{\left( {{L^2} + 1} \right)}^2} + 8b\gamma {{\bar L}^2}}  - b\left( {{L^2} + 1} \right)}}{{4{{\bar L}^2}}}$ and $0 < \gamma  \le 1$, there exists a constant $\varsigma$ with
$1 < \varsigma  \le \min \left\{ {\frac{{\frac{1}{2} + \alpha \beta  + \frac{\alpha }{2}}}{{1 + \gamma  + \frac{\alpha }{2}}},\frac{{\gamma  - \alpha  + 2\alpha {L^2}}}{{\frac{{2\alpha^2 {{\bar L}^2}}}{b} + 3\alpha {L^2}}}} \right\}$ and $1-\left( {1 - \gamma } \right)\varsigma \ge 0$. Combining this and \eqref{theorem4.4'}, we can obtain
\begin{equation}\label{t3.2.8'}
\begin{array}{*{20}{l}}
\begin{array}{l}
\left( {\frac{1}{2} + \alpha \beta } \right)D\left( {x',{x_{s + 1}}} \right) + D\left( {x',{x_s}} \right) + \alpha \left\langle {F\left( {{x_s}} \right) - F\left( {{x_{s + 1}}} \right),{x_{s + 1}} - x'} \right\rangle  - \varsigma {{\alpha '}_{s + 1}}\\
 + \left( {\gamma  - \alpha } \right)D\left( {{x_{s + 1}},{x_s}} \right) + \left( {1 - \gamma } \right)D\left( {{x_{s + 1}},{x_{s - 1}}} \right)
\end{array}\\
\begin{array}{l}
 \ge \left( {\frac{1}{2} + \alpha \beta  - \varsigma \left( {1 + \gamma } \right)} \right)D\left( {x',{x_{s + 1}}} \right) + \left( {1 - \varsigma \left( {1 - \gamma } \right)} \right)D\left( {x',{x_s}} \right)\\
 + \left( {\gamma  - \alpha  - \varsigma \left( {\frac{{2{\alpha ^2}{{\bar L}^2}}}{b} + \alpha {L^2}} \right)} \right)D\left( {{x_{s + 1}},{x_s}} \right)
\end{array}\\
{ - \alpha \left( {\varsigma  - 1} \right)\left\langle {F\left( {{x_s}} \right) - F\left( {{x_{s + 1}}} \right),{x_{s + 1}} - x'} \right\rangle }\\
{ \ge \left( {\frac{1}{2} + \alpha \beta  - \varsigma \left( {1 + \gamma } \right) - \frac{{\varsigma  - 1}}{2}\alpha } \right)D\left( {x',{x_{s + 1}}} \right) + \left( {1 - \varsigma \left( {1 - \gamma } \right)} \right)D\left( {x',{x_s}} \right)}\\
{ + \left( {\gamma  - \alpha  - \varsigma \left( {\frac{{2{\alpha ^2}{{\bar L}^2}}}{b} + \alpha {L^2}} \right) - 2\alpha {L^2}\left( {\varsigma  - 1} \right)} \right)D\left( {{x_{s + 1}},{x_s}} \right)}\ge 0.
\end{array}
\end{equation}
Combining (\ref{t3.2.6'}) and (\ref{t3.2.8'}), we have
\[{\mathbb{E}}\left[ {{a'_{s + 1}}} \right] \le \frac{1}{\varsigma }{\mathbb{E}}\left[ {{a'_s}} \right] \le  \cdots  \le \frac{1}{{{\varsigma ^{s + 1}}}}{\mathbb{E}}\left[ {{a'_0}} \right] = \frac{{2 }}{{{\varsigma ^{s + 1}}}}{D\left( {{{x'}},{x_0}} \right)}.\]
By \eqref{t3.2.7'}, we can get
\[{\mathbb{E}}\left[ {D\left( {{{x'}},{x_s}} \right)} \right] \le \frac{2}{{\left(1+\gamma-\alpha\right){\varsigma ^{s }}}}{D\left( {{{x'}},{x_0}} \right)}.\]
Using \eqref{2.1} and Jensen's inequality, we can get
\[\mathbb{E}\left[ {\left\| {x' - {x_s}} \right\|} \right] \le \sqrt {\frac{4}{{\left(1+\gamma-\alpha\right){\varsigma ^{s}}}}D\left( {x',{x_0}} \right)}. \]
We have completed this proof.
\end{proof}
\begin{remark}
To the best of our knowledge, no linear convergence result has been established for single-loop methods in Bregman settings to solve FSVI \eqref{1.1}. In this regard, Theorems \ref{T3.2} and \ref{T3.2'} offer a novel insight.
\end{remark}
\section{Numerical experiments}
\setcounter{equation}{0}

In this section, we present two numerical experiments to evaluate the performance of the proposed Algorithm \ref{Alg1}. All computations were carried out in MATLAB 2023(b) on a personal computer (PC) with 8.00 GB of RAM. Here, "time(s)" represents the CPU running time.
\begin{example}
\label{e4.1}
Firstly,  we consider the following matrix game
\[\mathop {\min }\limits_{x \in {\Delta ^n}} \mathop {\max }\limits_{y \in {\Delta ^n}} \left\langle {Ax,y} \right\rangle,\]
where ${\Delta ^n} = \left\{ {x \in {{\mathbb{R}}^n}:\sum\limits_{i = 1}^n {{x_i}}  = 1,{x_i} \ge 0} \right\}$ and $A \in {{\mathbb{R}}^{n \times n}}$ is a given matrix. For $z = \left( {x,y} \right) \in {R^{n + n}}$, we define $\left\| z \right\| = \sqrt {\left\| x \right\|_1^2 + \left\| y \right\|_1^2}$, where ${\left\| x \right\|_1} = \sum\limits_{i = 1}^n {\left| {{x_i}} \right|}~{\rm and} ~{\left\| y \right\|_1} = \sum\limits_{i = 1}^n {\left| {{y_i}} \right|} $. Correspondingly, for a vector $z^* = \left( {x^*,y^*} \right) \in {R^{n + n}}$ in the dual space, we define ${\left\| z^* \right\|} = \sqrt {\left\| x^* \right\|_\infty ^2 + \left\| y^* \right\|_\infty ^2} $, where ${\left\| x^* \right\|_\infty } = \mathop {\max }\limits_{1 \le i \le n} \left| {{x^*_i}} \right|~{\rm and} ~{\left\| y^* \right\|_\infty } = \mathop {\max }\limits_{1 \le i \le n} \left| {{y^*_i}} \right|$. For $z = \left( {x,y} \right) \in {\Delta ^n} \times {\Delta ^n}$, we use the negative
entropy $f\left( z \right) = \sum\limits_{i = 1}^n {{x_i}\log {x_i}}  + \sum\limits_{i = 1}^n {{y_i}\log {y_i}}  = \sum\limits_{i = 1}^{2n} {{z_i}\log {z_i}}$ to generate the KL divergence $D\left( {z,v} \right) = \sum\limits_{i = 1}^{2n} {{z_i}\log \frac{{{z_i}}}{{{v_i}}} + {v_i} - {z_i}}$. Note that $f$ is 1-strong convex on ${\Delta ^n} \times {\Delta ^n}$. According to \cite{AA}, this matrix game is equivalent to FSVI (\ref{1.1}) with
\[F\left( z \right) = F\left( {x,y} \right) = \left[ \begin{array}{l}
{A^T}y\\
 - Ax
\end{array} \right], g\left( z \right) = {{\rm{I}}_{{\Delta ^n}}}\left( x \right) + {{\rm{I}}_{{\Delta ^n}}}\left( y \right).\]
Note that $F$ is Lipschitz continuous with $L={\left\| A \right\|_{\max }}$, where ${\left\| A \right\|_{\max }} = \mathop {\max }\limits_{i,k} \left| {{A_{ik}}} \right|$. Following \cite{AA}, we select the following stochastic oracle: given $u = \left( {{u^x},{u^y}} \right)$ and $ v= \left( {{v^x},{v^y}} \right)$, for $z=\left( {{x},{y}} \right)$, we define
\[{F_\xi}\left( z \right) = \left( \begin{array}{l}
\frac{1}{{{r_i}}}{A_{i:}}{y_i}\\
 - \frac{1}{{{c_k}}}{A_{:k}}{x_k}
\end{array} \right)~{\rm and}~{\rm{Pr}}\left\{ {\xi = \left( {i,k} \right)} \right\} = {r_i}{c_k},\]
 where ${r_i} = \frac{{\left| {u_i^y - v_i^y} \right|}}{{{{\left\| {{u^y} - {v^y}} \right\|}_1}}}$, ${c_k} = \frac{{\left| {u_k^x - v_k^x} \right|}}{{{{\left\| {{u^x} - {v^x}} \right\|}_1}}},$"Pr" means "Probability", ${A_{i:}}$ is $i-$th row of $A$ and ${A_{:k}}$ is $k-$th column of $A$. As shown in \cite{CT},  ${F_\xi}$ satisfies Assumption (A3) with ${\bar L}=L$. This oracle was introduced in \cite{CT}.  In our experiments, the matrix $A$ is generated by a standard normal distribution. Following \cite{AA2}, we measure the residual using the duality gap, defined as $\mathop {\max }\limits_i {\left( {Ax} \right)_i} - \mathop {\min }\limits_j {\left( {{A^T}y} \right)_j}$, and take $n=100$.  The stopping criterion is that the number of iterations does not exceed $15,000$. We compare our proposed Algorithm \ref{Alg1} (denoted as Alg1) with three existing algorithms: Algorithm 2 from \cite{AA} (VR-MP), Algorithm 1 from \cite{PA} (OM-MB), and Algorithm 1 from \cite{AA2} (VR-FoRB). The parameters for each algorithm are set as follows: \\
\textbullet VR-FoRB: $p=\frac{8.15}{\sqrt{n}},~\tau=\frac{1-\sqrt{1-p}}{3L^2},~z_0=\left(x_0,y_0\right),~x_0=y_0=\left(\frac{1}{n},\cdots,\frac{1}{n}\right)^T$;\\
\textbullet VR-MP: $K=n,~p=\frac{1}{n},~\alpha=1-p,~\tau=\frac{\sqrt{p}}{2L}$,
$z_0=\left(x_0,y_0\right),~x_0=y_0=\left(\frac{1}{n},\cdots,\frac{1}{n}\right)^T$;\\
\textbullet OM-MB: $b=4,~K=\frac{n}{3b},~\gamma=\frac{1}{K},~\eta  = \min \left\{ {\frac{{\sqrt {\gamma b} }}{{2L}},\frac{1}{{8L}}} \right\}$,
$z_0=\left(x_0,y_0\right),~x_0=y_0=\left(\frac{1}{n},\cdots,\frac{1}{n}\right)^T$;\\
\textbullet Alg1: $\alpha  = \min \left\{ {\frac{{\gamma  - p}}{{2\left( {1 - p} \right)}},\frac{{\left( {1 - \gamma} \right)b}}{{\left( {1 - p} \right)\left( {2{{\bar L}^2} + b{L^2}} \right)}}} \right\}$, $p=\frac{1}{n}$, $\gamma=p+\frac{1}{\sqrt{n}},~z_0=\left(x_0,y_0\right),~x_0=y_0=\left(\frac{1}{n},\cdots,\frac{1}{n}\right)^T$.\\

\begin{figure}[htbp]
  \centering
  \begin{subfigure}[b]{0.3\textwidth}
    \centering
    \includegraphics[width=\textwidth]{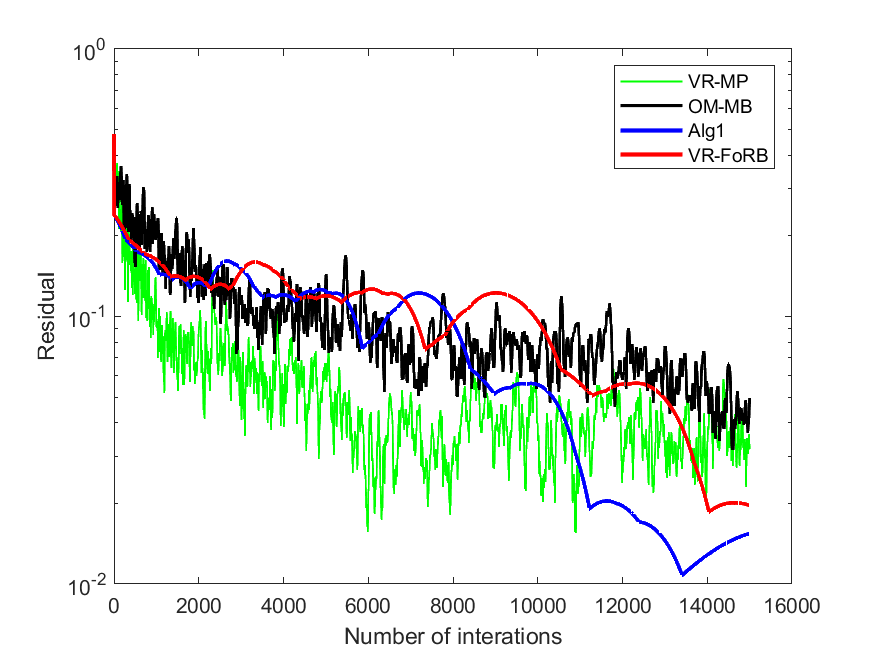}
    \label{n1fig:sub1}
    \caption{Batch size=2}
  \end{subfigure}
  \begin{subfigure}[b]{0.3\textwidth}
    \centering
    \includegraphics[width=\textwidth]{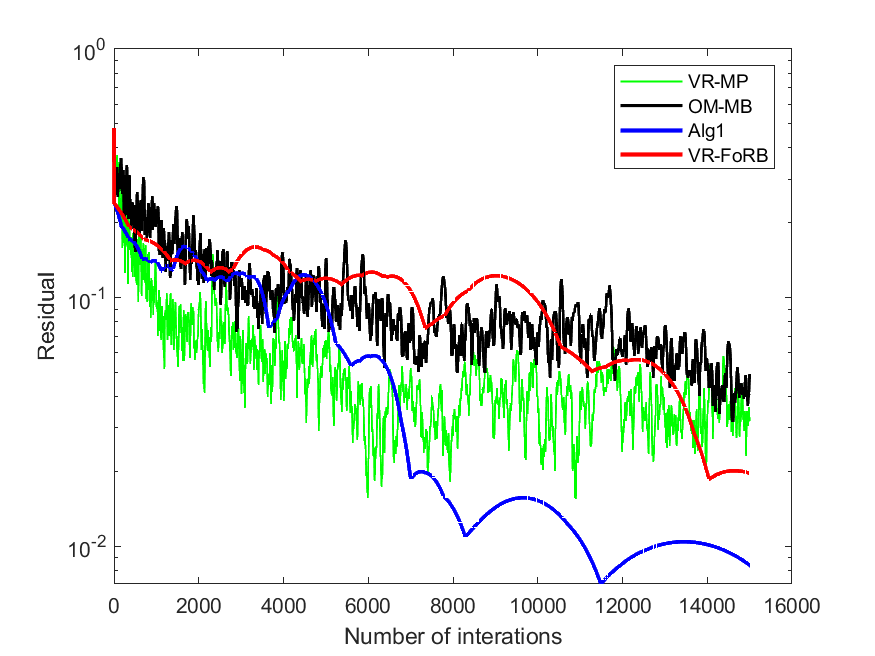}
    \label{n1fig:sub2}
    \caption{Batch size=8}
  \end{subfigure}
  \begin{subfigure}[b]{0.3\textwidth}
    \centering
    \includegraphics[width=\textwidth]{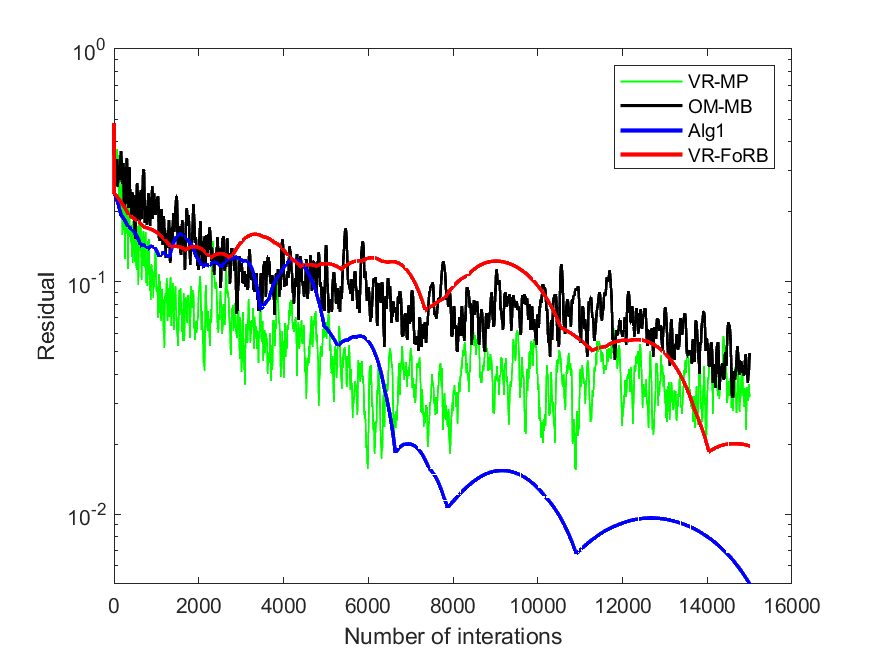}
    \label{n1fig:sub3}
    \caption{Batch size=11}
  \end{subfigure}
   \caption{Numerical behaviour of the residual for Example \ref{e4.1}}
  \label{figure 1}
\end{figure}
\begin{table}[h!]
\centering
\caption{Numerical behaviour of the residual and time for Example \ref{e4.1}}
\label{Table3}
\begin{tabular}{ccccccc}
\hline
Algorithm & VR-FoRB & VR-MP & OM-MB & Alg\ref{Alg1} ($b=2$)& Alg\ref{Alg1} ($b=8$) &Alg\ref{Alg1} ($b=11$) \\
\hline
Residual & 0.0186 & 0.0156 & 0.0316&0.0108&0.007&0.005\\
Time(s) & 10.8639 & 705.2499 & 131.0436 & 18.0333&65.6604&89.8881 \\
\hline
\end{tabular}
\end{table}
\par  As visualized in Figure \ref{figure 1} and Table \ref{Table3}, Algorithm \ref{Alg1} outperforms the other comparative methods across all tested batch sizes.
\end{example}
\begin{example}
\label{e4.2}
To evaluate the performance of our proposed method under a non-monotone setting, we consider the following two-player matrix game
\[\mathop {\min }\limits_{{{\left\| x \right\|}_2} \le 1} \mathop {\max }\limits_{{{\left\| y \right\|}_2} \le 1}  - \frac{\lambda }{2}\left\| x \right\|_2^2 + \left\langle {Ax,y} \right\rangle  + \frac{\lambda }{2}\left\| y \right\|_2^2,\]
where $x \in {\mathbb{R}}^n,y \in {\mathbb{R}}^n$, $\lambda>0$ is a constant and $A \in {{\mathbb{R}}^{n \times n}}$ is a given matrix.
Analogous to Example \ref{e4.1}, this problem can be transformed into FSVI \eqref{1.1} by letting
\[z = \left[ \begin{array}{l}
x\\
y
\end{array} \right],F\left( z \right) = \left[ \begin{array}{l}
{A^T}y - \lambda x\\
 - Ax - \lambda y
\end{array} \right],g\left( z \right) = {{\rm{I}}_K}\left( z \right),\]
where $K = \left\{ {x \in {{\mathbb{R}}^n}:{{\left\| x \right\|}_2} \le 1} \right\} \times \left\{ {y \in {{\mathbb{R}}^n}:{{\left\| y \right\|}_2} \le 1} \right\}$. For $z = \left( {x,y} \right) \in {R^{n+ n}}$, we define the norm $\left\| z \right\| = \sqrt {\left\| x \right\|_2^2 + \left\| y \right\|_2^2}$, where ${\left\| x \right\|_2} = \sqrt {\sum\limits_{i = 1}^n {x_i^2} } ~{\rm and} ~{\left\| y \right\|_2} = \sqrt {\sum\limits_{i = 1}^n {y_i^2} } $. Correspondingly, for a vector $z^* = \left( {x^*,y^*} \right) \in {R^{n + n}}$ in the dual space, its norm is defined as ${\left\| z^* \right\|} = \sqrt {\left\| x^* \right\|_2 ^2 + \left\| y^* \right\|_2 ^2} $. For $z=(x,y) \in K$, let $f\left( z \right) = \frac{1}{2}{\left\| z \right\|^2_2}$. The Bregman distance associated with $f$ is referred to as Euclidean distance given by $D\left( {z,v} \right) = \frac{1}{2}{\left\| {z - v} \right\|^2_2}$. As established in \cite{PT2}, the mapping $F$ is a non-monotone and Lipschitz continuous operator with $L=\sqrt{\lambda^2+\left\| A \right\|_2^2}$ and $F$ satisfies Assumption (A5) with $\rho  = \frac{\lambda }{{{\lambda ^2} + \left\| A \right\|_2^2}}$, where ${\left\| A \right\|_2} = \mathop {\max }\limits_{x \ne {\bf{0}}} \frac{{{{\left\| {Ax} \right\|}_2}}}{{{{\left\| x \right\|}_2}}}$. We employ the following stochastic oracle: for $z=\left( {{x},{y}} \right)$, define
\[{F_\xi }\left( z \right) = \left[ {\begin{array}{*{20}{l}}
{\frac{1}{{{r_i}}}{A_{i:}}{y_i}}\\
{ - \frac{1}{{{c_k}}}{A_{:k}}{x_k}}
\end{array}} \right]-\left[ \begin{array}{l}
\lambda x\\
\lambda y
\end{array} \right]~{\rm and}~{\rm{Pr}}\left\{ {\xi  = \left( {i,k} \right)} \right\} = {r_i}{c_k},\]
 where ${r_i} = \frac{{\left\| {{A_{i:}}} \right\|_2^2}}{{\left\| A \right\|_{{\rm{Frob}}}^2}},~{c_k} = \frac{{\left\| {{A_{:k}}} \right\|_2^2}}{{\left\| A \right\|_{{\rm{Frob}}}^2}},$ "Pr" means "Probability", ${A_{i:}}$ is $i-$th row of $A$, ${A_{:k}}$ is $k-$th column of $A$ and ${\left\| A \right\|_{{\rm{Frob}}}} = \sqrt {\sum\limits_{i,k} {A_{ik}^2} } $. According to \cite{AA}, ${F_\xi}$ satisfies Assumption (A3) with ${\bar L}=L$.

 \par In our experiments, the entries of $A$ is given by ${A_{ij}} = {\left( {\frac{{\left| {i - j} \right| + 1}}{{2n - 1}}} \right)^2},1 \le i,j \le n$. In view of \cite{AA},  we normalize $A$ such that ${\left\| A \right\|_2}=10$. Set $\lambda=0.01$. We compare our proposed Algorithm \ref{Alg1} with Algorithm 1 from \cite{AZ} (denoted as VR-FoRMAB).
 The parameters of the VR-FoRMAB is taken: $\gamma  = \beta  =1- \frac{1}{{2\sqrt n  + 1}},\theta  = 1,S = 15,q = 3,\sigma  = \frac{{5\lambda}}{{3L}}$. For our Algorithm \ref{Alg1}, we consider the following two cases:
 \\
\textbullet Alg11: $\gamma=0.94,~p=0.82,~b=15,~\alpha=\frac{8\rho}{\left(1-\frac{1-p}{1.63\left(\gamma-p\right)}\right)}\approx 0.01$,\\
\textbullet Alg12: $\gamma=1,~p=1,~b=15,~\alpha=0.015$.\\
Since $z=\mathbf{0}$ is the equilibrium of this problem, we focus on the change of the residual $\frac{{\left\| {{z_s}} \right\|}}{{\sqrt n }}$.  The stopping criterion is that the number of iterations does not exceed $15,000$.  Choose $z_0=\left(x_0,y_0\right),~x_0=y_0=\left(1,\cdots,1\right)^T$.
\begin{figure}[htbp]
  \centering
  \begin{subfigure}[b]{0.3\textwidth}
    \centering
    \includegraphics[width=\textwidth]{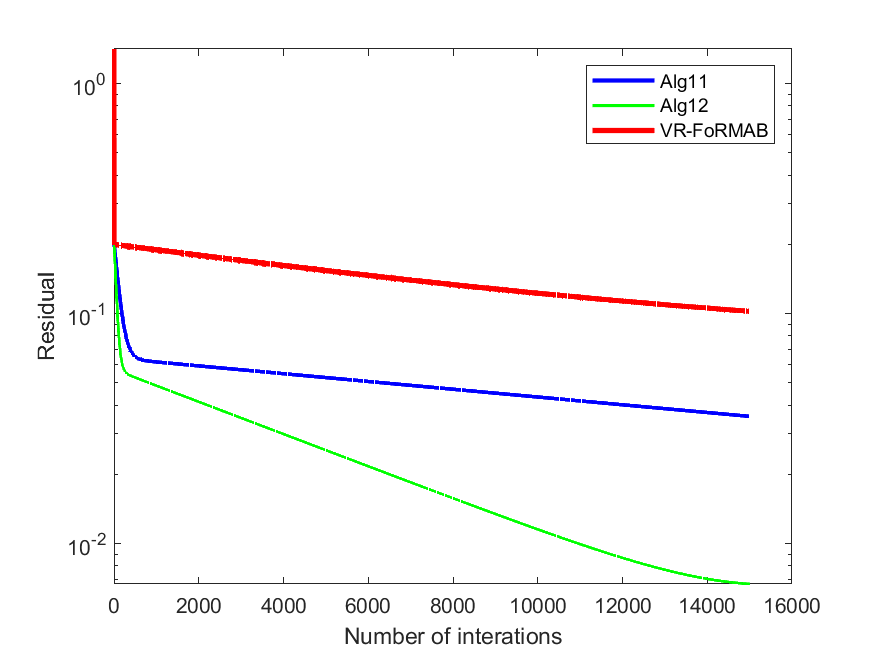}
    \label{n1fig:sub1}
    \caption{n=100}
  \end{subfigure}
  \begin{subfigure}[b]{0.3\textwidth}
    \centering
    \includegraphics[width=\textwidth]{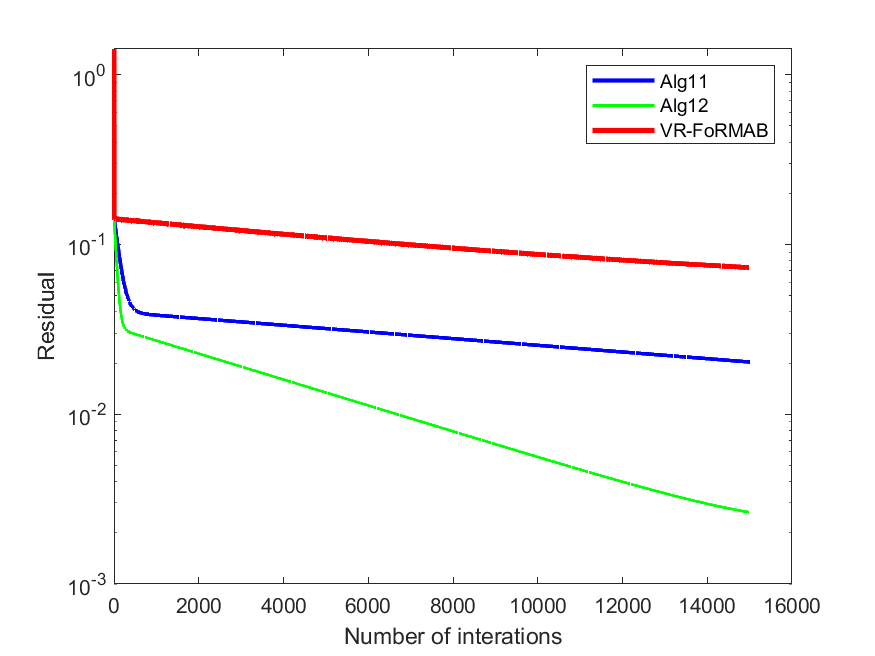}
    \label{n1fig:sub2}
    \caption{n=200}
  \end{subfigure}
  \begin{subfigure}[b]{0.3\textwidth}
    \centering
    \includegraphics[width=\textwidth]{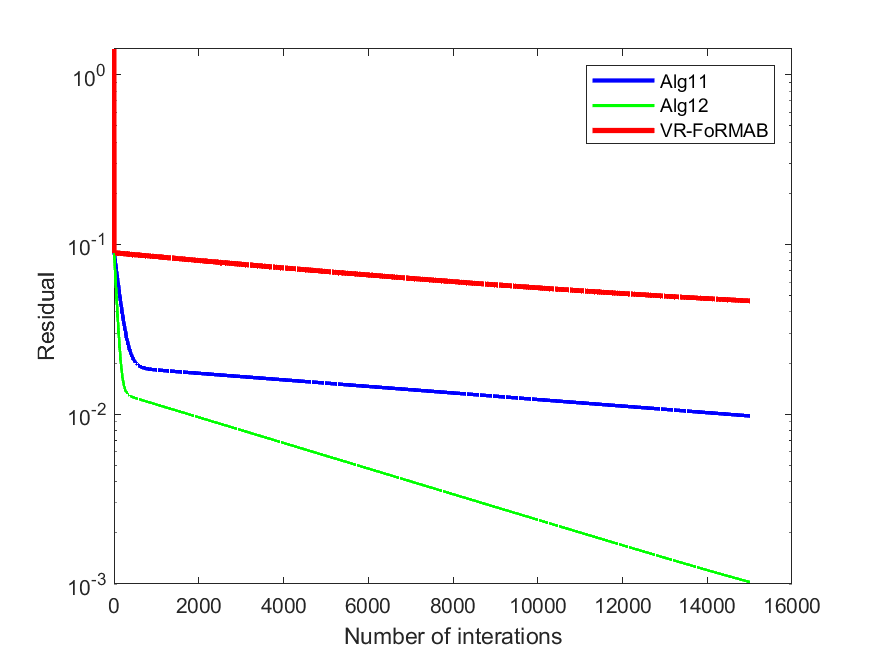}
    \label{n1fig:sub3}
    \caption{n=500}
  \end{subfigure}
   \caption{Numerical Behavior of the Residual for Example \ref{e4.2} under Different Dimensional Settings}
  \label{fig2}
\end{figure}
\par The results are presented in Figure \ref{fig2}. As illustrated, Algorithm \ref{Alg1} achieves superior performance over VR-FoRMAB across all tested dimensional settings.
\end{example}
\section{Conclusions}
\par
This work introduces a novel single-loop variance reduction algorithm, incorporating Bregman distance, to solve the FSVI \eqref{1.1}. Under mild and appropriate assumptions, we rigorously analyze the convergence of the proposed method, establishing its convergence rate and oracle complexity for Algorithm \ref{Alg1}. The principal advantages of our approach are summarized as follows:
\begin{itemize}
\item[(i)]
Our proposed method needs only a single computation of the Bregman proximal mapping over the feasible set per iteration, while inherently integrating both inertial acceleration and a batch size selection strategy. In contrast to VR-FoRMAB of \cite{AZ}, our algorithm avoids the use of the SPIDER technique proposed \cite{Fang}.
\item[(ii)]
 In the monotone setting, our algorithm achieves both almost sure convergence and an oracle complexity of $\mathcal{O}\left(\frac{\sqrt{M}}{\varepsilon }\right)$. For non-monotone cases with weak Minty solutions, the complexity is further enhanced to $\mathcal{O}\left({1}/{\varepsilon^2}\right)$. These results either match or surpass the current state-of-the-art bounds reported in existing literature.
 Notably, to the best of our knowledge, Theorems \ref{T3.2} and \ref{T3.2'} provide the first linear convergence results for a single-loop method within the Bregman framework for solving FSVI \eqref{1.1}, thereby addressing the gap in the existing research.
\item[(iii)]
 Through empirical evaluation, Algorithm \ref{Alg1} is shown to be more efficient than the existing algorithms reported in \cite{AA,AA2,AZ,PA}.
 \end{itemize}
 \begin{acknowledgements}
This study received financial support from Railway Foundation Joint Research Funds of the National
Natural Science Foundation of China (U2368216), the National Natural Science Foundation of China (11701479, 11701478), the Chinese Postdoctoral Science Foundation (2018M643434), the Fundamental Research Funds for the Central Universities (2682021ZTPY040) and the Science, Research and Innovation Promotion Funding (TSRI) (Grant no. FRB660012/0168).
\end{acknowledgements}

\end{document}